\documentclass[12pt,fullpage]{article}

\usepackage{amsfonts,xr,graphicx,amsmath,colortbl,epsfig,subfigure}

\def\to{\rightarrow}

\def\R{\mathbb{R}}
\def\N{\mathbb{N}}
\def\tilde{\widetilde}
\def\epsilon{\varepsilon}

\newcommand{\SE}{\setcounter{equation}{0} \section}
\newcommand{\be}{\begin{equation}}
\newcommand{\ee}{\end{equation}}
\newcommand{\baa}{\begin{array}}
\newcommand{\eaa}{\end{array}}
\newcommand{\ba}{\begin{eqnarray}}
\newcommand{\ea}{\end{eqnarray}}

\newtheorem{theo}{\bf Theorem}[section]
\newtheorem{lem}[theo]{\bf Lemma}
\newtheorem{pro}[theo]{\bf Proposition}
\newtheorem{cor}[theo]{\bf Corollary}
\newtheorem{defi}[theo]{\bf Definition}
\newtheorem{rem}[theo]{\bf Remark}

\begin{document}
\date{}
\title{\bf{Bistable transition fronts in $\R^N$}}
\author{Fran\c cois Hamel$\,$\thanks{The research leading to these results has received funding from the French ANR within the project PRE\-FERED and from the European Research Council under the European Union's Seventh Framework Programme (FP/2007-2013) / ERC Grant Agreement n.321186 - ReaDi -Reaction-Diffusion Equations, Propagation and Modelling. Part of this work was also carried out during visits by the author to the Departments of Mathema\-tics of the University of California, Berkeley and of Stanford University, the hospitality of which is thankfully acknowledged.}\\
\\
\footnotesize{Aix-Marseille Universit\'e \& Institut Universitaire de France}\\
\footnotesize{LATP (UMR CNRS 7353), 39 rue F. Joliot-Curie, F-13453 Marseille Cedex 13, France}}
\maketitle

\begin{abstract}
\noindent{}This paper is chiefly concerned with qualitative properties of some reaction-diffusion fronts. The recently defined notions of transition fronts generalize the standard notions of traveling fronts. In this paper, we show the existence and the uniqueness of the global mean speed of bistable transition fronts in $\R^N$. This speed is proved to be independent of the shape of the level sets of the fronts. The planar fronts are also characterized in the more general class of almost-planar fronts with any number of transition layers. These qualitative properties show the robustness of the notions of transition fronts. But we also prove the existence of new types of transition fronts in $\R^N$ that are not standard traveling fronts, thus showing that the notions of transition fronts are broad enough to include other relevant propagating solutions.
\end{abstract}

\tableofcontents


\SE{Introduction}\label{sec1}

This paper is concerned with some existence results and qualitative properties of generalized fronts, and with some estimates of their propagation speeds, for semilinear parabolic equations of the type
\be\label{eq}
u_t=\Delta u+f(u),\ \ (t,x)\in\R\times\R^N,
\ee
where $u_t=\frac{\partial u}{\partial t}$ and $\Delta$ denotes the Laplace operator with respect to the space variables $x\in\R^N$.

Throughout the paper, the reaction term $f:[0,1]\to\R$ is a $C^1$ function such that
\be\label{f}
f(0)=f(1)=0,\ \ f'(0)<0\ \hbox{ and }\ f'(1)<0.
\ee
Both zeroes $0$ and $1$ of $f$ are then stable. We define
\be\label{thetapm}
\theta^-=\min\big\{s\in(0,1);\ f(s)=0\big\}\hbox{ and }\theta^+=\max\big\{s\in(0,1);\ f(s)=0\big\}.
\ee
There holds $0<\theta^-\le\theta^+<1$. A particular important case corresponds to bistable nonlinearities $f$ which, in addition to~(\ref{f}), satisfy $\theta^-=\theta^+$, that is
\be\label{bistable}
\exists\,\theta\in(0,1),\ \ f<0\hbox{ on }(0,\theta)\hbox{ and }f>0\hbox{ on }(\theta,1).
\ee
A typical example of a function $f$ satisfying~(\ref{bistable}) is the cubic nonlinea\-rity $f_{\theta}(s)=s(1-s)(s-\theta)$ with $0<\theta<1$, and~(\ref{eq}) is then often referred to as Nagumo's or Huxley's equation. When~$\theta=1/2$, then equation~(\ref{eq}) with the function~$f(s)=8f_{1/2}((s+1)/2)=s-s^3$ with~$s\in[-1,1]$ corresponds to the celebrated Allen-Cahn equation arising in material sciences. More generally speaking, equation~(\ref{eq}) is also one of the most common reaction-diffusion equations arising in various mathematical models in biology or ecology.

In one of the main results of this paper (Theorem~\ref{thex} below), assumption~(\ref{bistable}) will be made together with~(\ref{f}). In all other results, the function $f$ is assumed to fulfil~(\ref{f}) only and can then be more general than the specific bistable type~(\ref{bistable}). 

The solution $u:\R\times\R^N\to[0,1]$ typically stands for a normalized density and it is understood as a classical solution of~(\ref{eq}). From the strong maximum principle, either $u$ is identically equal to~$0$, or it is identically equal to $1$, or it ranges in the open interval $(0,1)$. We only consider this last situation in the paper.

One of the most important aspects of these equations, which accounted for their success, is the description of propagation phenomena. By that we understand front-like time-global (also called entire, or eternal) solutions~$u$ of~(\ref{eq}) which connect the two stable stationary states $0$ and $1$ and which, in general, move as time runs (see precise definitions later). These solutions are an important class of solutions in that they usually describe the large-time behavior of the solutions of the associated Cauchy problems (some important large-time dynamics and stability results will actually be used in the proofs of the present paper). Much work has been devoted in the last decades to the study of standard front-like solutions for equations of the type~(\ref{eq}), and some of the main results in this very active field will be recalled below. On the other hand, new more general notions of propagation speeds and transition fronts, including the standard traveling fronts, have been introduced recently.

In this paper, we prove the existence and the uniqueness of the speed of any front among the class of transition fronts connecting $0$ and $1$ for problem~(\ref{eq}). We also establish some one-dimensional symmetry properties and various classification results related to the shape of the level sets of the fronts. All these qualitative properties show the robustness of the notions of transition fronts.

Furthermore, we prove the existence of new types of transition fronts in that are not known standard traveling fronts, thus showing that the notions of transition fronts are broad enough to include other relevant propagating solutions.

Before doing so, we first review the main existence and qualitative results for the standard traveling fronts.


\subsection{Standard traveling fronts}\label{sec11}

\subsubsection*{One-dimensional traveling fronts}

On the one-dimensional real line, standard traveling fronts are solutions of the type
$$u(t,x)=\phi_f(x-c_ft),$$
where $c_f\in\R$ is the propagation speed and $\phi_f:\R\to[0,1]$ is the propagation profile, such that
\be\label{eqphi}\left\{\baa{l}
\phi_f''+c_f\phi_f'+f(\phi_f)=0\hbox{ in }\R,\vspace{3pt}\\
\phi_f(-\infty)=1\hbox{ and }\phi_f(+\infty)=0.\eaa\right.
\ee
The profile $\phi_f$ is then a heteroclinic connection between the stable states $0$ and $1$. Such solutions $u(t,x)=\phi_f(x-c_ft)$ move with constant speed $c_f$ and they are invariant in the moving frame with speed~$c_f$.

If~$f$ satisfies~(\ref{bistable}) in addition to~(\ref{f}), then such fronts exist, see~\cite{aw,fm,k}. Under the sole condition~(\ref{f}), such fronts do not exist in general but further more precise conditions for the existence and non-existence have been given by Fife and McLeod in~\cite{fm}. For instance, if~$f$ satisfies~(\ref{f}) and if there are some real numbers $0<\theta^-<\mu<\theta^+<1$ such that
$$\left\{\baa{l}
f(\theta^-)=f(\mu)=f(\theta^+)=0,\ \ f'(\mu)<0,\vspace{3pt}\\
f<0\hbox{ on }(0,\theta^-)\cup(\mu,\theta^+)\ \hbox{ and }\ f>0\hbox{ on }(\theta^-,\mu)\cup(\theta^+,1),\eaa\right.$$
then there some fronts $(c^-_f,\phi^-_f)$ and $(c^+_f,\phi^+_f)$ connecting $0$ and $\mu$, and $\mu$ and $1$ respectively, since the function $f$ is of the bistable type on each of the subintervals $[0,\mu]$ and $[\mu,1]$; furthermore, in this case, it has been proved in~\cite{fm} that a front $(c_f,\phi_f)$ solving~(\ref{eqphi}) exists if and only if~$c^-<c^+$, and then the inequalities $c^-<c_f<c^+$ hold necessarily.

Coming back to the general case~(\ref{f}), if a front~$(c_f,\phi_f)$ solving~(\ref{eqphi}) exists, then the speed~$c_f$ is unique, it only depends on $f$ and it has the sign of~$\int_0^1f$, see~\cite{f,fm}. In particular, if $f$ is balanced, that is $\int_0^1f=0$, then~$c_f=0$. Moreover, the profile~$\phi_f$, if it exists, is unique up to shifts, it is such that $\phi'_f<0$ in~$\R$ and it can be assumed to be fixed with the normalization~$\phi_f(0)=1/2$. Lastly, when they exist, these fronts are globally stable in the sense that any solution of the Cauchy problem~$u_t=u_{xx}+f(u)$ for $t>0$ with an initial condition~$u(0,\cdot):\R\to[0,1]$ such that
$$\liminf_{x\to-\infty}u(0,x)>\theta^+\ge\theta^->\limsup_{x\to+\infty}u(0,x)$$
converges to the traveling front $\phi_f(x-c_ft+\xi)$ uniformly in $x\in\R$ as $t\to+\infty$, where $\xi$ is a real number which only depends on $u(0,\cdot)$ and $f$, see~\cite{fm,gr}. Let us mention here that the uniqueness of the speed $c_f$ is in sharp contrast with the case of positive nonlinearities $f>0$ on~$(0,1)$, for which the set of admissible speeds is a continuum~$[c^*_f,+\infty)$ with $c^*_f>0$, see e.g.~\cite{aw}.

{\it Throughout the paper $($except in the specific Proposition~$\ref{proplanar2}$ in Section~$\ref{sec32}$ below$)$, we assume~$(\ref{f})$ and the existence $($and then the uniqueness$)$ of a planar front~$(c_f,\phi_f)$ solving~$(\ref{eqphi})$.} In particular, we insist on the fact that all results of this paper hold if~$f$ is of the important bistable type~(\ref{bistable}). In one of the main results (Theorem~$\ref{thex}$), we actually assume additionally that $f$ has the bistable profile~$(\ref{bistable})$.

\subsubsection*{Standard planar traveling fronts in $\R^N$ with $N\ge 1$}

In any dimension $N\ge 1$, planar fronts
$$u(t,x)=\phi_f(x\cdot e-c_ft),$$
if any, are unique up to shifts, for any given unit vector~$e$ of $\R^N$, where the one-dimensional profile $\phi_f$ is as above. The level sets of such traveling fronts are parallel hyperplanes which are orthogonal to the direction of propagation~$e$. These fronts are invariant in the moving frame with speed~$c_f$ in the direction~$e$ and the unique speed $c_f$, if any, can then be referred in the sequel as the speed of planar fronts connecting $0$ and $1$ for problem~(\ref{eq}). Lastly, when~$f$ is of the bistable type~(\ref{bistable}), these planar fronts, which exist, are known to be stable with respect to some natural classes of perturbations, see~\cite{lx,mn,mnt,x2}.

\subsubsection*{Standard non-planar traveling fronts in $\R^N$ with $N\ge 2$}

When $N\ge 2$ and $f$ fulfills~(\ref{bistable}) with, say, $c_f>0$, there are other traveling fronts, which have non-planar level sets. That is, there are fronts whose profiles are still invariant in a moving frame with constant speed, but whose level sets are not hyperplanes anymore. Namely, taking~$x_N$ as the direction of propagation without loss of generality, calling $x'=(x_1,\ldots,x_{N-1})$ and~$|x'|=(x_1^2+\cdots+x_{N-1}^2)^{1/2}$ and letting $\alpha\in(0,\pi/2)$ be any given angle, equation~(\ref{eq}) admits ``conical-shaped" axisymmetric non-planar fronts of the type 
$$u(t,x)=\phi(|x'|,x_N-ct)$$
such that
\begin{figure}\label{figcurved}\centering
\includegraphics[scale=0.4]{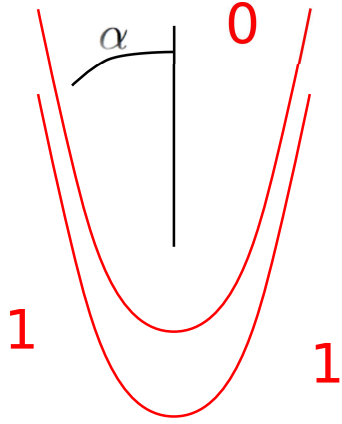}
\caption{Level sets of a conical-shaped curved front}
\end{figure}
\be\label{phipsi}\left\{\baa{l}
\phi(r,z)\to 1\hbox{ (resp. }0\hbox{) unif. as }z-\psi(r)\to-\infty\hbox{ (resp. }+\infty\hbox{)},\vspace{3pt}\\
c=\displaystyle{\frac{c_f}{\sin\alpha}}\hbox{ and }\psi'(+\infty)=\cot\alpha,\eaa\right.
\ee
for some $C^1$ function $\psi:[0,+\infty)\to\R$, see \cite{hmr2,nt1} and the joint figure. When $c_f<0$, the same result holds after changing the roles of the limits $0$ and $1$. Fronts $u(t,x)=\phi(x',x_N-ct)$ with non-axisymmetric shape, such as pyramidal fronts, are also known to exist, see~\cite{t1,t5}. A large literature has been devoted to the study of these axisymmetric and non-axisymmetric fronts in the recent years. For symmetry, uniqueness, stability and further qualitative properties of these traveling fronts, we refer to~\cite{hm,hmr2,hmr3,nt1,nt2,rrm,t2,t5} (see also~\cite{hasc} for the existence of conical-shaped fronts for some systems of reaction-diffusion equations and for angles $\alpha$ close to $\pi/2$).

When $N\ge 2$ and $f$ fulfills~(\ref{bistable}) with $c_f=0$, fronts $u(t,x)=\phi(|x'|,x_N-ct)$ with conical-shaped level sets cannot exist anymore, see~\cite{hm}. Nevertheless, for every $c\neq 0$, there exist some fronts $u(t,x)=\phi(|x'|,x_N-ct)$ such that $\phi(r,z)\to 1$ (resp. $0$) as $z\to-\infty$ (resp. $+\infty$) for every $r\ge 0$ and whose level sets have an exponential shape (if $N=2$) or a parabolic shape (if $N\ge 3$), see~\cite{cghnr}. The axisymmetry, up to shifts, of these fronts in dimension $N=2$ has been proved in~\cite{gu}. Furthermore, when $N\ge 3$ and $|c|\neq0$ is small enough,~(\ref{eq}) also admits fronts~$u(t,x)=\phi(|x'|,x_N-ct)$ such that $\phi(r,\pm\infty)=0$ for every~$r\ge 0$, $\sup_{x\in\R^N}\phi(|x'|,x_N)=1$ and the level set
$$E=\Big\{x\in\R^N;\ \phi(|x'|,x_N)=\frac{1}{2}\Big\}$$
is made of two non-Lipschitz graphs, see~\cite{dkw2} (other axisymmetric fronts exist for which~$E$ has only one connected component and has the shape of a catenoid, see~\cite{dkw2}). On the other hand, for any $N\ge 2$ (still with $c_f=0$), planar stationary fronts $u(t,x)=\phi_f(x\cdot e)$ obviously still exist, for any unit vector $e$. Stationary solutions $u(x)$ of~(\ref{eq}) such that $u_{x_N}<0$ and~$u(x',x_N)\to 1$ (resp.~$0$) as $x_N\to-\infty$ (resp. $+\infty$) are necessarily planar --hence of the type~$\phi_f(x\cdot e)$ up to shifts-- if $N\le8$~\cite{ac,gg,s} whereas they are not always planar when~$N\ge 9$~\cite{dkw1} (this problem is related to a celebrated conjecture by De Giorgi~\cite{dg} for~$f(s)=s(1-s)(s-1/2)$). Other non-monotone stationary saddle-shaped solutions or solutions whose level sets have multiple ends are also known to exist for some balanced bistable functions~$f$, see~\cite{acm,c,ct1,ct2,dfp,dkpw}.


\subsection{Notions of transition fronts and global mean speed}\label{sec12}

The above examples show that equation~(\ref{eq}) admits many types of traveling fronts. For all of them, the solutions~$u$ converge to the stable states~$0$ or~$1$ far away from their moving or stationary level sets, uniformly in time. This observation will be the key-point of the more general notion of transition fronts given in Definition~\ref{def1} below. Furthermore, another common property fulfilled by all the standard traveling fronts is that their level sets move at the global mean speed~$|c_f|$, only depending on~$f$, in a sense to be made more precise below. One of the main goals of the present paper is actually to prove that this property is shared by all transition fronts.

Let us now describe the general notions of transition fronts and global mean speed for problem~(\ref{eq}). First, for any two subsets $A$ and $B$ of $\R^N$ and for $x\in\R^N$, we set
\be\label{dAB}
d(A,B)=\inf\big\{|x-y|,\ (x,y)\in A\times B\big\}
\ee
and $d(x,A)=d(\{x\},A)$. The notions of transition fronts and global mean speeds are borrowed from~\cite{bh5} and are adapted here to the case of connections between the constant stationary states~$0$ and $1$ in~$\R^N$. They involve two families $(\Omega^-_t)_{t\in{\mathbb R}}$ and $(\Omega^+_t)_{t\in{\mathbb R}}$ of open nonempty subsets of~$\R^N$ such that
\begin{equation}\label{omegapm}
\forall\,t\in\R,\quad\left\{\baa{l}
\Omega^-_t\cap\Omega^+_t=\emptyset,\vspace{3pt}\\
\partial\Omega^-_t=\partial\Omega^+_t=:\Gamma_t,\vspace{3pt}\\
\Omega^-_t\cup\Gamma_t\cup\Omega^+_t=\R^N,\vspace{3pt}\\
\sup\big\{d(x,\Gamma_t);\ x\in\Omega^+_t\big\}=\sup\big\{d(x,\Gamma_t);\ x\in\Omega^-_t\big\}=+\infty\eaa\right.
\end{equation}
and
\be\label{unifgamma}\left\{\baa{l}
\inf\Big\{\sup\big\{d(y,\Gamma_t);\ y\in\Omega^+_t,\ |y-x|\le r\big\};\ t\in\R,\ x\in\Gamma_t\Big\}\to+\infty\vspace{3pt}\\
\inf\Big\{\sup\big\{d(y,\Gamma_t);\ y\in\Omega^-_t,\ |y-x|\le r\big\};\ t\in\R,\ x\in\Gamma_t\Big\}\to+\infty\eaa\right.\hbox{ as }r\to+\infty.\footnote{In ~\cite{bh5}, the condition $|y-x|=r$ was used instead of $|y-x|\le r$ (more precisely, the condition $d_{\Omega}(y,x)=r$ was used, where $d_{\Omega}$ denotes the geodesic distance in a domain $\Omega\subset\R^N$). In the case of the whole space $\R^N$, the condition~$|y-x|\le r$ in~(\ref{unifgamma}) of the present paper is broader than just $|y-x|=r$ and in some sense more natural. However, it is straightforward to check that all qualitative properties stated in~\cite{bh5} still hold with the condition~$|y-x|\le r$ or with $d_{\Omega}(y,x)\le r$ in a general domain $\Omega\subset\R^N$, instead of the corresponding condition~(1.5) of~\cite{bh5}.}
\ee
Notice that the condition~(\ref{omegapm}) implies in particular that the interface~$\Gamma_t$ is not empty for every $t\in\R$. As far as condition~(\ref{unifgamma}) is concerned, it is illustrated in the joint figure. Roughly speaking, it means that, when $r$ is large, for every $t\in\R$ and every point~$x\in\Gamma_t$, there are some points $y^+_{r,t,x}$ and $y^-_{r,t,x}$ in both $\Omega^+_t$ and $\Omega^-_t$ which are far from $\Gamma_t$ and are at a distance at most~$r$ from~$x$ (in particular, the points $y^{\pm}_{r,t,x}$ belong to the geodesic tubular neighborhoods of width~$r$ of the sets~$\Gamma_t$). Moreover, the sets $\Gamma_t$ are assumed to be made of a finite number of graphs: there is an integer $n\ge 1$ such that, for each $t\in\R$, there are $n$ open subsets~$\omega_{i,t}\subset\R^{N-1}$ (for~$1\le i\le n$), $n$ continuous maps $\psi_{i,t}:\omega_{i,t}\to\R$ and $n$ rotations $R_{i,t}$ of~$\R^N$, such that
\be\label{omegapmbis}
\Gamma_t\subset\mathop{\bigcup}_{1\le i\le n}R_{i,t}\Big(\big\{x\in\R^N;\ x'\in\omega_{i,t},\ x_N=\psi_{i,t}(x')\big\}\Big).
\ee
When $N=1$,~(\ref{omegapmbis}) means that $\Gamma_t$ has cardinal at most $n$, that is $\Gamma_t=\{x_{1,t},\cdots,x_{n,t}\}$, for each~$t\in\R$.
\begin{figure}\label{figinterfaces}\centering
\includegraphics[scale=0.5]{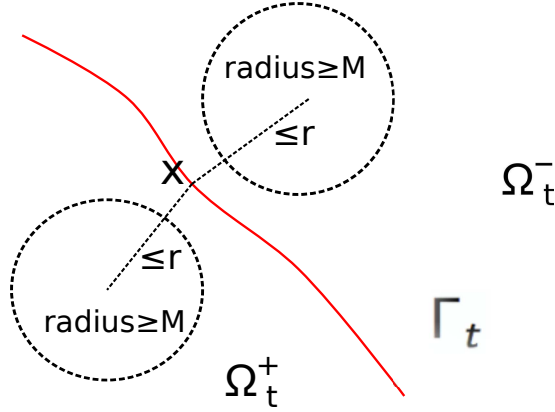}
\caption{Geometrical interpretation of the condition~(\ref{unifgamma})}
\end{figure}

\begin{defi}\label{def1} {\rm{\cite{bh5}}} For problem~$(\ref{eq})$, a transition front connecting $0$ and $1$ is a classical solution $u:\R\times\R^N\to(0,1)$ for which there exist some sets $(\Omega^{\pm}_t)_{t\in\R}$ and $(\Gamma_t)_{t\in\R}$ satisfying~$(\ref{omegapm})$,~$(\ref{unifgamma})$ and~$(\ref{omegapmbis})$, and, for every $\epsilon>0$, there exists $M\ge 0$ such that
\be\label{defunif}\left\{\baa{l}
\forall\, t\in\R,\ \forall\,x\in\Omega^+_t,\quad\big(d(x,\Gamma_t)\ge M\big)\Longrightarrow\big(u(t,x)\ge1-\epsilon\big),\vspace{3pt}\\
\forall\, t\in\R,\ \forall\,x\in\Omega^-_t,\quad\big(d(x,\Gamma_t)\ge M\big)\Longrightarrow\big(u(t,x)\le\epsilon\big).\eaa\right.
\ee
Furthermore, $u$ is said to have a global mean speed $\gamma\ (\ge 0)$ if
\be\label{defmeanspeed}
\frac{d(\Gamma_t,\Gamma_s)}{|t-s|}\to\gamma\ \hbox{ as }|t-s|\to+\infty.
\ee
\end{defi}

Before stating the main results of this paper, let us comment the notions given in Definition~\ref{def1} and let us connect them with the standard notions of traveling fronts. Firstly, it is easy to check that all moving or stationary fronts mentioned in Section~\ref{sec11} are transition fronts connecting~$0$ and~$1$ in the sense of Definition~\ref{def1}, for some suitable choices of sets $(\Omega^{\pm}_t)_{t\in\R}$. For instance, when~$f$ fulfils~(\ref{bistable}) with $c_f\neq 0$, the conical-shaped fronts $u(t,x)=\phi(|x'|,x_N-ct)$ satisfying~(\ref{phipsi}) are transition fronts connecting~$0$ and~$1$ with, say,
$$\Omega^+_t=\big\{x\in\R^N;\ x_N<\psi(|x'|)+ct\big\},\ \ \Omega^-_t=\big\{x\in\R^N;\ x_N>\psi(|x'|)+ct\big\}$$
and $\Gamma_t=\{x\in\R^N;\ x_N=\psi(|x'|)+ct\}$ for every $t\in\R$. However, Definition~\ref{def1} covers other transition fronts than the one described in Section~\ref{sec11} (see Theorem~\ref{thex} below). Definition~\ref{def1} was actually given in~\cite{bh5} for more general domains, equations or limiting states (instead of~$0$ and $1$). In the recent years, many papers have been devoted to the existence and stability of transition fronts for equations of the type~(\ref{eq}) in~$\R$ or in infinite cylinders with $x$-dependent \cite{mnrr,mrs,nrrz,nr,z1,z2} and~$t$-~or~$(t,x)$-dependent \cite{bh5,naro,s1,s2,s3,s4,s5,s6} bistable, combustion or monostable nonlinearities~$f$, as well as for~(\ref{eq}) in exterior \cite{bhm1} or cylindrical-type domains (see~\cite{bbc,cg,m,rrbk}, where blocking phenomena are also shown). We also mention \cite{abc,dgm,fz,he2,px,ww,x1,x3,x4} for the existence and qualitative properties of bistable pulsating fronts in periodic media. In a subsequent paper~\cite{h10}, we establish some bounds and estimates for the mean speed of transition fronts in heterogeneous media or in more general domains. In the present paper, for the sake of clarity and homogeneity of the presentation, we only focus on the case of the homogeneous equation~(\ref{eq}) in the whole space $\R^N$ but we prove that even this simple-looking problem already has many deep properties: new classification results and general estimates shared by all transition fronts are shown, not to mention the existence of new transition fronts.

For a given transition front $u$, the sets $(\Omega^{\pm}_t)_{t\in\R}$ and $(\Gamma_t)_{t\in\R}$ are not uniquely determined, in the sense that two families $(\Omega^{\pm}_t)_{t\in\R}$ (resp. $(\Gamma_t)_{t\in\R}$) and $(\widetilde{\Omega}^{\pm}_t)_{t\in\R}$ (resp. $(\widetilde{\Gamma}_t)_{t\in\R}$) may be associated to the same transition font $u$. Nevertheless, due to the key uniformity property in~(\ref{defunif}), the sets $\Gamma_t$ are located at a uniformly bounded distance of any given level set of $u$ and this boundedness property is intrinsic. Namely, under the assumption that~$\sup\big\{d(x,\Gamma_{t-\tau});\ t\in\R,\ x\in\Gamma_t\big\}<+\infty$ for some $\tau>0$ (the interfaces $\Gamma_t$ and $\Gamma_{t-\tau}$ are in some sense not too far from each other), it follows from Theorem~1.2 of~\cite{bh5} that
$$\forall\,\lambda\in(0,1),\quad\sup\big\{d(x,\Gamma_t);\ u(t,x)=\lambda,\ (t,x)\in\R\times\R^N\big\}<+\infty$$
and, for every $C\ge 0$, there is $\eta>0$ such that $\eta\le u(t,x)\le 1-\eta$ for all $(t,x)\in\R\times\R^N$ with~$d(x,\Gamma_t)\le C$. Roughly speaking, the transition zone between $0$ and $1$ is thus a neighborhood with uniformly bounded width of the (possibly moving) interfaces~$\Gamma_t$ (notice however that the transition zone between $0$ and $1$ can also be made of several possibly disconnected transition zones, each of them being a neighborhood of a graph, see for instance the aforementioned examples from~\cite{dkw2} and another example after Theorem~\ref{thplanar} below). The word transition in Definition~\ref{def1} thus corresponds to the intuitive idea of a spatial transition (we refer to~\cite{n} for the related but different notion of critical transition).

The global mean speed~$\gamma$ of a transition front, if any, corresponds to the limiting average speed of the minimal distance between the interfaces~$\Gamma_t$ and can then be viewed as a mean smallest normal speed of the interfaces~$\Gamma_t$. Since the sets $\Gamma_t$ are not uniquely determined, the notion of instantaneous normal speed of $\Gamma_t$ has no sense. However, due to the remarks of the previous paragraph, the notion of global mean speed~$\gamma$ given in~(\ref{defmeanspeed}) is meaningful. As a matter of fact, it is essential to say that, for a given transition front $u$, the global mean speed $\gamma$, if any, is uniquely determined and does not depend on the specific choice of the sets~$(\Omega^{\pm}_t)_{t\in\R}$ and~$(\Gamma_t)_{t\in\R}$, see Theorem~1.7 in~\cite{bh5}. For instance, when~$f$ fulfils~(\ref{bistable}) with~$c_f\neq 0$, the conical-shaped fronts~$u(t,x)=\phi(|x'|,x_N-ct)$ satisfying~(\ref{phipsi}) have a global mean speed~$\gamma$ and $\gamma=|c_f|$, whatever the angle~$\alpha\in(0,\pi/2)$ may be. For any such front~$u$, the speed~$c=c_f/\sin\alpha$ is the speed in the vertical direction $x_N$ of the frame in which the front is invariant, but the asymptotical smallest normal speed of the level sets of $u$ is equal to~$|c_f|$. It is also straightforward to check that, when~$f$ fulfils~(\ref{bistable}) with $c_f=0$, the fronts mentioned in Section~\ref{sec11} have global mean speed $\gamma=0$: this fact is obvious when the fronts are stationary, since the $\Gamma_t$ can all be chosen as any given (time-independent) level set, but this property also holds good for the exponentially-shaped or parabolic-shaped fronts since the level sets have an infinite slope in the $(x',x_N)$ coordinates as $|x'|\to+\infty$. These exponentially-shaped or parabolic-shaped fronts have zero global mean speed $\gamma=0$, but they are not stationary.

What is much stronger and not trivial at all is to show that, whatever the shape of the fronts and the value of $c_f$ may be, all transition fronts for~(\ref{eq}) have a global mean speed and this speed is equal to $|c_f|$: this will be one of the main results of this paper, see Theorem~\ref{thcf} below.


\SE{Main results}\label{sec2}

The first main results are concerned with some qualitative geometrical properties of the transition fronts connecting $0$ and $1$ for problem~(\ref{eq}), including some new classification Liouville-type results, and with some estimates of their global mean speed. More precisely, in the following subsections, we first give a characterization of the planar fronts among the more general class of almost-planar transition fronts. We then give a characterization of the mean speed of all transition fronts. Lastly, we deal with the existence of new non-standard transition fronts.


\subsection{Almost-planar and planar fronts}\label{sec21}

Planar fronts connecting $0$ and $1$ for~(\ref{eq}) are solutions of the type $\phi(x\cdot e-ct)$ with~$\phi(-\infty)=1$ and $\phi(+\infty)=0$. As recalled in Section~\ref{sec11}, the function $\phi=\phi_f$, if any, is unique up to shifts and the speed $c=c_f$, if any, is unique. These fronts have planar level sets and they are monotone with respect to the direction of propagation, at each time $t$. In particular, they fall within the more general class of almost-planar fronts introduced in~\cite{bh5}, and defined as follows.

\begin{defi}\label{defap}
A transition front $u$ in the sense of Definition~$\ref{def1}$ is called almost-planar if, for every $t\in\R$, the set $\Gamma_t$ can be chosen as the hyperplane
$$\Gamma_t=\big\{x\in\R^N;\ x\cdot e_t=\xi_t\big\}$$
for some vector $e_t$ of the unit sphere ${\mathbb{S}}^{N-1}$ and some real number $\xi_t$.
\end{defi}

In other words, the level sets of almost-planar fronts are in some sense close to hyperplanes, even if they are not a priori assumed to be planar. In~\cite{bh4}, we gave a characterization of the almost-planar fronts for which $e_t=e$ is a given constant vector and for which there exists~$\gamma\ge 0$ such that $|\xi_t-\xi_t|-\gamma|t-s|$ is bounded uniformly with respect to $(t,s)\in\R^2$: such fronts have to be planar fronts $\phi_f(\pm x\cdot e-c_ft)$, up to shifts, and $\gamma=|c_f|$ (see Theorem~3.1 in~\cite{bh4}).

\begin{figure}\centering
\includegraphics[scale=0.5]{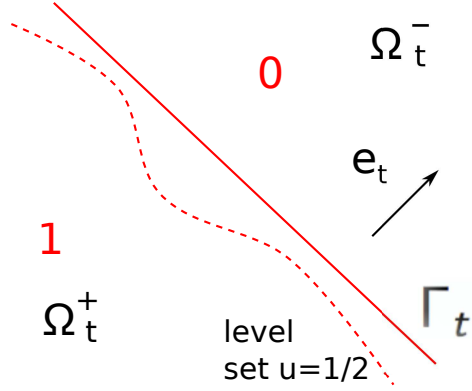}
\caption{Almost-planar fronts}
\end{figure}

In this paper, we first give a more general characterization of the planar fronts $\phi_f(x\cdot e-c_ft)$ for problem~(\ref{eq}), without assuming that the directions $e_t$ are a priori constant and without assuming any a priori bound on the positions $\xi_t$.

\begin{pro}\label{proplanar}
For problem~$(\ref{eq})$, any almost-planar transition front $u$ connecting $0$ and $1$ is planar, that is there exist a unit vector~$e$ of~$\R^N$ and a real number~$\xi$ such that
\be\label{planar}
u(t,x)=\phi_f(x\cdot e-c_ft+\xi)\hbox{ for all }(t,x)\in\R\times\R^N.
\ee
\end{pro}

The first step in the proof of Proposition~\ref{proplanar} is to show that the directions~$e_t$ appearing in Definition~\ref{defap} are equal to a constant vector~$e$ independent of time. As a consequence,~$u$ converges to $0$ and $1$ as~$x\cdot e-\xi_t\to\pm\infty$, uniformly in $t$ and in the spatial variables orthogonal to $e$. These properties have some similarities with the Gibbons conjecture about the one-dimensional symmetry of solutions $0\le v\le 1$ of elliptic equations $\Delta v+f(v)=0$ in~$\R^N$ with~$f$ satisfying~(\ref{f}) and $v(x)\to0$ and~$1$ as~$x\cdot e\to\pm\infty$ uniformly in the variables orthogonal to~$e$: for this latter problem, the solutions are proved to depend on~$x\cdot e$ only (and, necessarily,~$\int_0^1f=0$), see~\cite{bbg,bhm,fa,fv}. As for our parabolic problem~(\ref{eq}), the difference is that nothing is imposed a priori on the function $\xi_t$: to get the conclusion~(\ref{planar}), the one-dimensional stability of the planar front~$\phi_f$ (see Fife and McLeod~\cite{fm}) is used to get the boundedness of $t\mapsto\xi_t-c_ft$, together with the aforementioned parabolic Liouville type result of Berestycki and the author (Theorem~3.1 in~\cite{bh4}).

\begin{rem}{\rm As a matter of fact, the existence of the planar front~$(c_f,\phi_f)$, which is always assumed by default throughout the paper, can almost be dropped in Proposition~\ref{proplanar}. Namely, there is a $C^1([0,1])$ dense set of functions~$f$ satisfying~$(\ref{f})$ such that the existence of an almost-planar transition front connecting $0$ and $1$ for problem~$(\ref{eq})$ in $\R^N$ implies the existence of a planar one-dimensional front~$(c_f,\phi_f)$, and then the conclusion~(\ref{planar}): in other words, for these functions $f$, either there is a planar front~$(c_f,\phi_f)$ connecting $0$ and $1$, or there is no almost-planar transition front connecting $0$ and $1$ in $\R^N$. We refer to Section~\ref{sec32} and Proposition~\ref{proplanar2} below for more details. However, for the sake of the unity of the presentation, in Proposition~\ref{proplanar} as well as in all other results, we have chosen to keep the default assumption of the existence of a planar front~$(c_f,\phi_f)$. Lastly, we repeat that the assumption is fulfilled automatically if $f$ is of the bistable type~(\ref{bistable}).}
\end{rem}

It follows in particular from Proposition~\ref{proplanar} that the almost-planar fronts in any dimension~$N\ge 1$ have a (global mean) speed $\gamma=|c_f|$. Another immediate consequence of Proposition~\ref{proplanar} is a classification result in dimension $N=1$. For any non-empty set $E\subset\R^N$, let
$$\hbox{diam}(E)=\sup_{(x,y)\in E\times E}|x-y|$$
denote the Euclidean diameter of $E$.

\begin{cor} 
Let $u$ be any transition front connecting $0$ and $1$ for problem~$(\ref{eq})$ in $\R$. If~$\sup_{t\in\R}\,\hbox{diam}(\Gamma_t)<+\infty$, then $u$ is a classical traveling front $u(t,x)=\phi_f(\pm x-c_ft+\xi)$ for some $\xi\in\R$. In particular, if $u$ is almost-planar in the sense of Definition~$\ref{defap}$, then the same conclusion holds.
\end{cor}

Indeed, the boundedness of $\hbox{diam}(\Gamma_t)$ in dimension $N=1$ implies that $\Gamma_t$ can be reduced to a singleton without loss of generality, that is there is only one interface between the limi\-ting values~$0$ and~$1$. Notice that the boundedness of the width of the transition between~$0$ and $1$, which is one of the key-properties in Definition~\ref{def1}, is necessary for the conclusion of Proposition~\ref{proplanar} to hold in general, even in dimension $N=1$. For instance, if~$f$ is of the bistable type~(\ref{bistable}), there exist some solutions $u$ of~(\ref{eq}) in~$\R$ such that $0<u(t,x)<1$ in~$\R^2$,~$u(t,-\infty)=1$ and~$u(t,+\infty)=0$ for every $t\in\R$ and for which~(\ref{defunif}) is not satisfied for any families $(\Omega^{\pm}_t)_{t\in\R}$ and~$(\Gamma_t)_{t\in\R}$ satisfying~(\ref{omegapm}),~(\ref{unifgamma}) and~(\ref{omegapmbis}), see~\cite{mn1}. These solutions are indeed constructed in such a way that they are close to $\theta$ on very large intervals as $t\to-\infty$. Thus they cannot be transition fronts in the sense of Definition~\ref{def1}. 

Coming back to the transition fronts in $\R^N$ for any dimension $N\ge 1$, the conclusion of Proposition~\ref{proplanar} still holds when, at each time $t$, the transition between $0$ and $1$ is made of a finite number of bounded parallel strips, under the additional condition that the planar speed~$c_f$ is not zero. More precisely, the following result holds.

\begin{theo}\label{thplanar}
For problem~$(\ref{eq})$, let $u$ be a transition front connecting $0$ and $1$ such that, for every $t\in\R$, there are $e_t$ in $\mathbb{S}^{N-1}$ and $\xi_{1,t},\ldots,\xi_{n,t}$ in $\R$ such that
\be\label{xiit}
\Gamma_t=\bigcup_{1\le i\le n}\big\{x\in\R^N;\ x\cdot e_t=\xi_{i,t}\big\}.
\ee
If $c_f\neq 0$, then $u$ is a planar front of the type~$(\ref{planar})$.
\end{theo}

The condition that $c_f$ is not zero is actually necessary. Indeed, for some nonlinearities~$f$ such that $c_f=0$, there are transition fronts connecting $0$ and $1$ in $\R$ such that, say,
$$\Gamma_t=\big\{\xi_{1,t},\xi_{2,t},\xi_{3,t}\big\}\hbox{ for all }t<0,$$
with $\xi_{1,t}<\xi_{2,t}<\xi_{3,t}$ for every $t<0$, $\xi_{1,t}\to-\infty$, $\xi_{3,t}\to+\infty$, $|\xi_{1,t}|=o(|t|)$, $\xi_{3,t}=o(|t|)$ as~$t\to-\infty$, and
$$\Gamma_t=\big\{0\big\}\hbox{ for all }t\ge0.$$
These transition fronts, which can be derived from~\cite{er,ei}, describe the slow dynamics of some almost-stationary fronts. They have three interfaces as $t\to-\infty$, the leftmost and rightmost ones move toward the origin and disappear in finite time, and only one remains as $t\to+\infty$, in the sense that the solution converges to a finite shift of the stationary front $\phi_f(\pm x)$ as $t\to+\infty$. Notice that such solutions have global mean speed $\gamma=c_f=0$ in the sense of Definition~\ref{def1}, and that such fronts of course exist in any dimension $N\ge2$, by extending them in a trivial manner in the variables $x_2,\ldots,x_N$.


\subsection{Existence and uniqueness of the global mean speed among all transition fronts}\label{seccf}

Once we have characterized the (almost-)planar fronts for equation~(\ref{eq}), we now consider the general case of transition fronts whose level sets have arbitrary shapes. For this problem, as mentioned in Section~\ref{sec1}, the planar fronts $u(t,x)=\phi_f(x\cdot e-c_ft)$, the conical-shaped, pyramidal, exponentially-shaped or parabolic-shaped fronts $u(t,x)=\phi(x',x_N-ct)$ with $c_f\neq0$ or $c_f=0$, as well as the stationary fronts when $c_f=0$, share a common property: they all have a global mean speed and this mean speed is equal to $\gamma=|c_f|$, which depends on~$f$ only. The goal of the next theorem is to show both the existence and the uniqueness of the global mean speed of any transition front, whatever the shape of the level sets of the fronts may be and whatever the value of the planar speed $c_f$ may be.

\begin{theo}\label{thcf}
For problem~$(\ref{eq})$, any transition front connecting $0$ and $1$ has a global mean speed~$\gamma$. Furthermore, this global mean speed~$\gamma$ is equal to $|c_f|$.
\end{theo}

We point out the difference between this result and Theorem~1.7 of~\cite{bh5} recalled in Section~\ref{sec12}. Theorem~1.7 of~\cite{bh5} was concerned with the uniqueness of the global mean speed, if any, of a {\it given} transition front~$u$. Theorem~\ref{thcf} of the present paper not only shows the existence of a global mean speed for {\it any} transition front connecting $0$ and $1$, but it also shows the uniqueness of this global mean speed among {\it all} transition fronts. Notice that this existence and uniqueness result is in sharp contrast with the case of transition fronts of~$(\ref{eq})$ with other nonlinearities~$f$. For instance, if $f$ is positive and concave on $(0,1)$, then not only the admissible speeds of standard traveling fronts are not unique~\cite{aw}, but there are also some transition fronts connecting~$0$ and~$1$ which do not have any global mean speed, even in dimension $N=1$, see~\cite{hn2}.

\begin{rem}{\rm Other notions of distance could be used. For any two subsets $A$ and~$B$ of~$\R^N$, the quantity~$d(A,B)$ defined by~$(\ref{dAB})$ is the smallest geodesic distance between pairs of points in~$A$ and~$B$. Other notions are the distance $\widetilde{d}$ and the Hausdorff distance $\overline{d}$ defined by
\be\label{defdtilde}
\widetilde{d}(A,B)=\min\Big(\sup\big\{d(x,B);\,x\in A\big\},\,\sup\big\{d(y,A);\,y\in B\big\}\Big)
\ee
and
$$\overline{d}(A,B)=\max\Big(\sup\big\{d(x,B);\,x\in A\big\},\,\sup\big\{d(y,A);\,y\in B\big\}\Big).$$
There holds $d(A,B)\le\widetilde{d}(A,B)\le\overline{d}(A,B)$. It follows from the proof of Theorem~$\ref{thcf}$ (see Remark~$\ref{remhausdorff}$ below for the details) that, under the same assumptions, any transition front connecting $0$ and $1$ for equation~$(\ref{eq})$ has a global mean speed for the distance $\widetilde{d}$ and this global mean speed is equal to $|c_f|$, in the sense that
\be\label{defmeanspeedbis}
\frac{\widetilde{d}(\Gamma_t,\Gamma_s)}{|t-s|}\to|c_f|\ \hbox{ as }|t-s|\to+\infty.
\ee
For instance, for all the usual traveling fronts~$u(t,x)=\phi(x',x_N-ct)$ mentioned in Section~$\ref{sec1}$ with conical-shaped, pyramidal, exponential or parabolic level sets, the global mean speed defined by~$(\ref{defmeanspeedbis})$ exists and is equal to $|c_f|$, as for~$(\ref{defmeanspeed})$. On the other hand, all these fronts are invariant in the moving frame with speed $c$ in the direction $x_N$. The vertical speed of this specific frame can also be viewed as the asymptotic speed of the tip of the fronts, in the sense that
$$\frac{\overline{d}(\Gamma_t,\Gamma_s)}{|t-s|}\to|c|\ \hbox{ as }|t-s|\to+\infty.$$
Therefore, the Hausdorff distance gives rise to different global mean speeds, which provide another type of information about the evolution of the level sets but depend on the given transition front $($remember that $|c|$ can take all values in the interval $[|c_f|,+\infty)$ under assumption~$(\ref{bistable})$, whatever the value of~$c_f$ may be$)$. We think that the most natural notion of distance is the one defined in~$(\ref{defmeanspeed})$: it corresponds to the asymptotic smallest normal speed of the level sets. Moreover, as~$\widetilde{d}$ in~$(\ref{defmeanspeedbis})$, it has the advantage of depending only on~$f$ and thus being independent of the transition front, as shown in Theorem~$\ref{thcf}$.}
\end{rem}


\subsection{Existence of non-standard transition fronts}

The previous qualitative properties showed the strength of Definition~\ref{def1}, since the solutions of~(\ref{eq}) in the large class of transition fronts are proved to share some common features (existence and uniqueness of the global mean speed) as well as some further strong qualitative symmetry properties under some additional geometrical conditions. As far as the standard traveling fronts are concerned, all these well-known fronts $u(t,x)=\phi(x',x_N-ct)$, which were mentioned in Section~\ref{sec11} and which exist under the bistable condition~(\ref{bistable}), share another simple property, in addition to the existence and uniqueness of the global mean speed~$|c_f|$. Namely, as already emphasized, they are invariant in the frame moving with the speed $c$ in the direction $x_N$. However, what is even more intriguing is that there exist other transition fronts, which are not usual traveling fronts in the sense that there is no frame in which they are invariant as time runs. 

\begin{theo}\label{thex}
Let $N\ge 2$ and assume that $f$ is of the bistable type~$(\ref{bistable})$ with $c_f>0$. Then problem~$(\ref{eq})$ admits transition fronts $u$ connecting $0$ and $1$ which satisfy the following pro\-perty: there is no function~$\Phi:\R^N\to(0,1)$ $($independent of $t)$ for which there would be some families~$(R_t)_{t\in\R}$ and~$(x_t)_{t\in\R}$ of rotations and points in $\R^N$ such that $u(t,x)=\Phi\big(R_t(x-x_t)\big)$ for all~$(t,x)\in\R\times\R^N$.
\end{theo}

Theorem~\ref{thex} means that the class of transition fronts includes new types of solutions, even in the homogeneous space $\R^N$, thus showing the broadness of Definition~\ref{def1}. In dimension~$N=2$, the new transition fronts $u$ described in Theorem~\ref{thex} are constructed by mixing three planar fronts moving in three different directions: say, the direction $x_2$ and two directions which are symmetric with respect to the vertical axis $x_2$ (see the joint figure).  As $t\to-\infty$, the level sets of these solutions look like two symmetric oblique half-lines moving in the $x_2$ direction and separated by a larger and larger segment parallel to the $x_1$ axis. Then, as time increases, the medium segment disappears and, finally, the constructed solutions converge as~$t\to+\infty$ to a conical-shaped usual traveling front $\phi(x_1,x_2-ct)$. This scheme leads to the desired conclusion in dimension $N=2$ and then immediately in all dimensions $N\ge 3$ by trivially extending the two-dimensional solutions in the variables $x_3,\ldots,x_N$.
\begin{figure}\centering
\subfigure{\includegraphics[scale=0.35]{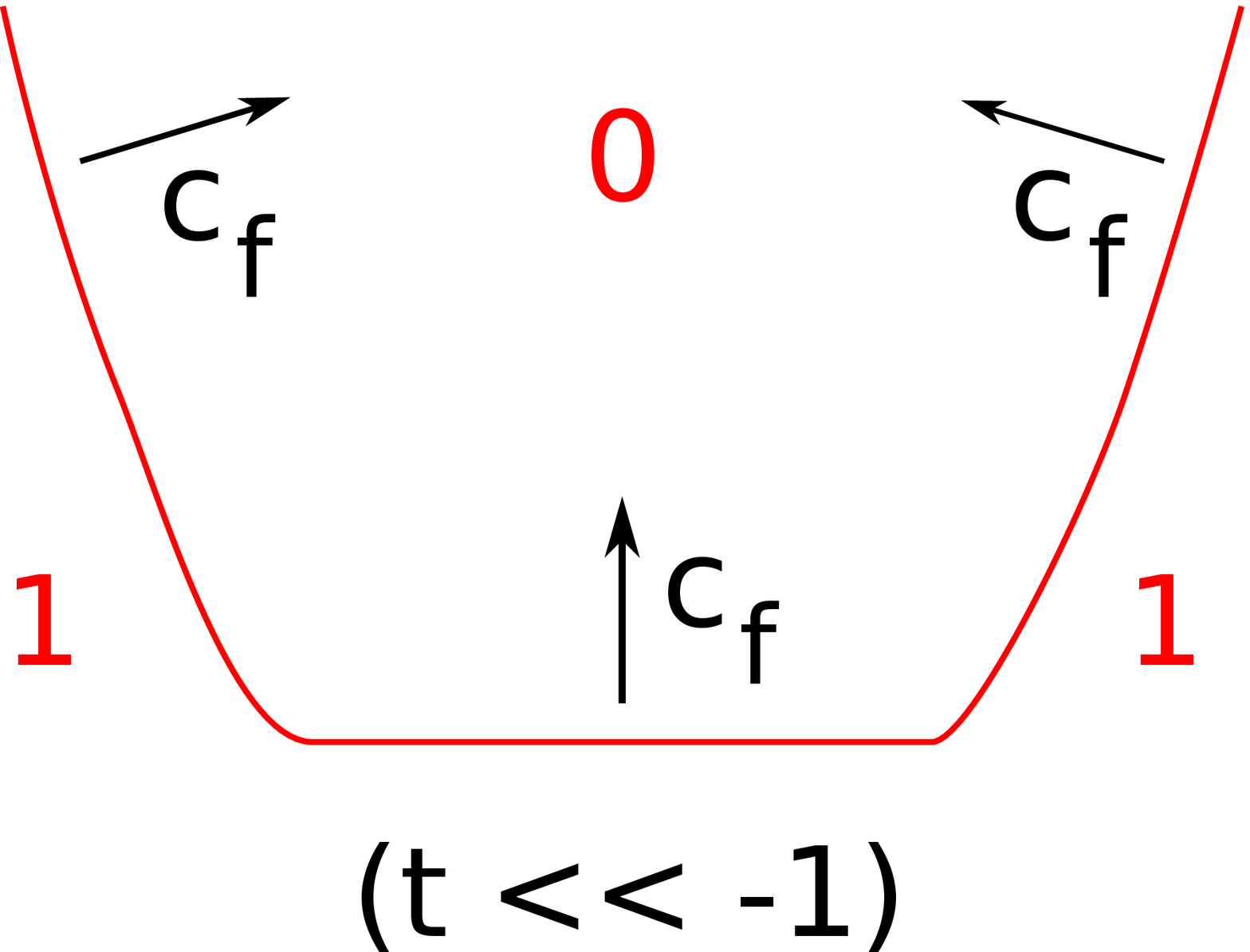}}
\hskip 2.5cm
\subfigure{\includegraphics[scale=0.35]{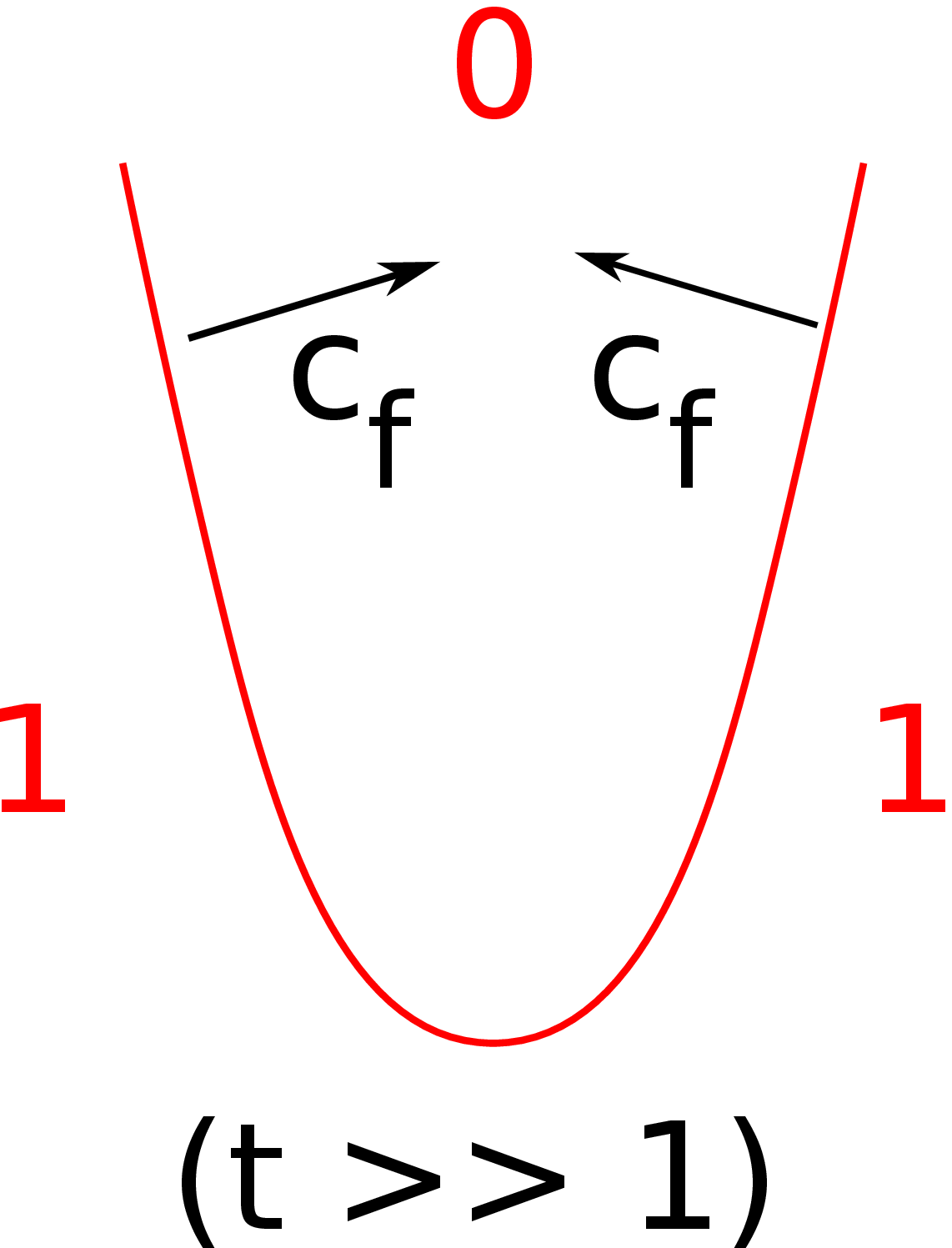}}
\caption{Example of a transition front given in Theorem~\ref{thex}}
\end{figure}
\hfill\break

\noindent{\bf{Outline of the paper.}} In the remaining sections, we perform the proof of the results. Section~\ref{sec3} is devoted to the Liouville-type results related to the characterization of planar fronts among the larger class of almost-planar fronts or fronts having finitely many parallel transition layers, that is we do the proof of Proposition~\ref{proplanar} and Theorem~\ref{thplanar}. Section~\ref{sec4} is devoted to the proof of Theorem~\ref{thcf}, that is we show the existence, the uniqueness and the characterization of the global mean speed among all transition fronts. Lastly, Section~\ref{sec5} is devoted to the proof of Theorem~\ref{thex}, that is the construction of transition fronts which are not standard traveling fronts.


\SE{Characterization of planar fronts}\label{sec3}

This section is devoted to the proof of Proposition~\ref{proplanar} and Theorem~\ref{thplanar}. That is, we characterize the planar fronts among the larger class of almost-planar fronts with a single or a finite number of interfaces between the limiting values $0$ and $1$. In other words, we prove the uniqueness of the transition fronts in this class. Firstly, we prove that, for the transition fronts whose level sets are almost-planar and orthogonal to the directions $e_t$ in the sense of Proposition~\ref{proplanar} and Theorem~\ref{thplanar}, the directions $e_t=e$ must be constant. Once this is done, the stability of the one-dimensional fronts and the one-dimensional symmetry of the almost-planar fronts moving with constant speed in a constant direction will lead to the conclusion of Proposition~\ref{proplanar}. For Theorem~\ref{thplanar}, one needs to exclude the case of transition fronts with 2 or more oscillations in the direction $e$. The proof will use the fact that the planar speed $c_f$ is assumed to be nonzero and that initial conditions which are above $\theta^+$ on a large set spread with speed $c_f$ at large times (if $c_f>0$).


\subsection{Proof of Proposition~\ref{proplanar}}\label{sec31}

We first begin with the following elementary lemma.

\begin{lem}\label{lem1}
Let $u:\R\times\R^N\to[0,1]$ be a solution of~$(\ref{eq})$ for which there are a real number~$t_0\in\R$ and a unit vector $e\in\mathbb{S}^{N-1}$ such that
\be\label{t0}
\inf_{x\in\R^N,\,x\cdot e\le-A}u(t_0,x)\to1\ \ \Big(\hbox{resp. }\sup_{x\in\R^N,\,x\cdot e\ge A}u(t_0,x)\to0\Big)\ \hbox{ as }A\to+\infty.
\ee
Then property~$(\ref{t0})$ holds at every time $t_1>t_0$ with the same vector $e$.
\end{lem}

\noindent{\bf{Proof.}} For any $\vartheta\in(0,1)$, let $\underline{v}_{\vartheta}$ and $\overline{v}_{\vartheta}$ be the solutions of the one-dimensional Cauchy problem
$$v_t=v_{yy}+f(v),\ \ t>0,\ y\in\R$$
with initial conditions
\be\label{defvalpha}
\underline{v}_{\vartheta}(0,y)=\left\{\baa{ll}\vartheta & \hbox{if }y\le0,\vspace{3pt}\\ 0 & \hbox{if }y>0,\eaa\right.\hbox{ and }\ \overline{v}_{\vartheta}(0,y)=\left\{\baa{ll}1 & \hbox{if }y\le0,\vspace{3pt}\\ \vartheta & \hbox{if }y>0.\eaa\right.
\ee
Let $\varrho_{\vartheta}:\R\to(0,1)$ denote the solution of the equation $\varrho_{\vartheta}'(t)=f(\varrho_{\vartheta}(t))$ with initial condition~$\varrho_{\vartheta}(0)=\vartheta$. It follows from the maximum principle and standard parabolic estimates that, for each $t>0$, $\underline{v}_{\vartheta}(t,\cdot)$ and $\overline{v}_{\vartheta}(t,\cdot)$ are decreasing in $\R$ and that $\underline{v}_{\vartheta}(t,-\infty)=\varrho_{\vartheta}(t)$, $\underline{v}_{\vartheta}(t,+\infty)=0$, $\overline{v}_{\vartheta}(t,-\infty)=1$ and $\overline{v}_{\vartheta}(t,+\infty)=\varrho_{\vartheta}(t)$.\par
Let now $u$, $e$ and $t_0<t_1$ be as in the lemma. We assume first that
$$\inf_{x\in\R^N,\,x\cdot e\le-A}\,u(t_0,x)\to1\hbox{ as }A\to+\infty.$$
Let $\epsilon\in(0,1)$ be arbitrary. There is $M\in\R$ such that
$$u(t_0,x)\ge\underline{v}_{1-\epsilon}(0,x\cdot e+M)\hbox{ for all }x\in\R^N,$$
whence $u(t_1,x)\ge\underline{v}_{1-\epsilon}(t_1-t_0,x\cdot e+M)$ for all $x\in\R^N$ from the maximum principle. Therefore,
$$\liminf_{A\to+\infty}\Big(\inf_{x\in\R^N,\,x\cdot e\le-A}u(t_1,x)\Big)\ge\underline{v}_{1-\epsilon}(t_1-t_0,-\infty)=\varrho_{1-\epsilon}(t_1-t_0).$$
Since this holds for all $\epsilon>0$ small enough, since $\varrho_{1-\epsilon}(t_1-t_0)\to1$ as $\epsilon\to0$ and since $u$ ranges in $(0,1)$, it follows that
$$\inf_{x\in\R^N,\,x\cdot e\le-A}u(t_1,x)\to1\ \hbox{ as }A\to+\infty.$$\par
Similarly, if
$$\sup_{x\in\R^N,\,x\cdot e\ge A}u(t_0,x)\to0\hbox{ as }A\to+\infty,$$
one gets that
$$\limsup_{A\to+\infty}\big(\sup_{x\in\R^N,\,x\cdot e\ge A}u(t_1,x)\big)\le\overline{v}_{\epsilon}(t_1-t_0,+\infty)=\varrho_{\epsilon}(t_1-t_0)$$
for all $\epsilon>0$ small enough, and the conclusion follows.\hfill$\Box$\break

From the previous lemma, the next result follows immediately.

\begin{cor}\label{cor1}
Let $u:\R\times\R^N\to(0,1)$ be a solution of~$(\ref{eq})$ such that, for every time $t\in\R$, there is a unit vector $e_t\in\mathbb{S}^{N-1}$ such that
\be\label{et}
\inf_{x\in\R^N,\,x\cdot e_t\le-A}u(t,x)\to1\ \hbox{ and }\ \sup_{x\in\R^N,\,x\cdot e_t\ge A}u(t,x)\to0\ \hbox{ as }A\to+\infty.
\ee
Then $e_t=e$ is independent of time $t$.
\end{cor}

\noindent{\bf{Proof of Proposition~\ref{proplanar}.}} Let $u$ be an almost-planar transition front connecting $0$ and $1$, in the sense of Definition~\ref{defap}, for problem~(\ref{eq}). That is, there exist some families~$(e_t)_{t\in\R}$ in~$\mathbb{S}^{N-1}$ and~$(\xi_t)_{t\in\R}$ in~$\R$ such that
$$\Gamma_t=\big\{x\in\R^N;\ x\cdot e_t=\xi_t\big\}$$
for every $t\in\R$. Up to changing $e_t$ into $-e_t$, it follows from~(\ref{omegapm}) and Definition~\ref{def1} that~(\ref{et}) holds for every $t\in\R$. Corollary~\ref{cor1} implies that $e_t=e$ is a constant vector, whence
\be\label{omegapmt}
\Omega^+_t=\big\{x\in\R^N;\ x\cdot e<\xi_t\big\}\ \hbox{ and }\ \Omega^-_t=\big\{x\in\R^N;\ x\cdot e>\xi_t\big\}
\ee
for all $t\in\R$.\par
We shall now prove that the function $\R\ni t\mapsto\xi_t-c_ft$ is bounded. To do so, we use the stability result of planar fronts of Fife and McLeod~\cite{fm} that we recalled in Section~\ref{sec11}. The last step of the proof will be based on a Liouville-type result of Berestycki and the author~\cite{bh4} on the characterization of the almost-planar transition fronts with constant direction and position of the order $c_ft$.\par
Let $\alpha$ and $\beta$ be two given real numbers such that
\be\label{alphabeta0}
0<\alpha<\theta^-\le\theta^+<\beta<1,
\ee
where we recall that $\theta^{\pm}$ are defined in~(\ref{thetapm}). From the previous observations and from Definition~\ref{def1}, there is $M\ge0$ such that
\be\label{defM}
\forall\,(t,x)\in\R\times\R^N,\ \left\{\baa{lcl}
x\cdot e-\xi_t\le-M & \Longrightarrow & \beta\le u(t,x)<1,\vspace{3pt}\\
x\cdot e-\xi_t\ge M & \Longrightarrow & 0<u(t,x)\le\alpha.\eaa\right.
\ee
Therefore, with the notations used in the proof of Lemma~\ref{lem1}, one infers that, for every $t_0\in\R$ and $x\in\R^N$,
$$\underline{v}_{\beta}(0,x\cdot e-\xi_{t_0}+M)\le u(t_0,x)\le\overline{v}_{\alpha}(0,x\cdot e-\xi_{t_0}-M).$$
Thus,
\be\label{comparaison}
\underline{v}_{\beta}(t-t_0,x\cdot e-\xi_{t_0}+M)\le u(t,x)\le\overline{v}_{\alpha}(t-t_0,x\cdot e-\xi_{t_0}-M)
\ee
for all $t>t_0$ and $x\in\R^N$, from the maximum principle.\par
On the other hand, from~(\ref{f}) and the existence of a planar front $(c_f,\phi_f)$ solving~(\ref{eqphi}), it follows from~\cite{fm} that there exist two real numbers $\underline{\xi}=\underline{\xi}(f,\beta)$ and $\overline{\xi}=\overline{\xi}(f,\alpha)$ depending only on $f$, $\alpha$ and $\beta$, such that
\be\label{valphabeta}
\sup_{y\in\R}\big|\underline{v}_{\beta}(s,y)-\phi_f(y-c_fs+\underline{\xi})\big|+\sup_{y\in\R}\big|\overline{v}_{\alpha}(s,y)-\phi_f(y-c_fs+\overline{\xi})\big|\to0\hbox{ as }s\to+\infty.
\ee
In particular, since $\phi_f(-\infty)=1$ and $\phi_f(+\infty)=0$, there exist $T>0$ and $A>0$ such that, for all $s\ge T$,
$$\left\{\baa{ll}
\underline{v}_{\beta}(s,y)>\alpha & \hbox{if }y\le c_fs-A,\vspace{3pt}\\
\overline{v}_{\alpha}(s,y)<\beta & \hbox{if }y\ge c_fs+A.\eaa\right.$$
Together with~(\ref{comparaison}), it follows that, for all $t_0<t_0+T\le t$,
\be\label{alpha}\left\{\baa{ll}
u(t,x)>\alpha & \hbox{if }x\cdot e-\xi_{t_0}+M\le c_f(t-t_0)-A,\vspace{3pt}\\
u(t,x)<\beta & \hbox{if }x\cdot e-\xi_{t_0}-M\ge c_f(t-t_0)+A.\eaa\right.
\ee
Properties~(\ref{defM}) and~(\ref{alpha}) imply that, for all $t_0<t_0+T\le t$,
\be\label{xitt0}
\xi_{t_0}-M+c_f(t-t_0)-A<\xi_t+M\ \hbox{ and }\ \xi_{t_0}+M+c_f(t-t_0)+A>\xi_t-M.
\ee
By fixing $t=0$, one gets that $\limsup_{t_0\to-\infty}|\xi_{t_0}-c_ft_0|\le|\xi_0|+2M+A$. For any arbitrary~$t\in\R$, letting $t_0\to-\infty$ in~(\ref{xitt0}) then leads to
$$|\xi_t-c_ft|\le|\xi_0|+4M+2A.$$
Therefore, Definition~\ref{def1} together with~(\ref{omegapmt}) implies that our solution $u:\R\times\R^N\to(0,1)$ of~(\ref{eq}) satisfies 
$$\inf_{(t,x)\in\R\times\R^N,\,x\cdot e-c_ft\le-A}u(t,x)\to1\ \hbox{ and }\ \sup_{(t,x)\in\R\times\R^N,\,x\cdot e-c_ft\ge A}u(t,x)\to0\ \hbox{ as }A\to+\infty.$$
It follows finally from Theorem~3.1 of~\cite{bh4} and the uniqueness of the planar fronts that there exists $\xi\in\R$ such that $u(t,x)=\phi_f(x\cdot e-c_ft+\xi)$ for all $(t,x)\in\R\times\R^N$. The proof of Proposition~\ref{proplanar} is thereby complete.\hfill$\Box$


\subsection{Planar and almost-planar transition fronts}\label{sec32}

In this section, we show that, without assuming a priori the existence of a one-dimensional planar front~$(c_f,\phi_f)$ solving~$(\ref{eqphi})$, the existence of an almost-planar transition front connec\-ting~$0$ and~$1$ for problem~$(\ref{eq})$ in $\R^N$ actually implies (and is then equivalent to) the existence of such a one-dimensional planar front~$(c_f,\phi_f)$, at least for a $C^1([0,1])$ dense set of functions~$f$. Namely, we prove the following result.

\begin{pro}\label{proplanar2}
For any $C^1([0,1])$ function $f$ satisfying~$(\ref{f})$ and for any $\epsilon>0$, there is a~$C^1([0,1])$ function $f_{\epsilon}$ satisfying~$(\ref{f})$ such that $\|f-f_{\epsilon}\|_{C^1([0,1])}\le\epsilon$ and for which the following holds: if $($and only if$)$ there exists an almost-planar transition front $u$ connecting $0$ and $1$ for problem~$(\ref{eq})$ with $f_{\epsilon}$, then there exists a planar front~$(c_{f_{\epsilon}},\phi_{f_{\epsilon}})$ solving~$(\ref{eqphi})$, and then~$u$ is a planar front of the type~$(\ref{planar})$, with $f_{\epsilon}$ instead of $f$. 
\end{pro}

\noindent{\bf{Proof.}} Let $f$ be any given $C^1([0,1])$ function satisfying~$(\ref{f})$ and let $\epsilon>0$ be arbitrary. Firstly, it is straightforward to check that there is a $C^1([0,1])$ function~$g$ satisfying~$(\ref{f})$ and
$$\|f-g\|_{C^1([0,1])}\le\frac{\epsilon}{2},$$
and for which there exist $k\in\N$ and some real numbers
$$0=\theta_0<\theta_1<\cdots<\theta_{2k-1}<\theta_{2k}=1$$
with
\be\label{thetai}\left\{\baa{ll}
g(\theta_i)=0 & \hbox{for all }0\le i\le 2k,\vspace{3pt}\\
g<0\hbox{ on }(\theta_{2i},\theta_{2i+1}) & \hbox{for all }0\le i\le k-1,\vspace{3pt}\\
g>0\hbox{ on }(\theta_{2i+1},\theta_{2i+2}) & \hbox{for all }0\le i\le k-1,\vspace{3pt}\\
g'(\theta_{2i})<0 & \hbox{for all }0\le i\le k,\vspace{3pt}\\
g'(\theta_{2i+1})>0 & \hbox{for all }0\le i\le k-1.\eaa\right.
\ee
In particular, the restriction of the function~$g$ on each interval $[\theta_{2i},\theta_{2i+2}]$ is of the bistable type. Notice also that
\be\label{thetapmbis}
\theta^-=\theta_1\hbox{ and }\theta^+=\theta_{2k-1}
\ee
for this function $g$, with the notation~$(\ref{thetapm})$.\par
As far as the existence of a planar connection~$(c_g,\phi_g)$ between $0$ and $1$ for this function~$g$ is concerned, two and only two cases occur, as follows from Fife and McLeod~\cite{fm}:
\begin{itemize}
\item (a)~either there is a (unique) planar front~$(c_g,\phi_g)$ solving~$(\ref{eqphi})$ with $g$ in place of~$f$,
\item (b)~or the pair~$(c_g,\phi_g)$ does not exist and there exist some integers~$l\in\{2,\cdots,k\}$ and~$0=i_0<i_1<\cdots<i_l=k$ such that there exists a (unique) planar front~$(\gamma_j,\varphi_j)$ connecting $\theta_{2i_{j-1}}$ and $\theta_{2i_j}$ for $g$ and for every~$1\le j\le l$, with
$$\gamma_1\ge\gamma_2\ge\cdots\ge\gamma_{l-1}\ge\gamma_l.$$
The fact that $\varphi_j$ connects $\theta_{2i_{j-1}}$ and $\theta_{2i_j}$ with the speed $\gamma_j$ for the function $g$ means that the pair~$(\gamma_j,\varphi_j)$ solves~$(\ref{eqphi})$ with the limits $\varphi_j(-\infty)=\theta_{2i_j}$ and $\varphi_j(+\infty)=\theta_{2i_{j-1}}$, and with~$(\gamma_j,\varphi_j,g)$ in place of $(c_f,\phi_f,f)$.
\end{itemize}\par
In case (a), we simply set $f_{\epsilon}=g$ and the conclusion of Proposition~\ref{proplanar2} follows trivially, since there is a planar connection $(c_{f_{\epsilon}},\phi_{f_{\epsilon}})$ solving~(\ref{eqphi}) with the function $f_{\epsilon}$.\par
Consider now the case (b). As it follows from~\cite{fm,t4}, for every $1\le j\le l$, there is a real number~$\eta_j>0$ such that for every
$$\widetilde{g}\in B_j:=\Big\{\widetilde{g}\in C^1([\theta_{2i_{j-1}},\theta_{2i_j}]);\ \|g-\widetilde{g}\|_{C^1([\theta_{2i_{j-1}},\theta_{2i_j}])}\le\eta_j,\ \widetilde{g}(\theta_{2i_{j-1}})=\widetilde{g}(\theta_{2i_j})=0\Big\}$$
there is a unique real number $\widetilde{\gamma}$ which is the speed of a planar connection between $\theta_{2i_{j-1}}$ and~$\theta_{2i_j}$ for the function $\widetilde{g}$; furthermore, the map $\widetilde{g}\mapsto\widetilde{\gamma}$ is continuous in $B_j$ endowed with the~$C^1([\theta_{2i_{j-1}},\theta_{2i_j}])$ norm and, if $\widetilde{g}$ and $\overline{g}$ belong to $B_j$ and satisfy $\widetilde{g}\le\not\equiv\overline{g}$ in $[\theta_{2i_{j-1}},\theta_{2i_j}]$, then the corresponding speeds satisfy $\widetilde{\gamma}<\overline{\gamma}$. Pick some arbitrary points $x_j\in(\theta_{2i_{j-1}},\theta_{2i_j})$ for every~$1\le j\le l$. By slightly changing the function $g$ locally around the points $x_j$ (namely by adding or subtracting some small nonnegative~$C^1$ functions supported in some small neighborhoods of~$x_j$), one infers straightforwardly that there is a $C^1([0,1])$ function $h$ satisfying~$(\ref{f})$ with $h$ instead of $f$, satisfying~$(\ref{thetai})$ with~$h$ instead of $g$, such that
$$\|g-h\|_{C^1([0,1])}\le\frac{\epsilon}{2},$$
and for which there are some planar fronts $(\widetilde{\gamma}_j,\widetilde{\varphi}_j)$ connecting~$\theta_{2i_{j-1}}$ and $\theta_{2i_j}$ for the function~$h$ and for every $1\le j\le l$, with
$$\widetilde{\gamma}_1>\widetilde{\gamma}_2>\cdots>\widetilde{\gamma}_{l-1}>\widetilde{\gamma}_l.$$\par
Finally, we will see that this last situation is incompatible with the existence of an almost-planar transition front connecting $0$ and $1$ for problem~$(\ref{eq})$ in $\R^N$ with $h$ in place of $f$. This will lead to the desired conclusion of Proposition~\ref{proplanar2}. Assume on the opposite that there is such an almost-planar transition front $u$ in the sense of Definition~\ref{defap}. Let $0<\alpha<\beta<1$ be given as in~(\ref{alphabeta0}) with the notation~(\ref{thetapmbis}), and let $\underline{v}_{\beta}$ and $\overline{v}_{\alpha}$ be the same as in the proof of Proposition~$\ref{proplanar}$ with the function $h$ instead of~$f$. Observe now that the proof of Lemma~$\ref{lem1}$ does not use the existence of a planar front solving~$(\ref{eqphi})$. Hence, Lemma~$\ref{lem1}$ and Corollary~$\ref{cor1}$ hold, that is $e_t=e$ is a constant vector, and one can also assume without loss of generality that all properties~$(\ref{et})$,~$(\ref{omegapmt})$,~$(\ref{defM})$ and~$(\ref{comparaison})$ hold with the functions $\underline{v}_{\beta}$ and $\overline{v}_{\alpha}$. On the other hand, it follows from~\cite{fm} that, instead of~$(\ref{valphabeta}),$\footnote{This was the place where the existence of $(c_f,\phi_f)$ played a role in the proof of Proposition~$\ref{proplanar}$.} there exist some real numbers $\underline{\xi}_j$ and~$\overline{\xi}_j$ for every~$1\le j\le l$, such that
$$\left\{\baa{rl}
\displaystyle{\mathop{\sup}_{y\le(\widetilde{\gamma}_l+\widetilde{\gamma}_{l-1})s/2}}\big|\underline{v}_{\beta}(s,y)-\widetilde{\varphi}_l(y-\widetilde{\gamma}_ls+\underline{\xi}_l)\big|\,\displaystyle{\mathop{\longrightarrow}_{s\to+\infty}}0, & \vspace{3pt}\\
\displaystyle{\mathop{\sup}_{(\widetilde{\gamma}_j+\widetilde{\gamma}_{j+1})s/2\le y\le(\widetilde{\gamma}_j+\widetilde{\gamma}_{j-1})s/2}}\big|\underline{v}_{\beta}(s,y)-\widetilde{\varphi}_j(y-\widetilde{\gamma}_js+\underline{\xi}_j)\big|\,\displaystyle{\mathop{\longrightarrow}_{s\to+\infty}}0 & \hbox{ for every }2\le j\le l-1,\vspace{3pt}\\
\displaystyle{\mathop{\sup}_{y\ge(\widetilde{\gamma}_1+\widetilde{\gamma}_2)s/2}}\big|\underline{v}_{\beta}(s,y)-\widetilde{\varphi}_1(y-\widetilde{\gamma}_1s+\underline{\xi}_1)\big|\,\displaystyle{\mathop{\longrightarrow}_{s\to+\infty}}0 & \eaa\right.$$
and
$$\left\{\baa{rl}
\displaystyle{\mathop{\sup}_{y\le(\widetilde{\gamma}_l+\widetilde{\gamma}_{l-1})s/2}}\big|\overline{v}_{\alpha}(s,y)-\widetilde{\varphi}_l(y-\widetilde{\gamma}_ls+\overline{\xi}_l)\big|\,\displaystyle{\mathop{\longrightarrow}_{s\to+\infty}}0, & \vspace{3pt}\\
\displaystyle{\mathop{\sup}_{(\widetilde{\gamma}_j+\widetilde{\gamma}_{j+1})s/2\le y\le(\widetilde{\gamma}_j+\widetilde{\gamma}_{j-1})s/2}}\big|\overline{v}_{\alpha}(s,y)-\widetilde{\varphi}_j(y-\widetilde{\gamma}_js+\overline{\xi}_j)\big|\,\displaystyle{\mathop{\longrightarrow}_{s\to+\infty}}0 & \hbox{ for every }2\le j\le l-1,\vspace{3pt}\\
\displaystyle{\mathop{\sup}_{y\ge(\widetilde{\gamma}_1+\widetilde{\gamma}_2)s/2}}\big|\overline{v}_{\alpha}(s,y)-\widetilde{\varphi}_1(y-\widetilde{\gamma}_1s+\overline{\xi}_1)\big|\,\displaystyle{\mathop{\longrightarrow}_{s\to+\infty}}0. & \eaa\right.$$
In other words, the one-dimensional functions~$\underline{v}_{\beta}$ and $\overline{v}_{\alpha}$ expand as an ordered family (or a terrace, with the terminology used in a more general framework in~\cite{dgm}) of traveling fronts with ordered speeds. Since $u$ is trapped between some finite shifts of~$\underline{v}_{\beta}$ and $\overline{v}_{\alpha}$ from~$(\ref{comparaison})$, this will be in contradiction with the uniform boundedness of the transition zone where $u$ is, say, between~$\alpha$ and~$\beta$. Indeed, owing to the aforementioned properties of the fronts $(\widetilde{\gamma}_j,\widetilde{\varphi}_j)$, there exist then some real numbers $T>0$ and $A>0$ such that, for all $s\ge T$,
$$\left\{\baa{ll}
\underline{v}_{\beta}(s,y)>\theta_1=\theta^->\alpha & \hbox{if }y\le\widetilde{\gamma}_1s-A,\vspace{3pt}\\
\overline{v}_{\alpha}(s,y)<\theta_{2k-1}=\theta^+<\beta & \hbox{if }y\ge\widetilde{\gamma}_ls+A.\eaa\right.$$
Together with property~$(\ref{comparaison})$ applied with $t_0=0$, one infers that, for all $t\ge T$,
$$\left\{\baa{ll}
u(t,x)>\alpha & \hbox{if }x\cdot e-\xi_0+M\le\widetilde{\gamma}_1t-A,\vspace{3pt}\\
u(t,x)<\beta & \hbox{if }x\cdot e-\xi_0-M\ge\widetilde{\gamma}_lt+A.\eaa\right.$$
With~(\ref{defM}), it follows that
$$\xi_0-M+\widetilde{\gamma}_1t-A<\xi_t+M\ \hbox{ and }\ \xi_0+M+\widetilde{\gamma}_lt+A>\xi_t-M$$
for all $t\ge T$. This leads to a contradiction as $t\to+\infty$, since $\widetilde{\gamma}_1>\widetilde{\gamma}_l$.\par
Therefore, in case (b) for $g$, we set $f_{\epsilon}=h$ and the conclusion of Proposition~\ref{proplanar2} follows since in this case there is no almost-planar transition front connecting $0$ and $1$ for problem~(\ref{eq}) with the function~$f_{\epsilon}$.\hfill$\Box$


\subsection{Proof of Theorem~\ref{thplanar}}

Although the statement of Theorem~\ref{thplanar} looks similar to that of Proposition~\ref{proplanar}, the proof is much more involved, not to mention that it requires necessarily that the planar speed $c_f$ be nonzero, as emphasized in Section~\ref{sec21}. The proof is based on a series of auxiliary lemmas establishing some bounds on the largest and/or smallest positive and/or negative parts of the positions~$\xi_{i,t}$ given in~(\ref{xiit}) of the interfaces along the direction $e_t=e$ (the direction~$e_t$ is easily seen to be independent of time). The bounds on the positions~$\xi_{i,t}$ rely crucially on the fact that the one-dimensional solutions of the Cauchy problem associated to~(\ref{eq}) in $\R$ with compactly supported initial conditions being above $\theta^+$ on a large set spread at the speed $c_f$ if $c_f>0$.\par 
Let $u$ be as in the statement of Theorem~\ref{thplanar} and assume that $c_f\neq 0$. Let $\widetilde{u}$, $g$ and $\phi_g$ be the functions defined by
$$\widetilde{u}(t,x)=1-u(t,x)\hbox{ for all }(t,x)\in\R\times\R^N,\ \ g(s)=-f(1-s)\hbox{ for all }s\in[0,1]$$
and $\phi_g(x)=1-\phi_f(-x)$ for all $x\in\R$. The function~$\widetilde{u}$ obeys the equation~(\ref{eq}) with $g$ instead of $f$, while $\phi_g(-\infty)=1>\phi_g(x)>\phi_g(+\infty)=0$ for all $x\in\R$ and
$$\phi_g''+c_g\phi_g'+g(\phi_g)=0\hbox{ in }\R$$
with $c_g=-c_f$. Therefore, even if it means replacing $u$ by $\widetilde{u}$, $f$ by $g$ and $c_f$ by $-c_f=c_g$, one can assume without loss of generality that $c_f>0$.\par
Definition~\ref{def1} and the assumptions made in Theorem~\ref{thplanar} imply that, for every $t\in\R$, both sets~$\Omega^+_t$ and~$\Omega^-_t$ contain a half-space. Therefore, up to changing $e_t$ into $-e_t$, one can assume without loss of generality that condition~(\ref{et}) is fulfilled for every $t\in\R$. It follows then from Corollary~\ref{cor1} that $e_t=e$ is a constant vector and that~(\ref{et}) holds with $e_t=e$ for every~$t\in\R$.\par
Even if it means reordering the real numbers $\xi_{i,t}$ given in~(\ref{xiit}), we denote, for every $t\in\R$,
\be\label{gammant}
\Gamma_t=\bigcup_{i=1}^{n_{t}}\big\{x\in\R^N;\ x\cdot e=\xi_{i,t}\big\},
\ee
with $n_t\in\N$ and $\xi_{1,t}<\cdots<\xi_{n_t,t}$. In particular,
$$\Omega^+_t\supset\big\{x\in\R^N;\ x\cdot e<\xi_{1,t}\big\}$$
and
\be\label{omega-t}
\Omega^-_t\supset\big\{x\in\R^N;\ x\cdot e>\xi_{n_t,t}\big\}.
\ee\par
In the following lemmas, we show some estimates for the positions (along the vector $e$) of the leftmost and rightmost interfaces $\xi_{1,t}$ and $\xi_{n_t,t}$, as well as $\xi^{\pm}_t$ defined in~(\ref{defxipmt}) below, as~$t\to-\infty$. These estimates are based in particular on the spreading properties of the solutions of the Cauchy problem associated to~(\ref{eq}) in $\R$ with sufficiently large initial conditions. They also use the fact that, from Definition~\ref{def1}, in some not-too-far neighborhoods of the hyperplanes~$\Gamma_t$, there are big regions where $u$ is close to $0$ and other ones where $u$ is close to $1$. Finally, putting together all the estimates, we show that the solution $u$ has essentially one interface, that has to move with the speed $c_f$, and then the solution~$u$ has to be a planar front.

\begin{lem}\label{lemn_{t}}
There holds
$$\liminf_{t\to-\infty}\,\big(\xi_{n_t,t}-c_ft\big)>-\infty\ \hbox{ and }\ \limsup_{t\to-\infty}\,\big(\xi_{1,t}-c_ft\big)<+\infty.$$
\end{lem}

\noindent{\bf{Proof.}} Assume that the first conclusion is not satisfied. Then there is a sequence $(t_k)_{k\in\N}$ of real numbers such that
$$t_k\to-\infty\hbox{ and }\xi_{n_{t_k},t_k}-c_ft_k\to-\infty\hbox{ as }k\to+\infty.$$
Let $\alpha\in(0,\theta^-)$ be given, where we recall that $\theta^-$ is defined in~(\ref{thetapm}). From Definition~\ref{def1} and~(\ref{omega-t}), there is $M_{\alpha}\ge 0$ such that
$$u(t,x)\le\alpha\hbox{ for all }(t,x)\in\R\times\R^N\hbox{ with }x\cdot e-\xi_{n_t,t}\ge M_{\alpha}.$$
Therefore, definition~(\ref{defvalpha}) of $\overline{v}_{\alpha}$ yields
$$u(t_k,x)\le\overline{v}_{\alpha}(0,x\cdot e-\xi_{n_{t_k},t_k}-M_{\alpha})\hbox{ for all }x\in\R^N,$$
whence
$$u(t,x)\le\overline{v}_{\alpha}(t-t_k,x\cdot e-\xi_{n_{t_k},t_k}-M_{\alpha})\hbox{ for all }x\in\R^N$$
and for all $t>t_k$. For any fixed $(t,x)\in\R\times\R^N$, since
$$\lim_{k\to+\infty}t_k=\lim_{k\to+\infty}\big(\xi_{n_{t_k},t_k}-c_ft_k\big)=-\infty,$$
it follows then from the existence of a planar front~$(c_f,\phi_f)$ and from~(\ref{valphabeta}) that, for some~$\overline{\xi}=\overline{\xi}(f,\alpha)\in\R$,
$$u(t,x)\le\limsup_{k\to+\infty}\,\phi_f\big(x\cdot e-\xi_{n_{t_k},t_k}-M_{\alpha}-c_f(t-t_k)+\overline{\xi}\big)=\phi_f(+\infty)=0.$$
This is impossible, since, as already emphasized in the introduction, $u$ is assumed to range in the open interval $(0,1)$.\par
Similarly, if there is a sequence $(t_k)_{k\in\N}$ of real numbers such that $t_k\to-\infty$ and~$\xi_{1,t_k}-c_ft_k\to+\infty$ as $k\to+\infty$, one would get that $u(t,x)\ge 1$ for all $(t,x)\in\R\times\R^N$, which is ruled out. Finally, $\limsup_{t\to-\infty}\big(\xi_{1,t}-c_ft\big)<+\infty$ and the proof of Lemma~\ref{lemn_{t}} is thereby complete.~\hfill$\Box$\break

Let us now introduce a few additional useful notations. Let $\beta$ be any given real number such that
$$\theta^+<\beta<1,$$
where we recall that $\theta^+$ is defined in~(\ref{thetapm}). Due to the existence of a planar front~$(c_f,\phi_f)$ solving~(\ref{eqphi}) with $c_f>0$, it follows from~\cite{fm} that there are $A>0$ and $\sigma\in\R$ such that the solution $w$ of the Cauchy problem
$$w_t=w_{yy}+f(w),\ \ t>0,\ y\in\R,$$
with initial condition
\be\label{defw0}
w(0,y)=\left\{\baa{ll}
\beta & \hbox{if }|y|\le A,\vspace{3pt}\\
0 & \hbox{if }|y|>A,\eaa\right.
\ee
is such that $w(t,y)\to1$ as $t\to+\infty$ locally uniformly in $y\in\R$, and moreover
\be\label{convw}
\sup_{y\in\R}\,\big|w(s,y)-\phi_f(y-c_fs+\sigma)-\phi_f(-y-c_fs+\sigma)+1\big|\to0\hbox{ as }s\to+\infty.
\ee\par
From Definition~\ref{def1}, there is $M_{\beta}\ge 0$ such that
\be\label{defMbis}
\forall\,(t,x)\in\R\times\R^N,\ \ \big(x\in\Omega^+_t\hbox{ and }d(x,\Gamma_t)\ge M_{\beta}\big)\Longrightarrow\big(\beta\le u(t,x)<1\big).
\ee
Without loss of generality, one can assume that $M_{\beta}\ge A$. From condition~(\ref{unifgamma}), there is then~$r>0$ such that
\be\label{defr}
\forall\,t\in\R,\ \forall\,x\in\Gamma_t,\ \exists\,y\in\Omega^+_t,\ |y-x|\le r,\ d(y,\Gamma_t)\ge 2M_{\beta}.
\ee\par
Lastly, let $\xi^{\pm}_{t}$ be defined in $[-\infty,+\infty]$ as
\be\label{defxipmt}
\xi^-_{t}=\sup\big\{x\cdot e;\ x\in\Gamma_t,\ x\cdot e\le 0\big\}\ \hbox{ and }\ \xi^+_{t}=\inf\big\{x\cdot e;\ x\in\Gamma_t,\ x\cdot e\ge 0\big\}.
\ee
If $\xi_{1,t}\le0$ (this holds for $t$ negative enough from the previous lemma), then $\xi^-_t$ is a real number and $\xi^-_t=\xi_{n^-_t\!,t}$ for some $1\le n^-_t\le n_t$, otherwise $\xi^-_t=-\infty$. If~$\xi_{n_t,t}\ge0$, then~$\xi^+_t$ is a real number and $\xi^+_t=\xi_{n^+_t\!,t}$ for some $1\le n^+_{t}\le n_t$, otherwise~$\xi^+_t=+\infty$.

\begin{lem}\label{lemaw}
One has
\be\label{xipmt}
\limsup_{t\to-\infty}\,\big(\xi^-_{t}-c_ft)<+\infty,\ \ \liminf_{t\to-\infty}\,\big(\xi^+_{t}+c_ft)>-\infty
\ee
and there is $T_1\in\R$ such that
\be\label{Et}
\Omega^-_t\supset\big\{x\in\R^N;\ \xi^-_{t}<x\cdot e<\xi^+_{t}\big\}=:E_t
\ee
for all $t\le T_1$.
\end{lem}

\noindent{\bf{Proof.}} Assume first that $\limsup_{t\to-\infty}\big(\xi^-_{t}-c_ft)=+\infty$. Then there is a sequence $(t_k)_{k\in\N}$ of real numbers such that
$$t_k\to-\infty\hbox{ and }\xi^-_{t_k}-c_ft_k\to+\infty\hbox{ as }k\to+\infty.$$
Since $\xi^-_{t_k}\in[-\infty,0]$ for all $k\in\N$ and since $\xi^-_{t_k}-c_ft_k\to+\infty$ as $k\to+\infty$, one can assume without loss of generality that $\xi^-_{t_k}\in\R$ and
$$\big\{x\in\R^N;\ x\cdot e=\xi^-_{t_k}\big\}\subset\Gamma_{t_k}\hbox{ for all }k\in\N.$$
For every $k\in\N$, let now $x_k=\xi^-_{t_k}\,e\in\Gamma_{t_k}$ and, from~(\ref{defr}), let $y_k\in\Omega^+_{t_k}$ such that~$|y_k-x_k|\le r$ and $d(y_k,\Gamma_{t_k})\ge 2M_{\beta}$. Set $\omega_k=y_k\cdot e$. There holds
\be\label{xi-tk}
\xi^-_{t_k}-r\le\omega_k\le\xi^-_{t_k}+r\hbox{ for all }k\in\N.
\ee
Furthermore, since $B(y_k,2M_{\beta})\subset\Omega^+_{t_k}$ and $\Gamma_{t_k}=\bigcup_{1\le i\le n_{t_k}}\big\{x\in\R^N;\ x\cdot e=\xi_{i,t_k}\big\}$, it follows that
$$\big\{x\in\R^N;\ \omega_k-2M_{\beta}<x\cdot e<\omega_k+2M_{\beta}\big\}\subset\Omega^+_{t_k}\hbox{ for all }k\in\N.$$
Property~(\ref{defMbis}) implies that, for every $k\in\N$,
\be\label{utk}
u(t_k,x)\ge\beta\hbox{ for all }x\in\R^N\hbox{ such that }\omega_k-M_{\beta}\le x\cdot e\le\omega_k+M_{\beta}.
\ee
The definition of $w(0,\cdot)$ in~(\ref{defw0}) and the inequality $M_{\beta}\ge A$ yield $u(t_k,x)\ge w(0,x\cdot e-\omega_k)$ for all $x\in\R^N$, whence
$$u(t,x)\ge w(t-t_k,x\cdot e-\omega_k)\ \hbox{ for all }t>t_k\hbox{ and }x\in\R^N,$$
for every $k\in\N$. For any fixed $(t,x)\in\R\times\R^N$, one gets from~(\ref{convw}) that
$$u(t,x)\ge\limsup_{k\to+\infty}\Big(\phi_f(x\cdot e-\omega_k-c_f(t-t_k)+\sigma)+\phi_f(-x\cdot e+\omega_k-c_f(t-t_k)+\sigma)-1\Big).$$
Notice now that $\omega_k-c_ft_k\to+\infty$ as $k\to+\infty$, since $\xi^-_{t_k}-c_ft_k\to+\infty$ as $k\to+\infty$ by assumption and since the sequence $(\omega_k-\xi^-_{t_k})_{k\in\N}$ is bounded from~(\ref{xi-tk}). On the other hand, since~$\xi^-_{t_k}\le 0<c_f$ and $t_k\to-\infty$ as $k\to+\infty$, one gets that $\xi^-_{t_k}+c_ft_k\to-\infty$ and~$\omega_k+c_ft_k\to-\infty$ as $k\to+\infty$. Finally,
$$\phi_f(x\cdot e-\omega_k-c_f(t-t_k)+\sigma)\to\phi(-\infty)=1\hbox{ as }k\to+\infty$$
and
$$\phi_f(-x\cdot e+\omega_k-c_f(t-t_k)+\sigma)\to\phi(-\infty)=1\hbox{ as }k\to+\infty,$$
whence $u(t,x)\ge 1$ for all $(t,x)\in\R\times\R^N$. One has then reached a contradiction. Therefore,~$\limsup_{t\to-\infty}\big(\xi^-_{t}-c_ft)<+\infty$.\par
Similarly, if one assumes that there is a sequence $(t_k)_{k\in\N}$ of real numbers such that $t_k\to-\infty$ and $\xi^+_{t_k}+c_ft_k\to-\infty$ as $k\to+\infty$ (which would imply in particular that $0\le\xi^+_{t_k}<+\infty$ for~$k$ large enough), one would get a similar contradiction (one can also apply the previous result to the function $(t,x)\mapsto u(t,\tilde{x})$, where $\tilde{x}=x-2(x\cdot e)e$ denotes the image of a point~$x$ by the orthogonal symmetry with respect to the hyperplane orthogonal to $e$ and containing the origin).\par
From Lemma~\ref{lemn_{t}}, it follows in particular that $\xi_{1,t}\to-\infty$ as $t\to-\infty$. Hence, there is~$T_0\in\R$ such that, for all $t\le T_0$, $\xi^-_{t}\in(-\infty,0]$ and either $E_t\subset\Omega^+_t$ or $E_t\subset\Omega^-_t$, under the notation given in~(\ref{Et}). Assume now by contradiction that there is a sequence $(t_k)_{k\in\N}$ of real numbers such that $t_k\to-\infty$ as $k\to+\infty$ and $E_{t_k}\subset\Omega^+_{t_k}$ for all $k\in\N$. Then~(\ref{xipmt}) and Definition~\ref{def1} imply in particular that
$$\inf_{x\in\R^N,\,|x\cdot e|\le A}\,u(t_k,x)\to1\hbox{ as }k\to+\infty.$$
Therefore, for $k$ large enough, $u(t_k,x)\ge w(0,x\cdot e)$ for all $x\in\R^N$, whence $u(t,x)\ge w(t-t_k,x\cdot e)$ for all $t>t_k$ and $x\in\R^N$. Finally, for any fixed $(t,x)\in\R\times\R^N$, one would get that
$$u(t,x)\ge\limsup_{k\to+\infty}\Big(\phi_f(x\cdot e-c_f(t-t_k)+\sigma)+\phi_f(-x\cdot e-c_f(t-t_k)+\sigma)-1\Big)=1$$
since $t_k\to-\infty$ as $k\to+\infty$. A contradiction has been reached, and the desired conclusion~(\ref{Et}) follows for $-t$ large enough.\hfill$\Box$\break

The following lemma is one of the key-points in the proof of Theorem~\ref{thplanar}. It states that it the positions $\xi^+_{t}$ are far on the right of the (very positive) position $c_f|t|$ along the direction~$e$, at least for a sequence of times $t_k$ converging to~$-\infty$, then they are actually pushed to~$+\infty$ and the solution $u$ has no interface far on the right of the (very negative) position $c_ft$ as $t\to-\infty$.

\begin{lem}\label{lemexp}
If $\limsup_{t\to-\infty}\big(\xi^+_{t}+c_ft\big)=+\infty$, then there are $T_2\in\R$, $\delta>0$ and $B\in\R$ such that
\be\label{upperexp}
u(t,x)\le\delta\,e^{-\delta(x\cdot e-c_ft-B)}\ \hbox{ for all }t\le T_2\hbox{ and }x\cdot e\ge c_ft+B.
\ee
\end{lem}

\begin{rem}{\rm 
With the same arguments as in the proof of Lemma~$\ref{lemexp}$, one can get the following result, which we state in a remark since it will actually not be used in the sequel: if~$\liminf_{t\to-\infty}\big(\xi^-_{t}-c_ft\big)=-\infty$, then there are $T'_2\in\R$, $\delta'>0$ and $B'\in\R$ such that
$$u(t,x)\le\delta'e^{\delta'(x\cdot e+c_ft+B')}\ \hbox{ for all }t\le T'_2\hbox{ and }x\cdot e\le -c_ft-B'.$$}
\end{rem}

The proof of Lemma~\ref{lemexp} is quite lengthy and is postponed in Section~\ref{sec33}. We prefer to go on the proof of Theorem~\ref{thplanar} with the following lemma.

\begin{lem}\label{lemphi}
If there is a sequence $(t_k)_{k\in\N}$ of real numbers such that $t_k\to-\infty$ as $k\to+\infty$ and $\liminf_{k\to+\infty}\big(\xi^-_{t_k}-c_ft_k\big)>-\infty$ $($resp. $\limsup_{k\to+\infty}\big(\xi^+_{t_k}+c_ft_k\big)<+\infty)$, then there is~$\eta\in\R$ such that
$$u(t,x)\ge\phi_f(x\cdot e-c_ft+\eta)\ \ \Big(\hbox{resp. }u(t,x)\ge\phi_f(-x\cdot e-c_ft+\eta)\Big)$$
for all $(t,x)\in\R\times\R^N$.
\end{lem}

Since $\phi_f(-\infty)=1$ and $u(t,x)\to0$ as $x\cdot e\to+\infty$ for every $t\in\R$ by~(\ref{omega-t}), the following corollary follows immediately.

\begin{cor}\label{corphi}
There holds $\xi^+_{t}+c_ft\to+\infty$ as $t\to-\infty$.
\end{cor}

\noindent{\bf{Proof of Lemma~\ref{lemphi}.}} Let $(t_k)_{k\in\N}$ be a sequence of real numbers such that $t_k\to-\infty$ as~$k\to+\infty$ and $\liminf_{k\to+\infty}\big(\xi^-_{t_k}-c_ft_k\big)>-\infty$. It follows then from Lemma~\ref{lemaw} that the sequence $(\xi^-_{t_k}-c_ft_k)_{k\in\N}$ is bounded. As in the proof of Lemma~\ref{lemaw}, there are $r>0$ and a sequence $(\omega_k)_{k\in\N}$ in $\R$ such that~(\ref{xi-tk}) and~(\ref{utk}) hold, whence $u(t_k,x)\ge w(0,x\cdot e-\omega_k)$ and~$u(t,x)\ge w(t-t_k,x\cdot e-\omega_k)$ for all $k\in\N$, $t>t_k$ and $x\in\R^N$. For any fixed $(t,x)\in\R\times\R^N$, one infers that
$$u(t,x)\ge\limsup_{k\to+\infty}\Big(\phi_f(x\cdot e-\omega_k-c_f(t-t_k)+\sigma)+\phi_f(-x\cdot e+\omega_k-c_f(t-t_k)+\sigma)-1\Big).$$
Since the sequences $(\xi^-_{t_k}-c_ft_k)_{k\in\N}$ and $(\omega_k-\xi^-_{t_k})_{k\in\N}$ are bounded, so is the sequence~$(\omega_k-c_ft_k)_{k\in\N}$ and one can assume, up to extraction of a sequence, that there is $\widetilde{\sigma}\in\R$ such that $\omega_k-c_ft_k\to\widetilde{\sigma}$ as $k\to+\infty$. On the other hand, $\omega_k+c_ft_k=\omega_k-c_ft_k+2c_ft_k\to-\infty$ since $t_k\to-\infty$ as~$k\to+\infty$. Finally, one concludes that, for any fixed $(t,x)\in\R\times\R^N$,
$$u(t,x)\ge\phi_f(x\cdot e-c_ft+\sigma-\widetilde{\sigma})+\phi_f(-\infty)-1=\phi_f(x\cdot e-c_ft+\sigma-\widetilde{\sigma}),$$
which gives the desired result with $\eta=\sigma-\widetilde{\sigma}$.\par
Similar arguments imply that, if $\limsup_{k\to+\infty}\big(\xi^+_{t_k}+c_ft_k\big)<+\infty$ for some sequence~$(t_k)_{k\in\N}$ converging to~$-\infty$, then $u(t,x)\ge\phi_f(-x\cdot e-c_ft+\eta)$ in $\R\times\R^N$ for some~$\eta\in\R$.~\hfill$\Box$\break

With all the previous lemmas, we are now ready to finish the proof of Theorem~\ref{thplanar}.\hfill\break

\noindent{\bf{End of the proof of Theorem~\ref{thplanar}.}} First, Lemma~\ref{lemexp} and Corollary~\ref{corphi} provide the existence of some $T_2\in\R$, $\delta>0$ and $B\in\R$ such that~(\ref{upperexp}) is satisfied.\par
We shall now prove that
\be\label{xintbis}
\limsup_{t\to-\infty}\,\big(\xi_{n_t,t}-c_ft\big)<+\infty.
\ee
Assume not. There is then a sequence $(t_k)_{k\in\N}$ of real numbers such that $t_k\to-\infty$ and~$\xi_{n_{t_k},t_k}-c_ft_k\to+\infty$ as $k\to+\infty$. Denote $x_k=\xi_{n_{t_k},t_k}e$. As in the proof of Lemma~\ref{lemaw}, there are $r>0$ and a sequence $(y_k)_{k\in\N}$ in $\R^N$ such that
$$y_k\in\Omega^+_{t_k},\ |y_k-x_k|\le r\hbox{ and }u(t_k,y_k)\ge\beta\hbox{ for all }k\in\N.$$
But the sequence $(y_k\cdot e-\xi_{n_{t_k},t_k})_{k\in\N}$ is bounded, whence $y_k\cdot e-c_ft_k\to+\infty$ as $k\to+\infty$. Therefore, $t_k\le T_2$ and $y_k\cdot e\ge c_ft_k+B$ for $k$ large enough. Finally,~(\ref{upperexp})  yields
$$u(t_k,y_k)\le\delta\,e^{-\delta(y_k\cdot e-c_ft_k-B)}$$
for $k$ large enough. The right-hand side converges to $0$ as $k\to+\infty$, whereas the left-hand side is bounded from below by $\beta>0$. One is led to a contradiction, and~(\ref{xintbis}) is proved.\par
Notice that~(\ref{xintbis}) implies in particular that $\xi_{n_t,t}<0$, $\xi^+_{t}=+\infty$ and $\xi^-_{t}=\xi_{n_t,t}$ for $t$ negative enough. Furthermore, together with Lemma~\ref{lemn_{t}}, one obtains that
\be\label{xi-xint}
\limsup_{t\to-\infty}\,\big|\xi^-_{t}-c_ft\big|=\limsup_{t\to-\infty}\,\big|\xi_{n_t,t}-c_ft\big|<+\infty.
\ee
Lemma~\ref{lemphi} provides then the existence of a real number $\eta$ such that
\be\label{uphi}
u(t,x)\ge\phi_f(x\cdot e-c_ft+\eta)\hbox{ for all }(t,x)\in\R\times\R^N.
\ee\par
On the other hand,~(\ref{xi-xint}) also yields the existence of a real number $\widetilde{\xi}$ and a sequence of real numbers $(t_k)_{k\in\N}$ such that $t_k\to-\infty$ and $\xi_{n_{t_k},t_k}-c_ft_k\to\widetilde{\xi}$ as $k\to+\infty$. As in the proof of Lemma~\ref{lemn_{t}}, there are then $\alpha\in(0,\theta^-)$ and $M_{\alpha}\ge0$ such that
$$u(t_k,x)\le\overline{v}_{\alpha}\big(0,x\cdot e-\xi_{n_{t_k},t_k}-M_{\alpha}\big)$$
whence $u(t,x)\le\overline{v}_{\alpha}(t-t_k,x\cdot e-\xi_{n_{t_k},t_k}-M_{\alpha})$ for all $k\in\N$, $t>t_k$ and $x\in\R^N$. Since~$\overline{v}_{\alpha}(s,y)-\phi_f(y-c_fs+\overline{\xi})\to0$ as $s\to+\infty$ uniformly in $y\in\R$, for some $\overline{\xi}\in\R$, one infers that, for any fixed $(t,x)\in\R\times\R^N$,
$$u(t,x)\le\limsup_{k\to+\infty}\phi_f(x\cdot e-\xi_{n_{t_k},t_k}-M_{\alpha}-c_f(t-t_k)+\overline{\xi})=\phi_f(x\cdot e-c_ft-\widetilde{\xi}-M_{\alpha}+\overline{\xi}).$$\par
As a conclusion,
$$\phi_f(x\cdot e-c_ft+\eta)\le u(t,x)\le\phi_f(x\cdot e-c_ft+\widetilde{\eta})$$
for all $(t,x)\in\R\times\R^N$, with $\widetilde{\eta}=\overline{\xi}-\widetilde{\xi}-M_{\alpha}$. As in the end of the proof of Proposition~\ref{proplanar}, Theorem~3.1 of~\cite{bh4} yields the existence of a real number $\xi$ such that
$$u(t,x)=\phi_f(x\cdot e-c_ft+\xi)\hbox{ for all }(t,x)\in\R\times\R^N.$$
The proof of Theorem~\ref{thplanar} is thereby complete.\hfill$\Box$

\begin{rem}{\rm 
It is immediate to see that the conclusion of Theorem~$\ref{thplanar}$ still holds even if the family of integers $(n_{t})_{t\in\R}$ appearing in~$(\ref{gammant})$ is not bounded. In other words, if $c_f\neq 0$ and if~$u$ satisfies all assumptions of Theorem~$\ref{thplanar}$ with the exception of~$(\ref{omegapmbis})$ and the boundedness of~$n_t$ in
$$\Gamma_t=\bigcup_{1\le i\le n_t}\big\{x\in\R^N;\ x\cdot e_t=\xi_{i,t}\big\},$$
then $u$ is still a planar front of the type $u(t,x)=\phi_f(x\cdot e-c_ft+\xi)$ for all $(t,x)\in\R\times\R^N$, for some unit vector $e$ of $\R^N$ and some real number $\xi$.}
\end{rem}


\subsection{Proof of Lemma~\ref{lemexp}}\label{sec33}

Let $(t_k)_{k\in\N}$ be a sequence of real numbers such that
\be\label{assumption}
t_k\to-\infty\hbox{ and }\xi^+_{t_k}+c_ft_k\to+\infty\hbox{ as }k\to+\infty.
\ee
Without loss of generality, one can assume in particular that
\be\label{xi+tk}
\xi^+_{t_k}+c_ft_k\ge0\hbox{ for all }k\in\N.
\ee
The strategy of the proof consists in constructing a sequence of supersolutions of~(\ref{eq}) which are approximately of the type $\delta\,e^{-\delta(x\cdot e-c_ft-B)}+\phi_f(-x\cdot e-c_f(t-t_k)+\xi^+_{t_k})$ (plus some shifts and some small exponential terms, as in the original proof of Fife and McLeod~\cite{fm}) for~$t_k\le t\le T_2$ and $x\cdot e-c_ft\ge B$. The parameters $T_2\in\R$, $\delta>0$ and $B\in\R$ will be chosen independently of $k$. The passage to the limit as $k\to+\infty$ in the supersolution will provide the desired conclusion, since $\xi^+_{t_k}+c_ft_k\to+\infty$.\par
Let us first choose some parameters. Remember that $f'(0)$ and $f'(1)$ are negative. Let~$\delta>0$ be such that
\be\label{defdelta}
0<\delta<\min\Big(1,\frac{|f'(0)|}{2},\frac{|f'(1)|}{2}\Big),\ \ f'\le\frac{f'(0)}{2}\hbox{ on }[0,3\delta]\ \hbox{ and }\ f'\le\frac{f'(1)}{2}\hbox{ on }[1-\delta,1].
\ee
Let $C>0$ be such that
\be\label{defC}
\phi_f\ge1-\delta\hbox{ on }(-\infty,-C]\ \hbox{ and }\ \phi_f\le\delta\hbox{ on }[C,+\infty).
\ee
Since $\phi'_f$ is negative and continuous on $\R$, there is $\kappa>0$ such that
\be\label{defkappa}
-\phi'_f\ge\kappa>0\hbox{ on }[-C,C].
\ee
Set $L=\max_{[0,1]}|f'|$ and let $\omega>0$ such that
\be\label{defomega}
\kappa\,\omega\ge L+\delta\ \hbox{ and }\ \kappa\,\omega\,c_f\ge L+\delta^2
\ee
(remember that $c_f$ is assumed to be positive in the proof of Theorem~\ref{thplanar}, without loss of generality). From Definition~\ref{def1} and Lemma~\ref{lemaw}, there are $T_2\le0$ and $B>0$ such that
\be\label{udelta}\left\{\baa{l}
c_ft+B<-c_ft-2B<\xi^+_{t}-B\ \hbox{ for all }t\le T_2,\vspace{3pt}\\
u(t,x)\le\delta\ \hbox{ for all }t\le T_2\hbox{ and }c_ft+B\le x\cdot e\le\xi^+_{t}-B,\eaa\right.
\ee
and
\be\label{defT2}
c_ft+2(\omega+B+C)\le 0\ \hbox{ for all }t\le T_2.
\ee\par
Without loss of generality, one can assume that $t_k<T_2$ for all $k\in\N$. For all $k\in\N$ and~$(t,x)\in\R\times\R^N$, let us now set
$$\overline{u}_k(t,x)=\min\Big(\phi_f\big(\zeta_k(t,x)\big)+\delta\,e^{-\delta(x\cdot e-c_ft-B)}+\delta\,e^{-\delta(t-t_k)},1\Big),$$
where
$$\zeta_k(t,x)=-x\cdot e-c_f(t-t_k)+\xi^+_{t_k}+\omega\,e^{-\delta(t-t_k)}-\omega-\omega\,e^{\delta c_f t}-B-C,$$
under the convention that $\zeta_k(t,x)=+\infty$ and $\overline{u}_k(t,x)=\min\big(\delta\,e^{-\delta(x\cdot e-c_ft-B)}+\delta\,e^{-\delta(t-t_k)},1\big)$ if~$\xi^+_{t_k}=+\infty$. Let us check that $\overline{u}_k$ is a supersolution of the equation~(\ref{eq}) satisfied by $u$, in the set
$$\mathcal{E}_k=\big\{(t,x)\in\R\times\R^N;\ t_k\le t\le T_2,\ x\cdot e\ge c_ft+B\big\}.$$\par
In what follows, $k$ denotes an arbitrary integer. At the time $t_k$, if follows from~(\ref{udelta}) and the definition of $\overline{u}_k$ that
$$u(t_k,x)\le\delta\le\overline{u}_k(t_k,x)\hbox{ for all }x\in\R^N\hbox{ such that }c_ft_k+B\le x\cdot e\le\xi^+_{t_k}-B.$$
On the other hand, if $x\cdot e\ge\xi^+_{t_k}-B$, then $\zeta_k(t_k,x)\le-\omega\,e^{\delta c_ft_k}-C\le-C$, whence
$$\overline{u}_k(t_k,x)\ge\min\big(\phi_f(\zeta_k(t_k,x))+\delta,1\big)\ge\min\big((1-\delta)+\delta,1\big)=1\ge u(t_k,x)$$
due to the first part of~(\ref{defC}). The above calculation holds whether $\xi^+_{t_k}$ and $\zeta_k(t,x)$ be real numbers or equal to~$+\infty$. As far as the boundary condition is concerned, if $t_k\le t\le T_2$ and~$x\cdot e=c_ft+B$, then
$$u(t,x)\le\delta\le\overline{u}_k(t,x)$$
from~(\ref{udelta}) and the definition of $\overline{u}_k$.\par
Let us then check that
$$N_k(t,x):=\overline{u}_k(t,x)-\Delta\overline{u}_k(t,x)-f(\overline{u}_k(t,x))\ge0\hbox{ for all }(t,x)\in\mathcal{E}_k\hbox{ such that }\overline{u}_k(t,x)<1.$$
This will be sufficient to ensure that $\overline{u}_k$ is a supersolution of~(\ref{eq}) in~$\mathcal{E}_k$, since $f(1)=0$.\par
In this paragraph, $(t,x)$ denotes any point in $\mathcal{E}_k$ such that $\overline{u}_k(t,x)<1$. From~(\ref{eqphi}), it is straightforward to check that
$$\baa{rcl}
N_k(t,x) & = & f(\phi_f(\zeta_k(t,x)))-f(\overline{u}_k(t,x))-\omega\,\delta\big(e^{-\delta(t-t_k)}+c_fe^{\delta c_ft}\big)\phi'_f(\zeta_k(t,x))\vspace{3pt}\\
& & +\delta^2(c_f-\delta)\,e^{-\delta(x\cdot e-c_ft-B)}-\delta^2e^{-\delta(t-t_k)}\vspace{3pt}\\
& \ge & f(\phi_f(\zeta_k(t,x)))-f(\overline{u}_k(t,x))-\omega\,\delta\big(e^{-\delta(t-t_k)}+c_fe^{\delta c_ft}\big)\phi'_f(\zeta_k(t,x))\vspace{3pt}\\
& & -\delta^3e^{-\delta(x\cdot e-c_ft-B)}-\delta^2e^{-\delta(t-t_k)},\eaa$$
whether $\xi^+_{t_k}$ and $\zeta_k(t,x)$ be real numbers or equal to~$+\infty$. If $\zeta_k(t,x)\le-C$, then~(\ref{defC}) implies that
$$1-\delta\le\phi_f(\zeta_k(t,x))\le\overline{u}_k(t,x)<1,$$
whence
$$f(\phi_f(\zeta_k(t,x)))-f(\overline{u}_k(t,x))\ge\frac{-f'(1)}{2}\,\big(\delta\,e^{-\delta(x\cdot e-c_ft-B)}+\delta\,e^{-\delta(t-t_k)}\big)$$
from~(\ref{defdelta}), and thus
$$N_k(t,x)\ge\delta\,\Big(\frac{-f'(1)}{2}-\delta^2\Big)\,e^{-\delta(x\cdot e-c_ft-B)}+\delta\,\Big(\frac{-f'(1)}{2}-\delta\Big)\,e^{-\delta(t-t_k)}\ge0$$
from~(\ref{defdelta}) and the negativity of $\phi'_f$. Similarly, if $\zeta_k(t,x)\ge C$, then $\phi_f(\zeta_k(t,x))\le\delta$. Thus,
$$0<\phi_f(\zeta_k(t,x))\le\overline{u}_k(t,x)\le3\delta$$
and
$$f(\phi_f(\zeta_k(t,x)))-f(\overline{u}_k(t,x))\ge\frac{-f'(0)}{2}\,\big(\delta\,e^{-\delta(x\cdot e-c_ft-B)}+\delta\,e^{-\delta(t-t_k)}\big)$$
from~(\ref{defdelta}), whence
$$N_k(t,x)\ge\delta\,\Big(\frac{-f'(0)}{2}-\delta^2\Big)\,e^{-\delta(x\cdot e-c_ft-B)}+\delta\,\Big(\frac{-f'(0)}{2}-\delta\Big)\,e^{-\delta(t-t_k)}\ge0,$$
again from~(\ref{defdelta}) and the negativity of $\phi'_f$.\par
Finally, if $(t,x)\in\mathcal{E}_k$ is such that $\overline{u}_k(t,x)<1$ and $-C\le\zeta_k(t,x)\le C$, then~(\ref{defkappa}) yields
$$-\phi'_f(\zeta_k(t,x))\ge\kappa>0,$$
while
$$f(\phi_f(\zeta_k(t,x)))-f(\overline{u}_k(t,x))\ge-L\,\big(\delta\,e^{-\delta(x\cdot e-c_ft-B)}+\delta\,e^{-\delta(t-t_k)}\big).$$
Therefore,
$$\baa{rcl}
N_k(t,x) & \ge & \delta\,(\kappa\,\omega-L-\delta)\,e^{-\delta(t-t_k)}+\delta\big(\kappa\,\omega\,c_f\,e^{\delta c_ft}-(L+\delta^2)\,e^{-\delta(x\cdot e-c_ft-B)}\big)\vspace{3pt}\\
& \ge & \delta\big(\kappa\,\omega\,c_f\,e^{\delta c_ft}-(L+\delta^2)\,e^{-\delta(x\cdot e-c_ft-B)}\big)\eaa$$
from~(\ref{defomega}). On the other hand, the inequality $\zeta_k(t,x)\le C$ implies that
$$\baa{rcl}
x\cdot e-c_ft-B & \ge & -2\,c_f\,t+c_ft_k+\xi^+_{t_k}+\omega\,e^{-\delta(t-t_k)}-\omega-\omega\,e^{\delta c_ft}-2B-2C\vspace{3pt}\\
& \ge & -2\,c_f\,t-2\,(\omega+B+C)\eaa$$
from~(\ref{xi+tk}) and the fact that $t\le T_2\le 0$. Thus,
$$-\delta\big(x\cdot e-c_ft-B)\le2\,\delta\,c_f\,t+2\,\delta\,(\omega+B+C)\le\delta\,c_f\,t$$
from~(\ref{defT2}), whence
$$N_k(t,x)\ge\delta\,(\kappa\,\omega\,c_f-L-\delta^2)\,e^{\delta c_ft}\ge 0$$
from~(\ref{defomega}).\par
As a conclusion, for every $k\in\N$, the function $\overline{u}_k$ is a supersolution, in the set~$\mathcal{E}_k$, of the equation~(\ref{eq}) satisfied by $u$. The maximum principle implies that
\be\label{ineq}
u(t,x)\le\overline{u}_k(t,x)\le\phi_f(\zeta_k(t,x))+\delta\,e^{-\delta(x\cdot e-c_ft-B)}+\delta\,e^{-\delta(t-t_k)}
\ee
for all $k\in\N$ and $(t,x)\in\mathcal{E}_k$. As a consequence, for any fixed $(t,x)\in\R\times\R^N$ such that $t\le T_2$ and $x\cdot e\ge c_f+B$, one has $(t,x)\in\mathcal{E}_k$ for $k$ large enough, while $t_k\to-\infty$ and $\zeta_k(t,x)\to+\infty$ as $k\to+\infty$ from our assumption~(\ref{assumption}). Passing to the limit as $k\to+\infty$ in~(\ref{ineq}) gives
$$u(t,x)\le\delta\,e^{-\delta(x\cdot e-c_ft-B)}$$
for all $t\le T_2$ and $x\cdot e\ge c_f+B$. The proof of Lemma~\ref{lemexp} is thereby complete.\hfill$\Box$


\SE{Characterization of the global mean speed of transition fronts}\label{sec4}

This section is devoted to the proof of Theorem~\ref{thcf} on the existence and the uniqueness of the global mean speed of the transition fronts connecting $0$ and $1$ for problem~(\ref{eq}). The gene\-ral idea can be summarized as follows. Any such transition front $u$ is above~$\theta^+$ (resp. below~$\theta^-$) on big sets in some neighborhoods of $\Gamma_t$, thanks to~(\ref{unifgamma}). We recall that~$0<\theta^-\le\theta^+<1$ are given in~(\ref{thetapm}). When $c_f>0$, a solution of the Cauchy problem associated to~(\ref{eq}) with a compactly supported initial condition above some $\beta>\theta^+$ on a large ball spreads at the speed~$c_f$ in all directions. On the other hand, when the initial condition is below some $\alpha<\theta^-$ on a large ball and is, say, equal to $1$ outside, then the region where the solution is above~$\alpha$ cannot move towards the center of the ball too much faster than~$c_f$. These two key-properties, which are stated in Lemmas~\ref{leminf} and~\ref{lemsup} below, will force the interfaces~$\Gamma_t$ of the transition front~$u$ to move at the global mean speed $\gamma=|c_f|$, in the sense of~(\ref{defmeanspeed}).


\subsection{Two key-lemmas}\label{seckeylemmas}

In the sequel, we fix to real numbers $\alpha$ and $\beta$ as in~(\ref{alphabeta0}), that is
\be\label{alphabeta}
0<\alpha<\theta^-\le\theta^+<\beta<1,
\ee
where $\theta^{\pm}$ are defined in~(\ref{thetapm}). For any $R>0$, let $v_R$ and $w_R$ denote the solutions of the Cauchy problems 
\be\label{vwR}\left\{\baa{rcll}
(v_R)_t & = & \Delta v_R+f(v_R), & t>0,\ x\in\R^N,\vspace{3pt}\\
(w_R)_t & = & \Delta w_R+f(w_R), & t>0,\ x\in\R^N,\eaa\right.
\ee
with initial conditions
\be\label{vR0}
v_R(0,x)=\left\{\baa{ll} \beta & \hbox{if }|x|<R,\vspace{3pt}\\ 0 & \hbox{if }|x|\ge R\eaa\right.
\ee
and
\be\label{wR0}
w_R(0,x)=\left\{\baa{ll} \alpha & \hbox{if }|x|<R,\vspace{3pt}\\ 1 & \hbox{if }|x|\ge R\eaa\right.
\ee

\begin{lem}\label{leminf}
Assume that $c_f>0$. There is $R>0$ such that the following holds: for all~$\epsilon\in(0,c_f]$, there is $T_{\epsilon}>0$ such that
\be\label{expansion}
v_R(t,x)\ge\beta\ \hbox{ for all }t\ge T_{\epsilon}\hbox{ and }|x|\le(c_f-\epsilon)t
\ee
and, in fact,
\be\label{convvR1}
v_R(t,\cdot)\to1\hbox{ uniformly in }\big\{x\in\R^N;\ |x|\le(c_f-\epsilon)t\big\}\hbox{ as }t\to+\infty.
\ee
\end{lem}

In the proof of Theorem~\ref{thcf}, Lemma~\ref{leminf} will provide a sharp lower bound for the speed of the interfaces $\Gamma_t$ of any transition front of~(\ref{eq}) connecting $0$ and $1$.\par
The following lemma, which is a kind of counterpart of Lemma~\ref{leminf}, will give the upper bound. That is, if~$c_f\ge 0$ and if the initial condition $w_R(0,\cdot)$ given in~(\ref{wR0}) is equal to (it could also be less than)~$\alpha$ on a very large ball with radius $R$ and equal to $1$ outside, then the region where $w_R$ is above~$\alpha$ is not filled at a speed too much larger than~$c_f$ on some interval of time whose size is related to the initial radius $R$.

\begin{lem}\label{lemsup}
Assume that $c_f\ge 0$. Then, for any $\epsilon>0$, there are some real numbers $T_{\epsilon}>0$ and $R_{\epsilon}\ge(c_f+\epsilon)T_{\epsilon}>0$ such that for all $R\ge R_{\epsilon}$, the solution $w_R$ of~$(\ref{vwR})$ and~$(\ref{wR0})$ satisfies
\be\label{wRbeta}
w_R(t,x)\le\alpha\ \hbox{ for all }T_{\epsilon}\le t\le\frac{R}{c_f+\epsilon}\hbox{ and }|x|\le R-(c_f+\epsilon)t.
\ee
\end{lem}

Before doing the proof of these two lemmas, let us first comment and compare them with some existing results of the literature. Lemma~\ref{leminf} could be viewed as an analog of the one-dimensional propagation result of Fife and McLeod~\cite{fm} used in~(\ref{convw}) in the proof of Theorem~\ref{thplanar}. Actually,~(\ref{convw}) is much more precise than Lemma~\ref{leminf} since~(\ref{convw}) implies in particular that the position of any given level set of the considered solution is $\pm c_ft+O(1)$ as $t\to+\infty$. Such an estimate cannot hold in higher dimension due to the curvature effects. Actually, the position of the level sets of the solutions $v_R$ of~(\ref{vwR}-\ref{vR0}) could likely be estimated more precisely, but the conclusion of Lemma~\ref{leminf} will be sufficient for the proof of Theorem~\ref{thcf}.\par
Notice also that if $f$ is of the bistable type~(\ref{bistable}) and if $\beta=\beta_f<1$ is sufficiently close to~$1$, then Lemma~\ref{leminf} follows directly from Theorem~6.2 of Aronson and Weinberger~\cite{aw}. The proof of Lemma~\ref{leminf} is inspired by~\cite{aw} but the subsolutions used in the proof of Lemma~\ref{leminf} are different and, as such, the result is new.\par
Lemmas~\ref{leminf} and~\ref{lemsup} also have some similarities with Theorem~$2$ of Chen~\cite{ch}. In~\cite{ch}, some problems with small diffusion have been considered, for bistable functions $f$ of the type~(\ref{bistable}) with $f'(\theta)>0$. In this case, after scaling and coming back to~(\ref{eq}), it follows from~\cite{ch} that, if $c_f>0$ and if $R>0$ is large, then $v_R(t,x)$ is larger than $1-O(R^{-k})$ (where~$k>0$ is given) for $t\ge O(\ln R)$ and $|x|\le R+(c_f-O(R^{-1}))t-O(1)$. This latter estimate is quantitatively more precise than~(\ref{expansion}) for $R\ge O(\epsilon^{-1})$. In Lemma~\ref{leminf} of the present paper,~$R$ can be chosen independently of $\epsilon$ and we only need $v_R$ to be not too far from $1$, without any precise rate of convergence, in some balls expanding with a speed close to $c_f$. The proof of~\cite{ch} is based on the construction of subsolutions with nonlinearities of the type $f-\lambda$. Even if the ideas of~\cite{ch} could likely be adapted here to the case of functions~$f$ satisfying only~(\ref{f}) together with the existence of a planar front~$(c_f,\phi_f)$, we are going to work directly with the function~$f$ in the proof of Lemma~\ref{leminf}. As far as Lemma~\ref{lemsup} is concerned, it also follows from~\cite{ch} that, if~$c_f>0$ with~(\ref{bistable}) and~$f'(\theta)>0$, and if $R>0$ is large, then $w_R(t,x)$ is smaller than $O(R^{-k})$ for~$O(\ln R)\le t\le R/(2c_f)$ and~$|x|\le R-c_ft-O(\ln R)$. The above pointwise estimate~(\ref{wRbeta}), that is $w_R(t,x)\le\alpha$, is less precise. However it holds until a time of the order~$R/c_f$ (instead of $R/(2c_f)$). In the proof of Theorem~\ref{thcf}, we will actually need the estimate~(\ref{wRbeta}) on a time interval of the order~$R/c_f$, in order to show that the global mean speed of any transition front for problem~(\ref{eq}) exists and is exactly equal to $c_f$.\par
Finally, even if Lemmas~\ref{leminf} and~\ref{lemsup} have some similarities with~\cite{aw,ch}, the assumptions and conclusions are different and, as such, the statements and the proofs are new to the best of our knowledge. Let us now turn to the proofs.\hfill\break

\noindent{\bf{Proof of Lemma~\ref{leminf}.}} We assume here that $c_f>0$. Observe first that~(\ref{convvR1}) is clearly stronger than~(\ref{expansion}). The strategy to prove~(\ref{convvR1}) is to construct some radially symmetric subsolutions in~$\R^N$ of the type $\phi_f(|x|-(c_f-\epsilon/2)t-R)$ plus some exponentially decaying terms, for~$0<\epsilon\le c_f$. In doing so, the solution $v_R$ will be close to~$1$ inside the balls of radii~$(c_f-\epsilon)t$ as $t$ is large.\par
{\it Step 1: choice of some parameters which are independent of~$\epsilon$.} As in the proof of Lemma~\ref{lemexp}, we first introduce some parameters which are independent of~$\epsilon$. Let $\delta>0$ be chosen such that
\be\label{defdelta2}
0<\delta<\min\Big(1,\frac{|f'(0)|}{2},\frac{|f'(1)|}{2}\Big),\ \ f'\le\frac{f'(0)}{2}\hbox{ on }[0,\delta]\ \hbox{ and }\ f'\le\frac{f'(1)}{2}\hbox{ on }[1-2\delta,1].
\ee
Since $\phi''_f(s)\sim-\sigma\,e^{\mu s}$ as $s\to-\infty$ with $\sigma>0$ and $\mu=(-c_f+(c_f^2-4f'(1))^{1/2})/2>0$, one can choose $C>0$ so that~(\ref{defC}) holds together with $\phi''_f\le0$ on $(-\infty,-C]$, that is
\be\label{defC2}
\phi_f\ge1-\delta\hbox{ on }(-\infty,-C],\ \ \phi''_f\le0\hbox{ on }(-\infty,-C]\ \hbox{ and }\ \phi_f\le\delta\hbox{ on }[C,+\infty).
\ee
Let $\kappa>0$ be chosen as in~(\ref{defkappa}), that is $-\phi'_f\ge\kappa$ on $[-C,C]$, and let $\omega>0$ be such that
\be\label{defomega2}
\kappa\,\omega\ge L+\delta,
\ee
where $L=\max_{[0,1]}|f'|$.\par
Let $\varrho_{\beta}$ be the solution of the ordinary differential equation~$\varrho_{\beta}'(t)=f(\varrho_{\beta}(t))$ with initial condition $\varrho_{\beta}(0)=\beta$. Since $\beta\in(\theta^+,1)$, there holds $\varrho_{\beta}(t)\to1$ as~$t\to+\infty$, and there is $T>0$ such that $\varrho_{\beta}(T)\ge1-\delta/2$. It follows from the maximum principle that
$$0\le \varrho_{\beta}(T)-v_R(T,x)\le\frac{e^{LT}}{(4\pi T)^{N/2}}\int_{|y|\ge R}e^{-\frac{|x-y|^2}{4T}}dy$$
for all $R>0$ and $x\in\R^N$. Therefore, if $0<B\le R$ and $|x|\le R-B$, one infers that
$$0\le \varrho_{\beta}(T)-v_R(T,x)\le\frac{e^{LT}}{(4\pi T)^{N/2}}\int_{|z|\ge B}e^{-\frac{|z|^2}{4T}}dz.$$
Thus, there exists $B>0$ such that, for all $R\ge B$ and $|x|\le R-B$, $\varrho_{\beta}(T)-v_R(T,x)\le\delta/2$, whence
\be\label{defB1}
v_R(T,x)\ge \varrho_{\beta}(T)-\frac{\delta}{2}\ge1-\delta\hbox{ for all }R\ge B\hbox{ and }|x|\le R-B.
\ee\par
{\it Step 2: choice of some functions $h_{\epsilon}$ which depend on $\epsilon$}. It is elementary to check that, for every $\epsilon>0$, there is a $C^2$ function $h_{\epsilon}:[0,+\infty)\to\R$ satisfying the following properties:
\be\label{defh0}\left\{\baa{l}
0\le h'_{\epsilon}\le 1\hbox{ on }[0,+\infty),\vspace{3pt}\\
h'_{\epsilon}=0\hbox{ on a neighborhood of }0,\vspace{3pt}\\
h_{\epsilon}(r)=r\hbox{ on }[H_{\epsilon},+\infty)\hbox{ for some }H_{\epsilon}>0,\vspace{3pt}\\
\displaystyle{\frac{(N-1)h'_{\epsilon}(r)}{r}}+h''_{\epsilon}(r)\le\displaystyle{\frac{\epsilon}{2}}\hbox{ on }[0,+\infty).\eaa\right.
\ee
Notice in particular that, necessarily,
\be\label{defh1}
r\le h_{\epsilon}(r)\le r+h_{\epsilon}(0)\hbox{ for all }r\ge 0.
\ee\par
{\it Step 3: proof of~$(\ref{convvR1})$ when $c_f/2\le\epsilon\le c_f$.} To do so, it is sufficient to show that~(\ref{convvR1}) holds with~$\epsilon=\epsilon_0:=c_f/2>0$, for some $R>0$. Let us set
\be\label{deftildeR}
R=B+H_{\epsilon_0}+\omega+2\,C>B>0.
\ee
We will see that the conclusion of Lemma~\ref{leminf} holds for all $0<\epsilon\le c_f$ with this value of $R$ (notice that $R$ is independent of $\epsilon$). Let us show in this step that~(\ref{convvR1}) holds with~$\epsilon=\epsilon_0$. For all $(t,x)\in\R\times\R^N$, we set
\be\label{defunderv}
\underline{v}(t,x)=\max\Big(\phi_f\big(\underline{\zeta}(t,x)\big)-\delta\,e^{-\delta(t-T)},0\Big),
\ee
where
\be\label{defunderzeta}
\underline{\zeta}(t,x)=h_{\epsilon_0}(|x|)-\Big(c_f-\frac{\epsilon_0}{2}\Big)\,(t-T)-\omega\,e^{-\delta(t-T)}-H_{\epsilon_0}-C.
\ee
Let us then check that $\underline{v}$ is a subsolution for the problem satisfied by $v_R$, for $t\ge T$ and $x\in\R^N$.\par
First, at the time $T$, it follows from~(\ref{defB1}),~(\ref{deftildeR}) and the definition of $\underline{v}$ that
$$v_R(T,x)\ge1-\delta\ge\underline{v}(T,x)\hbox{ for all }|x|\le{R}-B.$$
On the other hand, if $|x|\ge{R}-B$, then $h_{\epsilon_0}(|x|)\ge|x|\ge{R}-B=H_{\epsilon_0}+\omega+2C$ from~(\ref{defh1}) and~(\ref{deftildeR}), whence $\underline{\zeta}(T,x)\ge C$ and
$$\underline{v}(T,x)=\max\big(\phi_f(\underline{\zeta}(T,x))-\delta,0\big)\le\max\big(\delta-\delta,0\big)=0\le v_{{R}}(T,x)$$
from~(\ref{defC2}) and the fact that $v_R\ge0$ in $(0,+\infty)\times\R^N$. Thus,
$$v_{{R}}(T,x)\ge\underline{v}(T,x)\hbox{ for all }x\in\R^N.$$\par
Let us now check that
\be\label{defunderN}
\underline{N}(t,x)=\underline{v}_t(t,x)-\Delta\underline{v}(t,x)-f(\underline{v}(t,x))\le 0
\ee
for all $t\ge T$ and $x\in\R^N$ such that $\underline{v}(t,x)>0$. This will be sufficient to ensure that $\underline{v}$ is a subsolution, since $f(0)=0$ (notice that $\underline{v}(t,x)=\phi_f(\underline{\zeta}(t,x))-\delta e^{-\delta(t-T)}$ is of class $C^2$ in the set where it is positive, since $\phi_f$ is of class $C^2$ and $h$ vanishes in a neighbourhood of $0$).\par
In this paragraph, let $(t,x)$ be any point in $[T,+\infty)\times\R^N$ such that $\underline{v}(t,x)>0$. Since $\phi_f$ obeys~(\ref{eqphi}), there holds
$$\baa{rcl}
\underline{N}(t,x) & \!=\! & f(\phi_f(\underline{\zeta}(t,x)))-f(\underline{v}(t,x))+\omega\,\delta\,e^{-\delta(t-T)}\,\phi'_f(\underline{\zeta}(t,x))+\delta^2\,e^{-\delta(t-T)}\vspace{3pt}\\
& \!\! & +\Big(\displaystyle{\frac{\epsilon_0}{2}}-\displaystyle{\frac{(N\!-\!1)\,h'_{\epsilon_0}(|x|)}{|x|}}-h''_{\epsilon_0}(|x|)\Big)\,\phi_f'(\underline{\zeta}(t,x))+\big(1-(h'_{\epsilon_0}(|x|))^2\big)\,\phi''_f(\underline{\zeta}(t,x))\vspace{3pt}\\
& \!\le\! & f(\phi_f(\underline{\zeta}(t,x)))-f(\underline{v}(t,x))+\omega\,\delta\,e^{-\delta(t-T)}\,\phi'_f(\underline{\zeta}(t,x))+\delta^2\,e^{-\delta(t-T)}\vspace{3pt}\\
& \!\! & +\big(1-(h'_{\epsilon_0}(|x|))^2\big)\,\phi''_f(\underline{\zeta}(t,x))\eaa$$
from~(\ref{defh0}) and $\phi'_f\le 0$. If $\underline{\zeta}(t,x)\le-C$, then $1-\delta\le\phi_f(\underline{\zeta}(t,x))<1$ from~(\ref{defC2}), whence $1-2\delta\le\underline{v}(t,x)\le\phi_f(\underline{\zeta}(t,x))<1$ and
$$f(\phi_f(\underline{\zeta}(t,x)))-f(\underline{v}(t,x))\le\frac{f'(1)}{2}\,\delta\,e^{-\delta(t-T)}$$
from~(\ref{defdelta2}). Furthermore, $\phi_f''(\underline{\zeta}(t,x))\le0$ from~(\ref{defC2}), while $0\le h'_{\epsilon_0}(|x|)\le 1$ from~(\ref{defh0}). Therefore, if $\underline{\zeta}(t,x)\le-C$, then
$$\underline{N}(t,x)\le\delta\Big(\frac{f'(1)}{2}+\delta\Big)\,e^{-\delta(t-T)}+\omega\,\delta\,e^{-\delta(t-T)}\,\phi'_f(\underline{\zeta}(t,x))\le0$$
from~(\ref{defdelta2}) and the negativity of $\phi'_f$. On the other hand, if $\underline{\zeta}(t,x)\ge C$, then $\phi_f(\underline{\zeta}(t,x))\le\delta$, whence $0<\underline{v}(t,x)\le\phi_f(\underline{\zeta}(t,x))\le\delta$ and
$$f(\phi_f(\underline{\zeta}(t,x)))-f(\underline{v}(t,x))\le\frac{f'(0)}{2}\,\delta\,e^{-\delta(t-T)},$$
again from~(\ref{defdelta2}). The inequality~$\underline{\zeta}(t,x)\ge C$ also implies that $h_{\epsilon_0}(|x|)\ge 2\,C+H_{\epsilon_0}\ge H_{\epsilon_0}$, whence $h'_{\epsilon_0}(|x|)=1$ from~(\ref{defh0}). Thus, if $\underline{\zeta}(t,x)\ge C$, then
$$\underline{N}(t,x)\le\delta\Big(\frac{f'(0)}{2}+\delta\Big)\,e^{-\delta(t-T)}+\omega\,\delta\,e^{-\delta(t-T)}\,\phi'_f(\underline{\zeta}(t,x))\le0$$
from~(\ref{defdelta2}) and the negativity of $\phi'_f$. Finally, if $-C\le\underline{\zeta}(t,x)\le C$, then 
$$f(\phi_f(\underline{\zeta}(t,x)))-f(\underline{v}(t,x))\le L\,\delta\,e^{-\delta(t-T)},$$
where we recall that $L=\max_{[0,1]}|f'|$, while $\underline{\zeta}(t,x)\ge-C$ yields $h_{\epsilon_0}(|x|)\ge H_{\epsilon_0}$ and $h'_{\epsilon_0}(|x|)=1$. Hence, if $-C\le\underline{\zeta}(t,x)\le C$, then 
$$\underline{N}(t,x)\le\delta\,(L-\kappa\,\omega+\delta)\,e^{-\delta(t-T)}\le0$$
from~(\ref{defkappa}) and~(\ref{defomega2}).\par
As a conclusion, the maximum principle implies that, for all $t\ge T$ and~$x\in\R^N$,
\be\label{vtildeR}
1\ge v_{{R}}(t,x)\ge\underline{v}(t,x)\ge\phi_f\big(\underline{\zeta}(t,x)\big)-\delta\,e^{-\delta(t-T)}.
\ee
But
$$\max_{|x|\le(c_f-\epsilon_0)t}\underline{\zeta}(t,x)\le(c_f-\epsilon_0)\,t+h_{\epsilon_0}(0)-\Big(c_f-\frac{\epsilon_0}{2}\Big)(t-T)\to-\infty\hbox{ as }t\to+\infty,$$
from~(\ref{defh1}),~(\ref{defunderzeta}) and the positivity of $\epsilon_0$, $\omega$, $H_{\epsilon_0}$ and $C$. Since $\phi_f(-\infty)=1$, it follows from~(\ref{vtildeR}) that
\be\label{convvR1bis}
v_R(t,\cdot)\to1\hbox{ uniformly in }\big\{x\in\R^N;\ |x|\le(c_f-\epsilon_0)t\big\}\hbox{ as }t\to+\infty.
\ee\par
{\it Step 4: conclusion and proof of~$(\ref{convvR1})$ for all $0<\epsilon\le c_f$}. Property~(\ref{convvR1}) has already been proved for $\epsilon_0=c_f/2\le\epsilon\le c_f$, from~(\ref{convvR1bis}). Let now $\epsilon$ be arbitrary in $(0,\epsilon_0)$. With the notations used in Steps~1 and~2, set
\be\label{defunderR}
\underline{R}_{\epsilon}=H_{\epsilon}+\omega+2\,C>0
\ee
and, from~(\ref{convvR1bis}), let $\underline{T}_{\epsilon}\ge T$ such that
$$v_R(\underline{T}_{\epsilon},x)\ge 1-\delta\hbox{ for all }|x|\le\underline{R}_{\epsilon}.$$\par
Let $\underline{v}$ and $\underline{\zeta}$ be defined as in~(\ref{defunderv}) and~(\ref{defunderzeta}) where $T$ and $\epsilon_0$ are replaced by $\underline{T}_{\epsilon}$ and $\epsilon$. The same calculations as in Step~3 show that~(\ref{defunderN}) holds for all $(t,x)\in[\underline{T}_{\epsilon},+\infty)\times\R^N$ such that~$\underline{v}(t,x)>0$. The only difference now is the comparison of $v_R$ and $\underline{v}$ at time $\underline{T}_{\epsilon}$. If $|x|\le\underline{R}_{\epsilon}$, then $v_R(\underline{T}_{\epsilon},x)\ge1-\delta\ge\underline{v}(\underline{T}_{\epsilon},x)$. If $|x|\ge\underline{R}_{\epsilon}$, then
$$\underline{\zeta}(\underline{T}_{\epsilon},x)=h_{\epsilon}(|x|)-\omega-H_{\epsilon}-C\ge\underline{R}_{\epsilon}-\omega-H_{\epsilon}-C=C$$
from~(\ref{defh1}) and~(\ref{defunderR}), whence $\phi_f(\underline{\zeta}(\underline{T}_{\epsilon},x))\le\delta$ from~(\ref{defC2}) and $\underline{v}(\underline{T}_{\epsilon},x)=0\le v_R(\underline{T}_{\epsilon},x)$. Finally,
$$v_R(\underline{T}_{\epsilon},x)\ge\underline{v}(\underline{T}_{\epsilon},x)\hbox{ for all }x\in\R^N.$$
Therefore, it follows from the maximum principle that
$$v_R(t,x)\ge\underline{v}(t,x)\ge\phi_f(\underline{\zeta}(t,x))-\delta\,e^{-\delta(t-\underline{T}_{\epsilon})}\hbox{ for all }t\ge\underline{T}_{\epsilon}\hbox{ and }x\in\R^N.$$
As in Step~3, this leads to~(\ref{convvR1}). The proof of Lemma~\ref{leminf} is thereby complete.\hfill$\Box$\break

Let us now turn to the proof of Lemma~\ref{lemsup}. The strategy used is similar to that used in the proof of Lemma~\ref{leminf}, with the construction of suitable supersolutions instead of subsolutions. However, one needs to be more careful with the application of the maximum principle, since the estimates~(\ref{wRbeta}) only hold in bounded time intervals. Notice also that Lemma~\ref{lemsup} is not an immediate consequence of Lemma~\ref{leminf} after changing $v$ into $1-v$ and $f(s)$ into $-f(1-s)$, since this operation would change the sign of the speed $c_f$ (however, Corollary~\ref{corsup} below can be deduced from Lemma~\ref{lemsup} thanks to this operation).\hfill\break

\noindent{\bf{Proof of Lemma~\ref{lemsup}.}} The strategy is to construct some radially symmetric supersolutions in $\R^N$ of the type $\phi_f(-(c_f+\epsilon/2)t+R-|x|)$ plus some exponentially decaying terms. In doing so, the solution $w_R$ will be small inside the balls of radii $R-(c_f+\epsilon)t$.\par
As in the proofs of Lemmas~\ref{lemexp} and~\ref{leminf}, we first introduce some parameters which are independent of~$\epsilon$. Let $\delta>0$ be chosen so that
\be\label{defdelta3}
0<\delta<\min\Big(1,\frac{|f'(0)|}{2},\frac{|f'(1)|}{2}\Big),\ \ f'\le\frac{f'(0)}{2}\hbox{ on }[0,2\delta]\ \hbox{ and }\ f'\le\frac{f'(1)}{2}\hbox{ on }[1-\delta,1].
\ee
Since $\phi''_f(s)\sim\nu\,e^{-\lambda s}$ as $s\to+\infty$ with $\nu>0$ and $\lambda=\big(c_f+(c_f^2-4f'(0))^{1/2}\big)/2$, one can choose $C>0$ so that~(\ref{defC}) holds together with $\phi''_f\ge0$ and $\phi_f\le\alpha/2$ on $[C,+\infty)$, that is
\be\label{defCbis}
\phi_f\ge1-\delta\hbox{ on }(-\infty,-C],\ \ \phi_f\le\min\Big(\delta,\frac{\alpha}{2}\Big)\hbox{ on }[C,+\infty)\ \hbox{ and }\ \phi''_f\ge0\hbox{ on }[C,+\infty).
\ee
Let $\kappa>0$ and $\omega>0$ be chosen as in~(\ref{defkappa}) and~(\ref{defomega2}), that is
\be\label{kappaomega}
-\phi_f'\ge\kappa>0\hbox{ on }[-C,C]\ \hbox{ and }\ \kappa\,\omega\ge L+\delta.
\ee
Let $\varrho_{\alpha}$ be the solution of the ordinary differential equation~$\varrho_{\alpha}'(t)=f(\varrho_{\alpha}(t))$ with initial condition $\varrho_{\alpha}(0)=\alpha$. Since $\alpha\in(0,\theta^-)$, there holds $\varrho_{\alpha}(t)\to0$ as~$t\to+\infty$, and there is~$T>0$ such that $\varrho_{\alpha}(T)\le\delta/2$. As in Step~1 of the proof of Lemma~\ref{leminf}, it follows from the maximum principle that there exists $B>0$ such that $0\le w_R(T,x)-\varrho_{\alpha}(T)\le\delta/2$ for all $R\ge B$ and~$|x|\le R-B$, whence
\be\label{defB}
w_R(T,x)\le\delta\hbox{ for all }R\ge B\hbox{ and }|x|\le R-B.
\ee\par
Now, we pick an arbitrary $\epsilon>0$ and we introduce some quantities which depend on $\epsilon$. Let the function $h_{\epsilon}$ be as in~(\ref{defh0}) and~(\ref{defh1}). Then, we choose $T_{\epsilon}\ge T\,(>0)$ such that
\be\label{defT}
\delta\,e^{-\delta(t-T)}\le\frac{\alpha}{2}\ \hbox{ and }\ \frac{\epsilon\,t}{2}\ge h_{\epsilon}(0)+\omega+B+2C\ \hbox{ for all }t\ge T_{\epsilon},
\ee
and $R_{\epsilon}>0$ such that
\be\label{defR}
R_{\epsilon}\ge\max\big(B,(c_f+\epsilon)T_{\epsilon}\big)\hbox{ and }\frac{\epsilon\,R_{\epsilon}}{2(c_f+\epsilon)}\ge\omega+B+2C+H_{\epsilon}.
\ee
We shall now prove that the conclusion of Lemma~\ref{lemsup} holds with these choices of $T_{\epsilon}$ and $R_{\epsilon}$.\par
In the sequel, $R$ is an arbitrary real number such that
\be\label{RRepsilon}
R\ge R_{\epsilon}.
\ee
For all $(t,x)\in\R\times\R^N$, we set
$$\overline{w}(t,x)=\min\Big(\phi_f\big(\overline{\zeta}(t,x)\big)+\delta\,e^{-\delta(t-T)},1\Big),$$
where
$$\overline{\zeta}(t,x)=-h_{\epsilon}(|x|)-\Big(c_f+\frac{\epsilon}{2}\Big)\,(t-T)+\omega\,e^{-\delta(t-T)}-\omega+R-B-C.$$
Let us then check that $\overline{w}$ is a supersolution for the problem satisfied by $w_R$, in the set
$$\mathcal{E}=\Big[T,\frac{R}{c_f+\epsilon}\Big]\times\R^N.$$
At the time $T$, it follows from~(\ref{defB}),~(\ref{defR}),~(\ref{RRepsilon}) and the definition of $\overline{w}$ that
$$w_R(T,x)\le\delta\le\overline{w}(T,x)\hbox{ for all }|x|\le R-B.$$
On the other hand, if $|x|\ge R-B$, then $h_{\epsilon}(|x|)\ge|x|\ge R-B$ from~(\ref{defh1}), whence $\overline{\zeta}(T,x)\le-C$ and
$$\overline{w}(T,x)=\min\big(\phi_f(\overline{\zeta}(T,x))+\delta,1\big)\ge\min\big((1-\delta)+\delta,1\big)=1\ge w_R(T,x)$$
from~(\ref{defCbis}) and the fact that $w_R\le 1$ in $(0,+\infty)\times\R^N$. Thus,
$$w_R(T,x)\le\overline{w}(T,x)\hbox{ for all }x\in\R^N.$$\par
Let us now check that
$$\overline{N}(t,x)=\overline{w}_t(t,x)-\Delta\overline{w}(t,x)-f(\overline{w}(t,x))\ge 0\hbox{ for all }(t,x)\in\mathcal{E}\hbox{ such that }\overline{w}(t,x)<1.$$
This will be sufficient to ensure that $\overline{w}$ is a supersolution, since $f(1)=0$ (notice that, as for $\underline{v}$ in Lemma~\ref{leminf},~$\overline{w}$ is of class $C^2$ in the set where it is less than $1$).\par
In this paragraph, let $(t,x)$ be any point in~$\mathcal{E}$ such that $\overline{w}(t,x)<1$. Since $\phi_f$ obeys~(\ref{eqphi}), there holds
$$\baa{rcl}
\overline{N}(t,x) & = & f(\phi_f(\overline{\zeta}(t,x)))-f(\overline{w}(t,x))-\omega\,\delta\,e^{-\delta(t-T)}\,\phi'_f(\overline{\zeta}(t,x))-\delta^2\,e^{-\delta(t-T)}\vspace{3pt}\\
& & -\Big(\displaystyle{\frac{\epsilon}{2}}-\displaystyle{\frac{(N-1)\,h'_{\epsilon}(|x|)}{|x|}}-h''_{\epsilon}(|x|)\Big)\,\phi_f'(\overline{\zeta}(t,x))+\big(1-(h'_{\epsilon}(|x|))^2\big)\,\phi''_f(\overline{\zeta}(t,x))\vspace{3pt}\\
& \ge & f(\phi_f(\overline{\zeta}(t,x)))-f(\overline{w}(t,x))-\omega\,\delta\,e^{-\delta(t-T)}\,\phi'_f(\overline{\zeta}(t,x))-\delta^2\,e^{-\delta(t-T)}\vspace{3pt}\\
& & +\big(1-(h'_{\epsilon}(|x|))^2\big)\,\phi''_f(\overline{\zeta}(t,x))\eaa$$
from~(\ref{defh0}) and $\phi'_f\le 0$. If $\overline{\zeta}(t,x)\le-C$, then $1-\delta\le\phi_f(\overline{\zeta}(t,x))\le\overline{w}(t,x)<1$ from~(\ref{defCbis}), whence
$$f(\phi_f(\overline{\zeta}(t,x)))-f(\overline{w}(t,x))\ge\frac{-f'(1)}{2}\,\delta\,e^{-\delta(t-T)}$$
from~(\ref{defdelta3}). Furthermore, the inequalities $\overline{\zeta}(t,x)\le-C$ and $0<T\le t\le R/(c_f+\epsilon)$ yield
$$h_{\epsilon}(|x|)\ge-\Big(c_f+\frac{\epsilon}{2}\Big)\,(t-T)-\omega+R-B\ge\frac{\epsilon R}{2(c_f+\epsilon)}-\omega-B\ge H_{\epsilon}$$
because of~(\ref{defR}) and~(\ref{RRepsilon}) (notice that the term $-2C$ in~(\ref{defR}) is not needed here, but it will be later). The inequality $h_{\epsilon}(|x|)\ge H_{\epsilon}$ implies that $h'_{\epsilon}(|x|)=1$, due to the properties~(\ref{defh0}). Therefore, if $\overline{\zeta}(t,x)\le-C$, then
$$\overline{N}(t,x)\ge\delta\Big(\frac{-f'(1)}{2}-\delta\Big)\,e^{-\delta(t-T)}-\omega\,\delta\,e^{-\delta(t-T)}\,\phi'_f(\overline{\zeta}(t,x))\ge0$$
from~(\ref{defdelta3}) and the negativity of $\phi'_f$. On the other hand, if $\overline{\zeta}(t,x)\ge C$, then $\phi_f(\overline{\zeta}(t,x))\le\delta$ from~(\ref{defCbis}), whence $0<\phi_f(\overline{\zeta}(t,x))\le\overline{w}(t,x)\le2\delta$ and
$$f(\phi_f(\overline{\zeta}(t,x)))-f(\overline{w}(t,x))\ge\frac{-f'(0)}{2}\,\delta\,e^{-\delta(t-T)}$$
from~(\ref{defdelta3}). Since $\phi_f''\ge0$ on $[C,+\infty)$ from~(\ref{defCbis}), since $0\le h_{\epsilon}'\le 1$ on $[0,+\infty)$ and since~$\phi'_f\le 0$ on $\R$, one gets that, if $\overline{\zeta}(t,x)\ge C$, then
$$\overline{N}(t,x)\ge\delta\Big(\frac{-f'(0)}{2}-\delta\Big)\,e^{-\delta(t-T)}-\omega\,\delta\,e^{-\delta(t-T)}\,\phi'_f(\overline{\zeta}(t,x))\ge0$$
from~(\ref{defdelta3}). Lastly, if $-C\le\overline{\zeta}(t,x)\le C$, then 
$$f(\phi_f(\overline{\zeta}(t,x)))-f(\overline{w}(t,x))\ge-L\,\delta\,e^{-\delta(t-T)},$$
while $\overline{\zeta}(t,x)\le C$ yields
$$h_{\epsilon}(|x|)\ge-\Big(c_f+\frac{\epsilon}{2}\Big)\,(t-T)-\omega+R-B-2C\ge\frac{\epsilon R}{2(c_f+\epsilon)}-\omega-B-2C\ge H_{\epsilon}$$
from~(\ref{defR}) and (\ref{RRepsilon}). Thus, $h'_{\epsilon}(|x|)=1$ and
$$\overline{N}(t,x)\ge\delta\,(-L+\kappa\,\omega-\delta)\,e^{-\delta(t-T)}\ge0$$
from~(\ref{kappaomega}).\par
As a conclusion, the maximum principle implies that, for all $T\le t\le R/(c_f+\epsilon)$ and~$x\in\R^N$,
$$w_R(t,x)\le\overline{w}(t,x)\le\phi_f\big(\overline{\zeta}(t,x)\big)+\delta\,e^{-\delta(t-T)}.$$
For all $T_{\epsilon}\le t\le R/(c_f+\epsilon)$ and $|x|\le R-(c_f+\epsilon)t$, there holds $\delta\,e^{-\delta(t-T)}\le\alpha/2$ from~(\ref{defT}), while $h_{\epsilon}(|x|)\le|x|+h_{\epsilon}(0)\le R-(c_f+\epsilon)t+h_{\epsilon}(0)$ and
$$\baa{rcl}
\overline{\zeta}(t,x) & \ge & -R+(c_f+\epsilon)t-h_{\epsilon}(0)-\Big(c_f+\displaystyle{\frac{\epsilon}{2}}\Big)\,(t-T)+\omega\,e^{-\delta(t-T)}-\omega+R-B-C\vspace{3pt}\\
& \ge & \displaystyle{\frac{\epsilon\,t}{2}}-h_{\epsilon}(0)-\omega-B-C\vspace{3pt}\\
& \ge & C\eaa$$
from~(\ref{defT}). Thus, $\phi_f(\overline{\zeta}(t,x))\le\alpha/2$ from~(\ref{defCbis}). Finally, if $T_{\epsilon}\le t\le R/(c_f+\epsilon)$ and~$|x|\le R-(c_f+\epsilon)t$, then
$$w_R(t,x)\le\phi_f\big(\overline{\zeta}(t,x)\big)+\delta\,e^{-\delta(t-T)}\le\frac{\alpha}{2}+\frac{\alpha}{2}=\alpha.$$
The proof of Lemma~\ref{lemsup} is thereby complete.\hfill$\Box$\break

By changing $f(s)$ into $-f(1-s)$, $c_f$ into $-c_f$, $w_R$ into $1-v_R$, $\alpha$ into $1-\beta$ and $\beta$ into~$1-\alpha$, the following result holds:

\begin{cor}\label{corsup}
Assume that $c_f\le 0$. Then, for any $\epsilon>0$, there are some real numbers $\widetilde{T}_{\epsilon}>0$ and $\widetilde{R}_{\epsilon}\ge(|c_f|+\epsilon)\widetilde{T}_{\epsilon}>0$ such that for all $R\ge\widetilde{R}_{\epsilon}$, the solution $v_R$ of~$(\ref{vwR})$ and~$(\ref{vR0})$ satisfies
\be\label{vRbeta}
v_R(t,x)\ge\beta\ \hbox{ for all }\widetilde{T}_{\epsilon}\le t\le\frac{R}{|c_f|+\epsilon}\hbox{ and }|x|\le R-(|c_f|+\epsilon)t.
\ee
\end{cor}


\subsection{Proof of Theorem~\ref{thcf}}

Let $u$ be any transition front connecting $0$ and $1$ for problem~(\ref{eq}). As in the proof of Theorem~\ref{thplanar}, even if it means changing~$u(t,x)$ into $\widetilde{u}(t,x)=1-u(t,x)$, $f(s)$ into $g(s)=-f(1-s)$ and $c_f$ into $-c_f$, one can assume without loss of generality that
$$c_f\ge0.$$\par
In order to prove that
$$\frac{d(\Gamma_t,\Gamma_s)}{|t-s|}\to c_f\hbox{ as }|t-s|\to+\infty,$$
we prove one inequality for the $\liminf$ and another one for the $\limsup$.\hfill\break\par
{\it Step 1: the lower estimate}. We first claim that
\be\label{claiminf}
\liminf_{|t-s|\to+\infty}\frac{d(\Gamma_t,\Gamma_s)}{|t-s|}\ge c_f.
\ee
Since there is nothing to prove when $c_f=0$, we only consider the case $c_f>0$. Let $\alpha$ and $\beta$ be given as in~(\ref{alphabeta}). From Definition~\ref{def1}, there is $M\ge 0$ such that
\be\label{defMter}\left\{\baa{l}
\forall\, t\in\R,\ \forall\,x\in\overline{\Omega^+_t},\quad\big(d(x,\Gamma_t)\ge M\big)\Longrightarrow\big(\beta\le u(t,x)<1\big),\vspace{3pt}\\
\forall\, t\in\R,\ \forall\,x\in\overline{\Omega^-_t},\quad\big(d(x,\Gamma_t)\ge M\big)\Longrightarrow\big(0<u(t,x)\le\alpha\big).\eaa\right.
\ee
Let $R>0$ be as in Lemma~\ref{leminf}. Since the functions $v_\rho$ are nondecreasing with respect to the parameter $\rho>0$, one can assume without loss of generality that $R\ge M$. From~(\ref{unifgamma}), there is~$r>0$ such that
\be\label{unifgammabis}
\forall\,t\in\R,\ \forall\,x\in\Gamma_t,\ \exists\,y^{\pm}_t\in\Omega^{\pm}_t,\ |x-y^{\pm}|\le r\hbox{ and }d(y^{\pm},\Gamma_t)\ge 2R.
\ee\par
Let now $\epsilon\in(0,c_f]$ and let $T_{\epsilon}>0$ be as in Lemma~\ref{leminf}. Let us assume by contradiction that
\be\label{absurde}
\liminf_{|t-s|\to+\infty}\frac{d(\Gamma_t,\Gamma_s)}{|t-s|}<c_f-2\,\epsilon.
\ee
There are then two sequences $(t_k)_{k\in\N}$ and $(s_k)_{k\in\N}$ of real numbers such that $|t_k-s_k|\to+\infty$ as $k\to+\infty$ and
$$d(\Gamma_{t_k},\Gamma_{s_k})<(c_f-2\,\epsilon)\,|t_k-s_k|\ \hbox{ for all }k\in\N.$$
Without loss of generality, one can assume that $t_k<s_k$ for all $k\in\N$. By definition of the distance $d(\Gamma_{t_k},\Gamma_{s_k})$, there are then two sequences $(x_k)_{k\in\N}$ and $(z_k)_{k\in\N}$ in $\R^N$ such that
$$x_k\in\Gamma_{t_k},\ \ z_k\in\Gamma_{s_k}\ \hbox{ and }\ |x_k-z_k|<(c_f-2\epsilon)\,(s_k-t_k)\ \hbox{ for all }k\in\N.$$\par
On the one hand, it follows from~(\ref{unifgammabis}) that there exists a sequence $(y^+_k)_{k\in\N}$ of points in~$\R^N$ such that
$$y^+_k\in\Omega^+_{t_k},\ \ |x_k-y^+_k|\le r\ \hbox{ and }\ d(y^+_k,\Gamma_{t_k})\ge 2R\ \hbox{ for all }k\in\N.$$
Property~(\ref{defMter}) implies then that, for every $k\in\N$ and every $y\in{B(y^+_k,R)}$, one has $y\in\Omega^+_{t_k}$ and $d(y,\Gamma_{t_k})\ge R\ge M$, whence $u(t_k,y)\ge\beta$. Therefore, $u(t_k,x)\ge v_R(0,x-y^+_k)$ and
$$u(t,x)\ge v_R(t-t_k,x-y^+_k)\ \hbox{ for all }k\in\N,\ t>t_k\hbox{ and }x\in\R^N$$
from the maximum principle. The conclusion of Lemma~\ref{leminf} implies that, for every $k\in\N$,
\be\label{betamin}
u(t,x)\ge\beta\ \hbox{ for all }t\ge t_k+T_{\epsilon}\hbox{ and }|x-y^+_k|\le(c_f-\epsilon)\,(t-t_k).
\ee\par
On the other hand,~(\ref{unifgammabis}) provides the existence of a sequence $(y^-_k)_{k\in\N}$ of points in $\R^N$ such that
$$y^-_k\in\Omega^-_{s_k},\ \ |z_k-y^-_k|\le r\ \hbox{ and }\ d(y^-_k,\Gamma_{s_k})\ge 2R\ (\ge M)\ \hbox{ for all }k\in\N.$$
In particular,
\be\label{alphamax}
u(s_k,y^-_k)\le\alpha\hbox{ for all }k\in\N,
\ee
due to~(\ref{defMter}).\par
Let us now check that one can choose $t=s_k$ and $x=y^-_k$ in~(\ref{betamin}) for $k$ large enough. Indeed, $s_k\ge t_k+T_{\epsilon}$ for $k$ large enough since $s_k-t_k\to+\infty$ as $k\to+\infty$. Furthermore,
$$|y^-_k-y^+_k|\le|y^-_k-z_k|+|z_k-x_k|+|x_k-y^+_k|\le r+(c_f-2\epsilon)\,(s_k-t_k)+r$$
for all $k\in\N$, whence
$$|y^-_k-y^+_k|\le(c_f-\epsilon)\,(s_k-t_k)\hbox{ for }k\hbox{ large enough},$$
since $s_k-t_k\to+\infty$ as~$k\to+\infty$ and $\epsilon>0$. Finally, one can apply~(\ref{betamin}) with $t=s_k$ and $x=y^-_k$ for $k$ large enough. Thus, $u(s_k,y^-_k)\ge\beta$ for $k$ large enough. This is in contradiction with~(\ref{alphamax}), since $\alpha<\beta$. Finally, our assumption~(\ref{absurde}) cannot hold. In other words,
$$\liminf_{|t-s|\to+\infty}\frac{d(\Gamma_t,\Gamma_s)}{|t-s|}\ge c_f-2\,\epsilon.$$
Since $\epsilon>0$ can be arbitrarily small, the claim~(\ref{claiminf}) follows.\hfill\break\par
{\it Step 2: the upper estimate}. Let us here show that
\be\label{claimsup}
\limsup_{|t-s|\to+\infty}\frac{d(\Gamma_t,\Gamma_s)}{|t-s|}\le c_f.
\ee
Let first $\alpha$ and $\beta$ be fixed as in~(\ref{alphabeta}), let $M\ge0$ be as in~(\ref{defMter}), and let $r\ge 0$ be such that
\be\label{unifgammater}
\forall\,t\in\R,\ \forall\,x\in\Gamma_t,\ \exists\,y^{\pm}_t\in\Omega^{\pm}_t,\ |x-y^{\pm}|\le r\hbox{ and }d(y^{\pm},\Gamma_t)\ge M.
\ee\par
Let $\epsilon>0$ be an arbitrary positive real number. Assume by contradiction that
\be\label{assume}
\limsup_{|t-s|\to+\infty}\frac{d(\Gamma_t,\Gamma_s)}{|t-s|}>c_f+3\,\epsilon.
\ee
There are then two sequences $(t_k)_{k\in\N}$ and $(s_k)_{k\in\N}$ of real numbers such that $|t_k-s_k|\to+\infty$ as $k\to+\infty$ and
$$d(\Gamma_{t_k},\Gamma_{s_k})>(c_f+3\,\epsilon)\,|t_k-s_k|\ \hbox{ for all }k\in\N.$$
Without loss of generality, one can assume that $t_k<s_k$ for all $k\in\N$. For each $k\in\N$, pick a point $z_k$ on $\Gamma_{s_k}$. From~(\ref{unifgammater}), there are two sequences $(y^{\pm}_k)_{k\in\N}$ of points in $\R^N$ such that
$$y^{\pm}_k\in\Omega^{\pm}_{s_k},\ \ |z_k-y^{\pm}_k|\le r\ \hbox{ and }\ d(y^{\pm}_k,\Gamma_{s_k})\ge M\ \hbox{ for all }k\in\N.$$
Thus,~(\ref{defMter}) implies that
\be\label{alphabetabis}
0<u(s_k,y^-_k)\le\alpha<\beta\le u(s_k,y^+_k)<1\hbox{ for all }k\in\N.
\ee
On the other hand, since $d(z_k,\Gamma_{t_k})>(c_f+3\,\epsilon)\,(s_k-t_k)>0$, there holds
$$\hbox{either }\ B\big(z_k,(c_f+3\,\epsilon)\,(s_k-t_k)\big)\subset\Omega^+_{t_k}\ \hbox{ or }\ B\big(z_k,(c_f+3\,\epsilon)\,(s_k-t_k)\big)\subset\Omega^-_{t_k}.$$\par
Let us assume by contradiction that, up to extraction of a subsequence,
\be\label{subset}
B\big(z_k,(c_f+3\,\epsilon)\,(s_k-t_k)\big)\subset\Omega^+_{t_k}\hbox{ for all }k\in\N.
\ee
Two cases shall be considered: $c_f>0$ and $c_f=0$. Consider first the former. Let $R>0$ be given as in Lemma~\ref{leminf}. Since $s_k-t_k\to+\infty$ as $k\to+\infty$, there holds, for $k$ large enough,~${B(z_k,R)}\subset\Omega^+_{t_k}$ together with $d(y,\Gamma_{t_k})\ge M$ for all $y\in{B(z_k,R)}$. Thus, it follows from~(\ref{defMter}) that, for $k$ large enough,
$$u(t_k,y)\ge\beta\hbox{ for all }y\in{B(z_k,R)},$$
whence $u(t_k,x)\ge v_R(0,x-z_k)$ for all $x\in\R^N$ and
\be\label{uvR}
u(t,x)\ge v_R(t-t_k,x-z_k)\hbox{ for all }t>t_k\hbox{ and }x\in\R^N
\ee
from the maximum principle. Let $T_{\epsilon'}>0$ be given by Lemma~\ref{leminf} with $\epsilon'=c_f/2\in(0,c_f]$. It follows from~(\ref{uvR}) and Lemma~\ref{leminf} that, for $k$ large enough,
$$u(t,x)\ge\beta\hbox{ for all }t\ge t_k+T_{\epsilon'}\hbox{ and }|x-z_k|\le(c_f-\epsilon')\,(t-t_k)=\frac{c_f}{2}\,(t-t_k).$$
Since $c_f>0$ and $s_k-t_k\to+\infty$ as $k\to+\infty$, there holds $s_k\ge t_k+T_{\epsilon'}$ and~$|y^-_k-z_k|\le r\le(c_f/2)\,(s_k-t_k)$ for $k$ large enough. Therefore, $u(s_k,y^-_k)\ge\beta$ for $k$ large enough, which is in contradiction with~(\ref{alphabetabis}). As a consequence, the assumption~(\ref{subset}) is ruled out if $c_f>0$.\par
Consider now the case $c_f=0$. Since $s_k-t_k\to+\infty$ as $k\to+\infty$, assumption~(\ref{subset}) implies that, for $k$ large enough, ${B(z_k,2\,\epsilon\,(s_k-t_k))}\subset\Omega^+_{t_k}$ together with $d(y,\Gamma_{t_k})\ge M$ for all $y\in{B(z_k,2\,\epsilon\,(s_k-t_k))}$. It follows from~(\ref{defMter}) that, for $k$ large enough, $u(t_k,y)\ge\beta$ for all~$y\in{B(z_k,2\,\epsilon\,(s_k-t_k))}$. Hence, for $k$ large enough, $u(t_k,x)\ge v_{2\epsilon(s_k-t_k)}(0,x-z_k)$ for all~$x\in\R^N$, whence
\be\label{uvRk}
u(t,x)\ge v_{2\epsilon(s_k-t_k)}(t-t_k,x-z_k)\hbox{ for all }t>t_k\hbox{ and }x\in\R^N
\ee
from the maximum principle. On the other hand, Corollary~\ref{corsup} provides the existence of~$\widetilde{T}_{\epsilon}$ and~$\widetilde{R}_{\epsilon}\ge\epsilon\,\widetilde{T}_{\epsilon}>0$ such that~(\ref{vRbeta}) holds with $c_f=0$ for all $R\ge\widetilde{R}_{\epsilon}$. In particular, since~$s_k-t_k\to+\infty$ as $k\to+\infty$, there holds, for $k$ large enough, $2\,\epsilon\,(s_k-t_k)\ge\widetilde{R}_{\epsilon}$ and
$$\left\{\baa{l}
\widetilde{T}_{\epsilon}\le s_k-t_k\le2\,(s_k-t_k)=\displaystyle{\frac{2\,\epsilon\,(s_k-t_k)}{\epsilon}},\vspace{3pt}\\
|y^-_k-z_k|\le r\le\epsilon\,(s_k-t_k)=2\,\epsilon\,(s_k-t_k)-\epsilon\,(s_k-t_k).\eaa\right.$$
Therefore, the conclusion~(\ref{vRbeta}) can be applied with $c_f=0$, $R=2\,\epsilon\,(s_k-t_k)$, $t=s_k-t_k$ and~$x=y^-_k-z_k$, for $k$ large enough. Finally, it follows from~(\ref{vRbeta}) and~(\ref{uvRk}) that
$$u(s_k,y^-_k)\ge v_{2\epsilon(s_k-t_k)}(s_k-t_k,y^-_k-z_k)\ge\beta$$
for $k$ large enough, which is in contradiction with~(\ref{alphabetabis}).\par
As a conclusion, the assumption~(\ref{subset}) is impossible (even for a subsequence). Hence,
$$B\big(z_k,(c_f+3\,\epsilon)\,(s_k-t_k)\big)\subset\Omega^-_{t_k}$$
for $k$ large enough. Since $s_k-t_k\to+\infty$ as $k\to+\infty$, it follows that, for $k$ large enough, $B\big(z_k,(c_f+2\epsilon)(s_k-t_k)\big)\subset\Omega^-_{t_k}$ and $d(y,\Gamma_{t_k})\ge M$ for all $y\in B\big(z_k,(c_f+2\epsilon)(s_k-t_k)\big)$. Therefore, $u(t_k,y)\le\alpha$ for all $y\in B\big(z_k,(c_f+2\epsilon)(s_k-t_k)\big)$, for $k$ large enough, due to~(\ref{defMter}). Hence, for $k$ large enough, $u(t_k,x)\le w_{(c_f+2\epsilon)(s_k-t_k)}(0,x-z_k)$ for all $x\in\R^N$, and
$$u(t,x)\le w_{(c_f+2\epsilon)(s_k-t_k)}(t-t_k,x-z_k)\hbox{ for all }t>t_k\hbox{ and }x\in\R^N$$
from the maximum principle. Let now $T_{\epsilon}>0$ and $R_{\epsilon}\ge(c_f+\epsilon)T_{\epsilon}>0$ be given by Lemma~\ref{lemsup}, so that~(\ref{wRbeta}) is valid for all $R\ge R_{\epsilon}$. In particular, since $s_k-t_k\to+\infty$ as $k\to+\infty$, there holds, for $k$ large enough, $(c_f+2\epsilon)(s_k-t_k)\ge R_{\epsilon}$ and
$$\left\{\baa{l}
T_{\epsilon}\le s_k-t_k\le\displaystyle{\frac{(c_f+2\epsilon)\,(s_k-t_k)}{c_f+\epsilon}},\vspace{3pt}\\
|y^+_k-z_k|\le r\le\epsilon\,(s_k-t_k)=(c_f+2\epsilon)\,(s_k-t_k)-(c_f+\epsilon)\,(s_k-t_k).\eaa\right.$$
Therefore, the conclusion~(\ref{wRbeta}) can be applied with $R=(c_f+2\epsilon)\,(s_k-t_k)$, $t=s_k-t_k$ and~$x=y^+_k-z_k$, for $k$ large enough. Finally, there holds
$$u(s_k,y^+_k)\le w_{(c_f+2\epsilon)(s_k-t_k)}(s_k-t_k,y^+_k-z_k)\le\alpha$$
for $k$ large enough, which is in contradiction with~(\ref{alphabetabis}).\par
To sum up, we have shown that the assumption~(\ref{assume}) is impossible. Since $\epsilon>0$ can be arbitrarily small, the conclusion~(\ref{claimsup}) follows. The proof of Theorem~\ref{thcf} is thereby complete.~\hfill$\Box$

\begin{rem}\label{remhausdorff}{\rm The second part of the proof of Theorem~$\ref{thcf}$ actually shows that
\be\label{ineq1}
\limsup_{s-t\to+\infty}\Big(\sup_{x\in\Gamma_s}\frac{d(x,\Gamma_t)}{s-t}\Big)\le|c_f|
\ee
for any transition front $u$ connecting $0$ and $1$ for equation~$(\ref{eq})$. This follows imme\-diately from Step~2 of the previous proof when $c_f\ge0$. If $c_f<0$, one can change $u(t,x)$ into~$\widetilde{u}(t,x)=1-u(t,x)$, $f(s)$ into $g(s)=-f(1-s)$, $c_f$ into $-c_f$, $\Omega^{\pm}_t$ into $\Omega^{\mp}_t$, but one can keep the interfaces $\Gamma_t$ for the front~$\widetilde{u}$, so~$(\ref{ineq1})$ is still valid. On the other hand, there holds
\be\label{ineq2}
\limsup_{|t-s|\to+\infty}\frac{\widetilde{d}(\Gamma_t,\Gamma_s)}{|t-s|}\le\limsup_{s-t\to+\infty}\Big(\sup_{x\in\Gamma_s}\frac{d(x,\Gamma_t)}{s-t}\Big),
\ee
where $\widetilde{d}(A,B)$ has been defined in~$(\ref{defdtilde})$. Since $d(A,B)\le\widetilde{d}(A,B)$ for any two subsets~$A$ and~$B$ of~$\R^N$, it finally follows from~$(\ref{ineq1}$-$\ref{ineq2})$ and from the conclusion of Theorem~$\ref{thcf}$ that any transition front $u$ connecting $0$ and $1$ for equation~$(\ref{eq})$ satisfies
$$\frac{\widetilde{d}(\Gamma_t,\Gamma_s)}{|t-s|}\to|c_f|\ \hbox{ as }|t-s|\to+\infty,$$
that is it has a global mean speed equal to $|c_f|$ for the distance $\widetilde{d}$.}
\end{rem}


\SE{Existence of non-standard transition fronts}\label{sec5}

This section is devoted to the proof of Theorem~\ref{thex} on the existence of non-standard transition fronts for problem~(\ref{eq}), that is transition fronts which are not invariant in any moving frame. We first consider the two-dimensional case~$\R^2$. Let us explain heuristically the general strategy before going into the details of the proof. For a better understanding, we refer to the figure shown in Section~\ref{sec2} after the statement of Theorem~\ref{thex}.\par
The first key-ingredient, from~\cite{hmr2,nt1}, is the existence, under the condition~(\ref{bistable}) with~$c_f>0$, of two-dimensional traveling fronts of the type $\phi(x_1,x_2-ct)$ whose level sets are asymptotic to two half-lines having an angle $\alpha\in(0,\pi/2)$ with respect to the $x_2$-axis, see the figure in Section~\ref{sec11}. Consider such a front with $c=c_f/\sin\alpha$ and $\pi/4<\alpha<\pi/2$ and rotate it of angle $\pi/2-\alpha$ clockwise. Namely, we consider the front $\phi(R^{-1}(x_1,x_2-ct))$ where~$R$ denotes the rotation of angle $\pi/2-\alpha$ clockwise. As far as the level sets of this new front are concerned, the right asymptotic half-line becomes parallel to the $x_1$-axis, the other one being then very far away from the $x_2$-axis at very negative times.\par
The next step is to take the restriction of this front on the half-plane $H=\{x_1< 0\}$. One will check that this ``left" traveling front is almost a solution of the same equation~(\ref{eq}) in~$H$ with Neumann boundary conditions on~$\partial H$ for very negative times. One will then solve this Neumann boundary value problem and symmetrize the solution with respect to $\partial H$. Finally, the obtained solution is shown to behave as three moving planar fronts at very negative times, and then as a $V$-shaped classical traveling front $\widetilde{\phi}(x_1,x_2-\widetilde{c}t)$ made of two moving planar fronts for very positive times.


\subsection{Proof of Theorem~\ref{thex}}\label{sec61}

In this section, we carry out the proof of Theorem~\ref{thex}. We leave the proof of some auxiliary lemmas in Section~\ref{sec62}. Throughout the proof of Theorem~\ref{thex}, we assume that $f$ satisfies~(\ref{bistable}) with $c_f>0$, in addition to~(\ref{f}). We repeat that the existence (and uniqueness) of $c_f$ is guaranteed by~(\ref{bistable}). We first consider the case $N=2$ and construct two-dimensional transition fronts satisfying the conclusion of Theorem~\ref{thex}. The conclusion in higher dimensions will be then obtained immediately by trivially extending the constructed two-dimensional fronts in the variables $x_3,\ldots,x_N$.\hfill\break\par

{\it{Step 1: an auxiliary $V$-shaped front.}} Fix an angle $\alpha$ such that
$$\frac{\pi}{4}<\alpha<\frac{\pi}{2}.$$
From~\cite{hmr2,hmr3,nt1}, there exists a unique traveling front $\phi(x_1,x_2-ct)$ of~(\ref{eq}) in~$\R^2$ satisfying the following properties: $0<\phi<1$ in $\R^2$, $\phi$ is of class $C^2(\R^2)$, $c=c_f/\sin\alpha$,
\be\label{limphi}\left\{\baa{rcl}
\displaystyle\mathop{\liminf}_{A\to+\infty}\Big(\mathop{\inf}_{x_2\le|x_1|\cot\alpha-A}\phi(x_1,x_2)\Big) = 1,\vspace{3pt}\\
\displaystyle\mathop{\limsup}_{A\to+\infty}\Big(\mathop{\sup}_{x_2\ge|x_1|\cot\alpha+A}\phi(x_1,x_2)\Big) = 0\eaa\right.
\ee
and $\phi$ is asymptotically planar along the directions $(\pm\sin\alpha,\cos\alpha)$ in the sense that there exist some positive constants $\rho_1$ and $\omega_1$ such that
\be\label{rhoomega}
0\le\phi(x_1,x_2)-\max\big(\phi_f(x_1\cos\alpha+x_2\sin\alpha),\phi_f(-x_1\cos\alpha+x_2\sin\alpha)\big)\le\rho_1\,e^{-\omega_1\sqrt{x_1^2+x_2^2}}
\ee
for all $(x_1,x_2)\in\R^2$. Since $\phi_f(s)$ converges exponentially fast to $0$ and $1$ as $s\to\pm\infty$, the function $\phi$ converges then exponentially fast to $0$ and $1$ as $x_2-|x_1|\cot\alpha\to\pm\infty$. From the Schauder interior estimates, it follows then that there exist some positive constants $\rho_2$ and $\omega_2$ such that
\be\label{rho2omega2}
\big|\nabla\phi(x_1,x_2)\big|\le\rho_2\,e^{-\omega_2|x_2-|x_1|\cot\alpha|}\ \hbox{ for all }(x_1,x_2)\in\R^2.
\ee
Similar arguments yield the existence of some positive constants $\rho_3$ and $\omega_3$ such that
$$\big|\nabla\phi(x_1,x_2)-\nabla\big(\phi_f(-x_1\cos\alpha+x_2\sin\alpha)\big)\big|\le\rho_3\,e^{-\omega_3\sqrt{x_1^2+x_2^2}}\ \hbox{ for all }x_1\ge0,\ x_2\in\R,$$
whence
\be\label{rho3omega3}\left\{\baa{rcl}
\!\big|\phi_{x_1}(x_1,x_2)\!+\!\cos\alpha\,\phi_f'(-x_1\cos\alpha\!+\!x_2\sin\alpha)\big|\le\rho_3\,e^{-\omega_3\sqrt{x_1^2+x_2^2}}\vspace{3pt}\\
\!\big|\phi_{x_2}(x_1,x_2)\!-\!\sin\alpha\,\phi_f'(-x_1\cos\alpha\!+\!x_2\sin\alpha)\big|\le\rho_3\,e^{-\omega_3\sqrt{x_1^2+x_2^2}}\eaa\right.\hbox{for all }x_1\ge0,\ x_2\in\R.
\ee
Lastly, it follows from~\cite{hm} that
\be\label{supstrip}
\forall\,A\ge 0,\quad\sup_{-A\le x_2-|x_1|\cot\alpha\le A}\phi_{x_2}(x_1,x_2)<0
\ee
and that $\phi$ is decreasing in any direction $(\cos\varphi,\sin\varphi)$ such that $\pi/2-\alpha<\varphi<\pi/2+\alpha$. In particular, the function $\phi$ is nonincreasing along the directions $(\pm\sin\alpha,\cos\alpha)$.\hfill\break\par

{\it{Step 2: the rotated $V$-shaped front.}} Let us now rotate the function $\phi$ with angle $\alpha-\pi/2$ clockwise. Namely, we define
\be\label{defpsi}
\psi(x_1,x_2)=\phi(x_1\sin\alpha-x_2\cos\alpha,x_1\cos\alpha+x_2\sin\alpha)
\ee
for all $(x_1,x_2)\in\R^2$. The function $\psi$ is decreasing in any direction $(\cos\varphi,\sin\varphi)$ with~$0<\varphi<2\alpha$. In particular, $\psi$ is nonincreasing in the horizontal direction $(1,0)$ and it converges to the planar front $\phi_f(x_2)$ along this direction. In other words, one of the asymptotic branches of the level sets of $\psi$ corresponds to the half-line $\R_+(1,0)$. The other branch is the half $\R_+(\cos(2\alpha),\sin(2\alpha))$ and it belongs to the left half-plane $\{x_1\le 0\}$ since $\alpha$ is chosen such that $\pi/4<\alpha<\pi/2$. Since $\phi(x_1,x_2-ct)$ solves~(\ref{eq}) in $\R^2$, the $C^2(\R\times\R^2)$ function $\underline{v}$ defined in $\R\times\R^2$ by
\be\label{defundervbis}
\underline{v}(t,x_1,x_2)=\psi(x_1-ct\cos\alpha,x_2-ct\sin\alpha)=\phi(x_1\sin\alpha-x_2\cos\alpha,x_1\cos\alpha+x_2\sin\alpha-ct)
\ee
satisfies~(\ref{eq}) in $\R^2$ too. The function $\underline{v}$ is a traveling front which is invariant in the moving frame with speed $c$ in the direction $(\cos\alpha,\sin\alpha)$, in the sense that
$$\underline{v}(t+\tau,x_1+c\tau\cos\alpha,x_2+c\tau\sin\alpha)=\underline{v}(t,x_1,x_2)$$
for all $(t,x_1,x_2)\in\R\times\R^2$ and for all $\tau\in\R$. At any time $t\in\R$, any level set of $\underline{v}$, that is the set
$$\big\{(x_1,x_2)\in\R^2,\ \underline{v}(t,x_1,x_2)=\lambda\big\}$$
for a given value $\lambda\in(0,1)$ and at a given time $t\in\R$, is asymptotic to some finite shifts of the two half-lines
$$(ct\cos\alpha,ct\sin\alpha)+\R_+(\cos(2\alpha),\sin(2\alpha))\ \hbox{ and }\ (ct\cos\alpha,ct\sin\alpha)+\R_+(1,0),$$
and the first one (i.e. the left one) is very far from the $x_2$-axis for very negative times. More precisely,~(\ref{rhoomega}) implies that
$$\baa{l}
0\le\underline{v}(t,x_1,x_2)-\max\big(\phi_f(x_1\sin(2\alpha)-x_2\cos(2\alpha)-ct\sin\alpha),\phi_f(x_2-ct\sin\alpha)\big)\vspace{3pt}\\
\qquad\qquad\qquad\qquad\qquad\qquad\qquad\qquad\qquad\le\rho_1\,e^{-\omega_1\sqrt{(x_1\sin\alpha-x_2\cos\alpha)^2+(x_1\cos\alpha+x_2\sin\alpha-ct)^2}}\eaa$$
for all $(t,x_1,x_2)\in\R\times\R^2$.\hfill\break\par

{\it{Step 3: a Neumann boundary value problem in a half-space and construction of some sub- and supersolutions.}} In the remaining steps, the strategy consists in constructing a solution of~(\ref{eq}) in $\R^2$ which looks like the function $\underline{v}(t,x_1,x_2)$ for very negative times in the half-plane
$$H=\big\{(x_1,x_2)\in\R^2,\ x_1<0\big\}.$$
To do so, we will work in the half-plane $H$ with Neumann boundary conditions on $\partial H$ and we will then extend the constructed solution by orthogonal symmetry with respect to~$\partial H$.\par
Let us then consider the problem
\be\label{eqhalf}\left\{\baa{rcll}
v_t & = & \Delta v+f(v), & (t,x_1,x_2)\in\R\times H,\vspace{3pt}\\
v_{x_1} & = & 0, & (t,x_1,x_2)=(t,0,x_2)\in\R\times\partial H.\eaa\right.
\ee\par
Remember that the function $\psi$ defined in~(\ref{defpsi}) is nonincreasing along the direction $(1,0)$, that is $\psi_{x_1}(x_1,x_2)\le0$ in $\R^2$. Therefore, $\underline{v}_{x_1}(t,x_1,x_2)\le0$ in $\R\times\R^2$. In particular, the function~$\underline{v}$ is a subsolution of~(\ref{eqhalf}).\par
Problem~(\ref{eqhalf}) also admits a supersolution which looks like the function $\underline{v}$ for very negative times, up to some exponentially small terms, as shown in the following lemma.

\begin{lem}\label{lemsuperv}
There exist some constants $\sigma>0$, $\delta>0$ and $T<0$ such that the function $\overline{v}$ defined in $\R\times\overline{H}$ by
\be\label{defsuperv}
\overline{v}(t,x_1,x_2)=\min\big(\underline{v}(t+\sigma\,e^{\delta t},x_1,x_2)+\delta\,e^{\delta(x_1+t)},1\big)
\ee
is a supersolution of~$(\ref{eqhalf})$ for $t\le T$.
\end{lem}

In order not to lengthen the proof of Theorem~\ref{thex}, the proof of Lemma~\ref{lemsuperv} is postponed in Section~\ref{sec62}.\hfill\break\par

{\it{Step 4: construction of a solution $v$ of~$(\ref{eqhalf})$ in $H$.}} Observe first that
\be\label{undervt}
\underline{v}_t(t,x_1,x_2)=-c\,\phi_{x_2}(x_1\sin\alpha-x_2\cos\alpha,x_1\cos\alpha+x_2\sin\alpha-ct)>0
\ee
for all $(t,x_1,x_2)\in\R\times\R^2$, whence $\underline{v}(t,x_1,x_2)<\overline{v}(t,x_1,x_2)$ in $\R\times\overline{H}$ (remember also that~$\underline{v}<1$ in $\R\times\R^2$). For any $n\in\N$ such that $n>|T|$, let $v^n$ be the solution of the Cauchy problem associated to~(\ref{eqhalf}) for times $t>-n$, with initial condition
$$v^n(-n,x_1,x_2)=\underline{v}(-n,x_1,x_2)\ \hbox{ for all }(x_1,x_2)\in H$$
at time $t=-n$. From Step~3 and the above observations, the maximum principle implies that
$$0<\underline{v}(t,x_1,x_2)\le v^n(t,x_1,x_2)\le\overline{v}(t,x_1,x_2)\le 1\ \hbox{ for all }-n<t\le T\hbox{ and }(x_1,x_2)\in\overline{H},$$
and that
\be\label{vvn}
0<\underline{v}(t,x_1,x_2)\le v^n(t,x_1,x_2)\le 1\hbox{ for all }t>-n\hbox{ and }(x_1,x_2)\in\overline{H}.
\ee
In particular, for every $(t,x_1,x_2)\in\R\times\overline{H}$, the sequence $(v^n(t,x_1,x_2))_{n>\max(|T|,|t|)}$ is nondecreasing. Furthermore, it follows from~(\ref{undervt}) and~(\ref{vvn}) that
$$v^n(-n+h,x_1,x_2)\ge\underline{v}(-n+h,x_1,x_2)>\underline{v}(-n,x_1,x_2)=v^n(-n,x_1,x_2)$$
for all $n>|T|$, $h>0$ and $(x_1,x_2)\in\overline{H}$, whence $v^n$ is increasing with respect to time $t$ in $\overline{H}$, from the maximum principle.\par
From monotone convergence and standard parabolic estimates up to the boundary, the functions~$v^n$ converge then as $n\to+\infty$ in $C^{1,2}_{loc}(\R\times\overline{H})$ to a solution $v$ of~(\ref{eqhalf}) such that
$$0<\underline{v}(t,x_1,x_2)\le v(t,x_1,x_2)\le\overline{v}(t,x_1,x_2)\le 1\ \hbox{ for all }t\le T\hbox{ and }(x_1,x_2)\in\overline{H}$$
and $0<\underline{v}\le v\le 1$ in $\R\times\overline{H}$ (the strong maximum principle also yields $0<v<1$ in $\R\times\overline{H}$, since $\overline{v}(t,x_1,x_2)\to0<1$ as $t\to-\infty$ for each fixed $(x_1,x_2)\in\overline{H}$). Moreover, $v_t\ge0$ in~$\R\times\overline{H}$ with even the strict inequality $v_t>0$ in $\R\times\overline{H}$ from the strong maximum principle, since the previous inequalities prevent $v$ from being independent of time.\hfill\break\par

{\it{Step 5: construction of a solution $u$ of~$(\ref{eq})$ in $\R^2$.}} Define~$u$ in $\R\times\R^2$ as
$$u(t,x_1,x_2)=\left\{\baa{ll}
v(t,x_1,x_2) & \hbox{ for all }t\in\R,\ x_1\le 0,\ x_2\in\R,\vspace{3pt}\\
v(t,-x_1,x_2) & \hbox{ for all }t\in\R,\ x_1> 0,\ x_2\in\R.\eaa\right.$$
Since $v$ satisfies~(\ref{eqhalf}) in the half-plane $H$ with Neumann boundary conditions, it follows that~$u$ is a classical time-global solution of~(\ref{eq}) in the whole plane~$\R^2$. Furthermore, $0<u<1$ in~$\R\times\R^2$, together with
\be\label{undervu1}
\underline{v}(t,-|x_1|,x_2)\le u(t,x_1,x_2)\hbox{ for all }(t,x_1,x_2)\in\R\times\R^2
\ee
and
\be\label{undervu2}
\underline{v}(t,-|x_1|,x_2)\le u(t,x_1,x_2)\le\overline{v}(t,-|x_1|,x_2)\hbox{ for all }t\le T\hbox{ and }(x_1,x_2)\in\R^2.
\ee
Therefore, it follows from~(\ref{rhoomega}),~(\ref{defundervbis}),~(\ref{undervu1}) and the equality $c=c_f/\sin\alpha$ that
\be\label{ineqleft}\baa{r}
\max\big(\phi_f(-|x_1|\sin(2\alpha)-x_2\cos(2\alpha)-c_ft),\phi_f(x_2-c_ft)\big)\le u(t,x_1,x_2)\qquad\qquad\vspace{3pt}\\
\hbox{ for all }(t,x_1,x_2)\in\R\times\R^2\eaa
\ee
and from~(\ref{rhoomega}),~(\ref{defundervbis}),~(\ref{defsuperv}) and~(\ref{undervu2}) that
\be\label{ineqright}\baa{rcl}
u(t,x_1,x_2) & \!\!\!\!\!\!\le\!\!\!\!\!\! & \!\!\max\!\big(\phi_f(-|x_1|\sin(2\alpha)\!-\!x_2\cos(2\alpha)\!-\!c_ft\!-\!c_f\sigma e^{\delta t}),\phi_f(x_2\!-\!c_ft\!-\!c_f\sigma e^{\delta t})\big)\vspace{3pt}\\
& \!\!\!\!\!\!\!\! & \!\!+\rho_1\,e^{-\omega_1\sqrt{(|x_1|\sin\alpha+x_2\cos\alpha)^2+(|x_1|\cos\alpha-x_2\sin\alpha+ct+c\sigma e^{\delta t})^2}}+\delta e^{\delta(t-|x_1|)}\vspace{3pt}\\
& & \qquad\qquad\qquad\qquad\qquad\qquad\qquad\qquad\hbox{ for all }t\le T\hbox{ and }(x_1,x_2)\in\R^2.\eaa
\ee\par

{\it{Step 6: the solution $u$ is a transition front connecting $0$ and $1$.}} To show this property, we need to introduce some families $(\Omega^{\pm}_t)_{t\in\R}$ and $(\Gamma_t)_{t\in\R}$, drawn on the joint figure below, satisfying the properties of Definition~\ref{def1}. For $t\le 0$, set
\be\label{defPlrt}\left\{\baa{ll}
P^l_t=(ct\cos\alpha,ct\sin\alpha)=(ct\cos\alpha,c_ft), & L^l_t=P^l_t+\R_+(\cos(2\alpha),\sin(2\alpha)),\vspace{3pt}\\
P^r_t=(-ct\cos\alpha,ct\sin\alpha)=(-ct\cos\alpha,c_ft), & L^r_t=P^r_t+\R_+(-\cos(2\alpha),\sin(2\alpha))\eaa\right.
\ee
and
\be\label{defgammat-}
\Gamma_t=L^l_t\cup[P^l_t,P^r_t]\cup L^r_t\ \hbox{ for all }t\le 0,
\ee
where the superscript $l$ (resp. $r$) stands for left (resp. right). Define
\be\label{defgammat+}
\Gamma_t=\Big\{(x_1,x_2)\in\R^2;\ x_2=|\!\tan(2\alpha)|\,|x_1|+\frac{c_ft}{|\!\cos(2\alpha)|}\Big\}\ \hbox{ for all }t>0.
\ee
Therefore, for every $t\in\R$, $\Gamma_t$ can be written as a graph $\Gamma_t=\{(x_1,x_2)\in\R^2;\ x_2=\varphi_t(x_1)\}$ for some Lipschitz-continuous function $\varphi_t:\R\to\R$. We finally define, for all $t\in\R$,
\be\label{defomegapmt}
\Omega^+_t=\big\{(x_1,x_2)\in\R^2;\ x_2<\varphi_t(x_1)\big\}\hbox{ and }\Omega^-_t=\big\{(x_1,x_2)\in\R^2;\ x_2>\varphi_t(x_1)\big\}.
\ee
It is immediate to see that the families $(\Omega^{\pm}_t)_{t\in\R}$ and $(\Gamma_t)_{t\in\R}$ satisfy the general properties~(\ref{omegapm}),~(\ref{unifgamma}) and~(\ref{omegapmbis}).
\begin{figure}\centering
\subfigure{\includegraphics[scale=0.5]{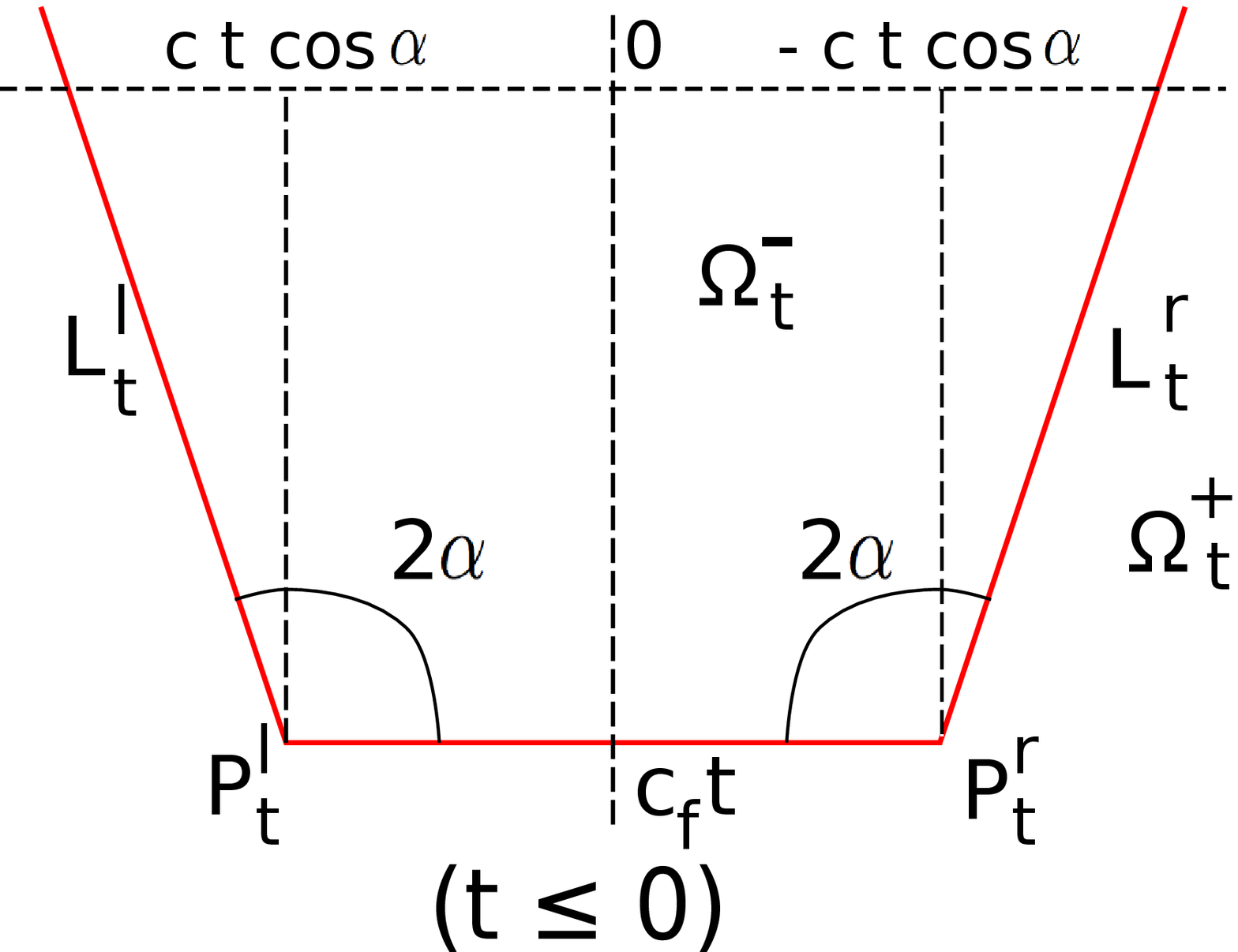}}
\hskip 1cm
\subfigure{\includegraphics[scale=0.4]{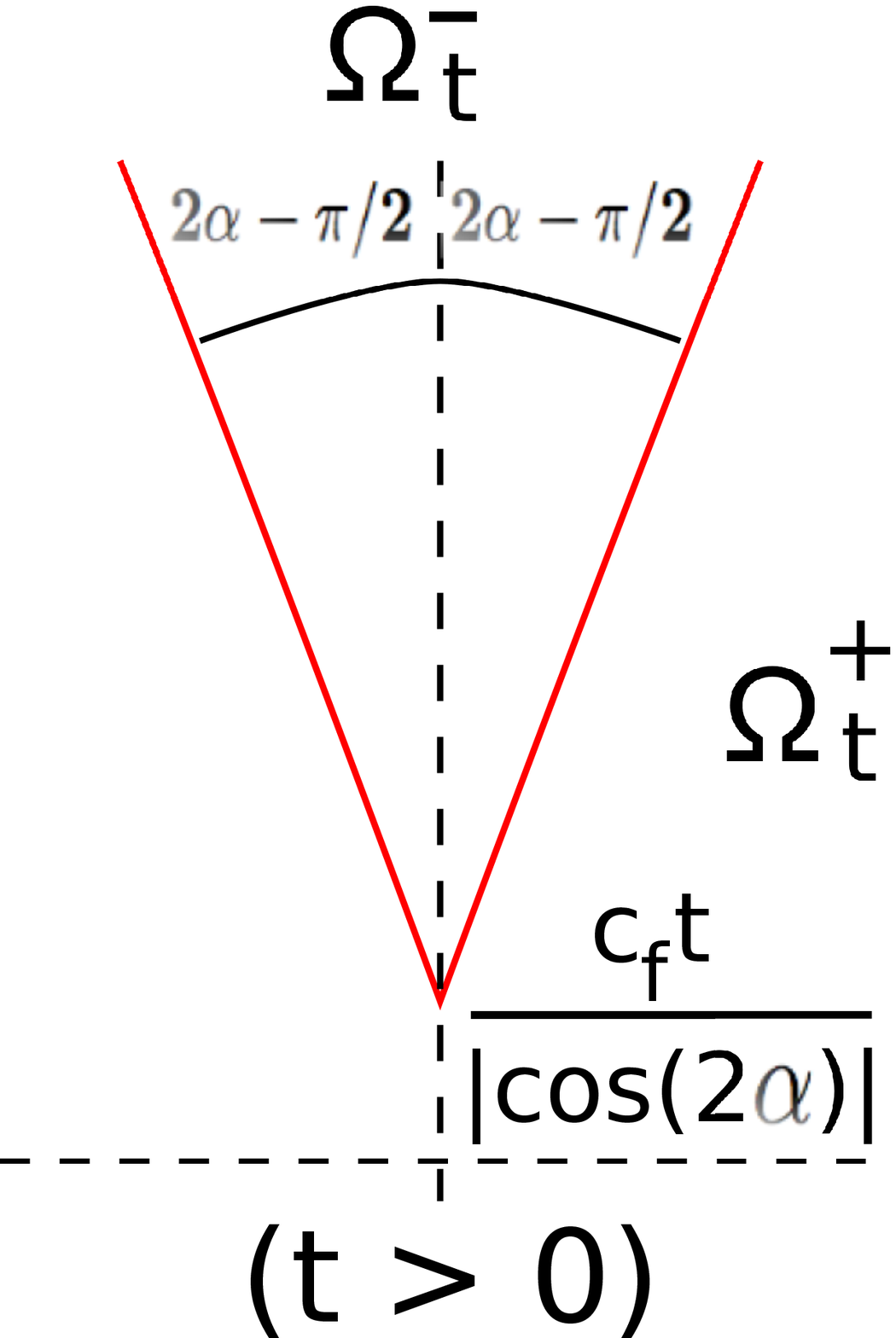}}
\caption{Profiles of the sets $\Gamma_t$ given in~(\ref{defgammat-}) and~(\ref{defgammat+}) for $t\le 0$ and $t>0$}
\end{figure}

\begin{lem}\label{lemtransition}
The function $u$ is a transition front connecting $0$ and $1$ for problem~$(\ref{eq})$ in~$\R^2$ with this choice of sets $(\Omega^{\pm}_t)_{t\in\R}$ and $(\Gamma_t)_{t\in\R}$.
\end{lem}

We point out that, in the course of the proof of Lemma~\ref{lemtransition}, the following interesting pro\-perty is shown: the solution~$u$ converges uniformly in $\R^2$ as $t\to+\infty$ to a traveling front of the type~$\tilde{\phi}(x_1,x_2-\tilde{c}t)$ solving~(\ref{eq}) and~(\ref{deftildephi}) below, with vertical speed $\tilde{c}=c_f\,/\,|\cos(2\alpha)|$.\hfill\break\par

{\it{Step 7: the solution $u$ satisfies the conclusion of Theorem~$\ref{thex}$.}} Assume by contradiction that there exist a function $\Phi:\R^2\to(0,1)$ and some families $(R_t)_{t\in\R}$ and~$(X_t)_{t\in\R}=(x_{1,t},x_{2,t})_{t\in\R}$ of rotations and points in $\R^2$ such that
$$u(t,x_1,x_2)=\Phi(R_t(x_1-x_{1,t},x_2-x_{2,t}))\hbox{ for all }(t,x_1,x_2)\in\R\times\R^2.$$
It would then follow from Lemma~\ref{lemtransition} that there is $M\ge 0$ such that
$$R_t(\Gamma_t-X_t)\subset\big\{(x_1,x_2)\in\R^2;\ d((x_1,x_2),R_s(\Gamma_t-X_s))\le M\big\}\ \hbox{ for all }(t,s)\in\R^2.$$
This is clearly in contradiction with the definitions~(\ref{defgammat-}) and~(\ref{defgammat+}) of the sets~$\Gamma_t$. Therefore, the function $u$ satisfies all properties of Theorem~\ref{thex} in $\R^2$.\hfill\break\par

{\it{Step 8: the case of the space $\R^N$ with $N\ge 3$.}} We start from the solution $u$ of~(\ref{eq}) in $\R^2$ constructed in the previous steps and we extend it trivially in $\R^N$ as
$$\tilde{u}(t,x_1,\cdots,x_N)=u(t,x_1,x_2)\ \hbox{ for all }(t,x_1,\cdots,x_N)\in\R\times\R^N.$$
From the previous steps, the function $\tilde{u}$ is obviously a transition front for problem~(\ref{eq}) in~$\R^N$ with the sets
$$\widetilde{\Omega}^{\pm}_t=\big\{(x_1,\cdots,x_N)\in\R^N;\ (x_1,x_2)\in\Omega^{\pm}_t\big\}$$
for all $t\in\R$. This solution fulfills the desired conclusion and the proof of Theorem~\ref{thex} is thereby complete.\hfill$\Box$


\subsection{Proof of the auxiliary lemmas}\label{sec62}

In this section, we do the proof of Lemmas~\ref{lemsuperv} and~\ref{lemtransition} stated in the previous section. We begin with Lemma~\ref{lemsuperv}.\hfill\break\par

\noindent{\bf{Proof of Lemma~\ref{lemsuperv}.}} As in the paper of Fife and McLeod~\cite{fm} and as in the lemmas of Section~\ref{seckeylemmas}, the general idea is to slightly perturb the function $\underline{v}$ by some exponentially small terms in order to make it a supersolution of~(\ref{eqhalf}). Here, we shall also deal with the boundary conditions on $\partial H$, that is we have to prove that $\overline{v}_{x_1}(t,0,x_2)\ge0$ for all $t\le T$ and $x_2\in\R$. We shall make use of the stability of the limiting states $0$ and $1$, as well as the uniform strict monotonicity~(\ref{supstrip}) of the $V$-shaped front $\phi$ along its level sets. Remember that $f(1)=0$. Thus, in order to prove that the function $\overline{v}$ defined in Lemma~\ref{lemsuperv} is a supersolution of~(\ref{eqhalf}) for $t\le T$, it is sufficient to show that
$$\overline{v}_t\ge\Delta\overline{v}+f(\overline{v})$$
for $(t,x_1,x_2)\in(-\infty,T]\times\overline{H}$ and $\overline{v}_{x_1}\ge0$ on $(-\infty,T]\times\partial H$ in the region where $\overline{v}<1$. Note that in this region, the function $\overline{v}$ is of class~$C^2$.\par
Let us first choose some parameters. Remember that the positive real numbers $\rho_2$, $\omega_2$, $\rho_3$ and $\omega_3$ are given in~(\ref{rho2omega2}) and~(\ref{rho3omega3}). Call
\be\label{defomega4}
\omega_4=\frac{\min(\omega_2c,\omega_3c_f\cos\alpha)}{2}>0.
\ee
We fix a real number $\delta>0$ such that
\be\label{defdelta4}
0<\delta<\min(1,\omega_4)\hbox{ and }f'\le0\hbox{ on }[0,2\delta]\hbox{ and }[1-\delta,1].
\ee
From~(\ref{limphi}), let $C>0$ be such that
\be\label{defC4}\left\{\baa{ll}
\phi(x_1,x_2)\ge1-\delta & \hbox{for all }x_2\le|x_1|\cot\alpha-C,\vspace{3pt}\\
\phi(x_1,x_2)\le\delta & \hbox{for all }x_2\ge|x_1|\cot\alpha+C.\eaa\right.
\ee
From~(\ref{supstrip}), let $\kappa>0$ be such that
\be\label{defkappa4}
\sup_{-C\le x_2-|x_1|\cot\alpha\le C}\phi_{x_2}(x_1,x_2)=-\kappa<0,
\ee
and choose $\sigma>0$ so that
\be\label{sigma4}
\sigma c\kappa\ge L=\max_{[0,1]}|f'|.
\ee
Call
\be\label{defrho4}
\rho_4=(\sin\alpha+\cos\alpha)\,\max\big(\rho_2e^{\omega_2c\sigma},\rho_3\big)>0.
\ee
Let finally $T<0$ be such that
\be\label{defT4}
T\le-2\sigma<0\hbox{ and }\delta^2e^{\delta t}\ge\rho_4e^{\omega_4t}\hbox{ for all }t\le T.
\ee\par
Let us now estimate $\overline{v}_{x_1}$ on the boundary $\partial H$. In this paragraph, we fix a point $(t,0,x_2)$ on $(-\infty,T]\times\partial H$ such that $\overline{v}(t,0,x_2)<1$. From~(\ref{defundervbis}) and~(\ref{defsuperv}), there holds
\be\label{vx1}\baa{rcl}
\overline{v}_{x_1}(t,0,x_2) & = & \underline{v}_{x_1}(t+\sigma e^{\delta t},0,x_2)+\delta^2e^{\delta t}\vspace{3pt}\\
& = & \sin\alpha\,\phi_{x_1}\big(\!-x_2\cos\alpha,x_2\sin\alpha-ct-c\sigma e^{\delta t}\big)\vspace{3pt}\\
& & +\cos\alpha\,\phi_{x_2}\big(\!-x_2\cos\alpha,x_2\sin\alpha-ct-c\sigma e^{\delta t}\big)+\delta^2e^{\delta t}.\eaa
\ee
We shall estimate this quantity when $|x_2-c_ft|\ge(c_f/2)|t|$ and $|x_2-c_ft|\le(c_f/2)|t|$. Consider first the case when $|x_2-c_ft|\ge(c_f/2)|t|$ and $x_2\le 0$. There holds
$$\big|x_2\sin\alpha-ct-c\sigma e^{\delta t}-|x_2|\cos\alpha\cot\alpha\big|=\frac{\big|x_2-c_ft-c_f\sigma e^{\delta t}\big|}{\sin\alpha}\ge-\frac{ct}{2}-c\sigma$$
since $c=c_f/\sin\alpha$ and $t\le T<0$, whence
\be\label{ineqexp}
\overline{v}_{x_1}(t,0,x_2)\ge-(\sin\alpha+\cos\alpha)\,\rho_2\,e^{\omega_2ct/2+\omega_2c\sigma}+\delta^2e^{\delta t}
\ee
from~(\ref{rho2omega2}) and~(\ref{vx1}). If $|x_2-c_ft|\ge(c_f/2)|t|$ and $x_2\ge 0$, then
$$\baa{rcl}
x_2\sin\alpha-ct-c\sigma e^{\delta t}-|x_2|\cos\alpha\cot\alpha & = & \displaystyle\frac{(\sin^2\alpha-\cos^2\alpha)\,x_2}{\sin\alpha}-ct-c\sigma e^{\delta t}\vspace{3pt}\\
& \ge & \displaystyle-ct-c\sigma\ge-\frac{ct}{2}-c\sigma\ge-\frac{cT}{2}-c\sigma\ge0\eaa$$
since $\pi/4\!\le\!\alpha\!<\!\pi/2$ and $t\!\le\!T\!\le\!-2\sigma\!<\!0$, whence~(\ref{ineqexp}) holds. If~$|x_2-c_ft|\!\le\!(c_f/2)|t|\!=\!-(c_f/2)t$, then $x_2\le(c_f/2)t\le0$ and~(\ref{rho3omega3}) and~(\ref{vx1}) yield
\be\label{ineqexp2}\baa{rcl}
\overline{v}_{x_1}(t,0,x_2) & \!\!\!\ge\!\!\! & -\sin\alpha\cos\alpha\,\phi_f'\big(x_2-c_ft-c_f\sigma e^{\delta t}\big)+\cos\alpha\sin\alpha\,\phi_f'\big(x_2-c_ft-c_f\sigma e^{\delta t}\big)\vspace{3pt}\\
& \!\!\!\!\!\! & -(\sin\alpha+\cos\alpha)\,\rho_3\,e^{-\omega_3\sqrt{(x_2\cos\alpha)^2+(x_2\sin\alpha-ct-c\sigma e^{\delta t})^2}}+\delta^2e^{\delta t}\vspace{3pt}\\
& \!\!\!\ge\!\!\! & -(\sin\alpha+\cos\alpha)\,\rho_3\,e^{-\omega_3|x_2|\cos\alpha}+\delta^2e^{\delta t}\vspace{3pt}\\
& \!\!\!\ge\!\!\! & -(\sin\alpha+\cos\alpha)\,\rho_3\,e^{(\omega_3c_ft\cos\alpha)/2}+\delta^2e^{\delta t}.\eaa
\ee
Finally, for all $(t,0,x_2)\in(-\infty,T]\times\partial H$ such that $\overline{v}(t,0,x_2)<1$, there holds
$$\overline{v}_{x_1}(t,0,x_2)\ge-\rho_4e^{\omega_4t}+\delta^2e^{\delta t}\ge0,$$
from~(\ref{defomega4}),~(\ref{defrho4}),~(\ref{defT4}),~(\ref{ineqexp}) and~(\ref{ineqexp2}).\par
As a last step, let us check that $\overline{v}$ is a supersolution of the parabolic equation~(\ref{eqhalf}) inside~$H$. In this paragraph, $(t,x_1,x_2)$ denotes a point in $\R\times\overline{H}$ such that $\overline{v}(t,x_1,x_2)<1$. Since $\underline{v}$ satisfies~(\ref{eq}) in $\R^2$ and since $\delta<1$, one gets from~(\ref{defsuperv}) that
\be\label{overN}\baa{rcll}
\overline{N}(t,x_1,x_2) & \!\!\!:=\!\!\! & \overline{v}_t(t,x_1,x_2)-\Delta\overline{v}(t,x_1,x_2)-f(\overline{v}(t,x_1,x_2))\vspace{3pt}\\
& \!\!\!=\!\!\! & \underline{v}_t(t+\sigma e^{\delta t},x_1,x_2)+\sigma\delta\underline{v}_t(t+\sigma e^{\delta t},x_1,x_2)\,e^{\delta t}+\delta^2e^{\delta(x_1+t)}\vspace{3pt}\\
& \!\!\!\!\!\! & -\Delta\underline{v}(t+\sigma e^{\delta t},x_1,x_2)-\delta^3e^{\delta(x_1+t)}-f(\overline{v}(t,x_1,x_2))\vspace{3pt}\\
& \!\!\!\ge\!\!\! & f\big(\underline{v}(t+\sigma e^{\delta t},x_1,x_2)\big)-f(\overline{v}(t,x_1,x_2))+\sigma\delta\underline{v}_t(t+\sigma e^{\delta t},x_1,x_2)\,e^{\delta t}.\eaa
\ee
Call
$$\zeta_1(x_1,x_2)=x_1\sin\alpha-x_2\cos\alpha\hbox{ and }\zeta_2(t,x_1,x_2)=x_1\cos\alpha+x_2\sin\alpha-ct-c\sigma e^{\delta t},$$
that is $\underline{v}(t+\sigma e^{\delta t},x_1,x_2)=\phi(\zeta_1(x_1,x_2),\zeta_2(t,x_1,x_2))$. If $\zeta_2(t,x_1,x_2)\le|\zeta_1(x_1,x_2)|\cot\alpha-C$, then
$$1>\overline{v}(t,x_1,x_2)>\underline{v}(t+\sigma e^{\delta t},x_1,x_2)=\phi(\zeta_1(x_1,x_2),\zeta_2(t,x_1,x_2))\ge1-\delta$$
from~(\ref{defC4}), whence
\be\label{overNbis}
\overline{N}(t,x_1,x_2)\ge f\big(\underline{v}(t+\sigma e^{\delta t},x_1,x_2)\big)-f(\overline{v}(t,x_1,x_2))+\sigma\delta\underline{v}_t(t+\sigma e^{\delta t},x_1,x_2)\,e^{\delta t}\ge0
\ee
from~(\ref{defdelta4}),~(\ref{overN}) and the positivity of $\underline{v}_t$ in~(\ref{undervt}). If $\zeta_2(t,x_1,x_2)\ge|\zeta_1(x_1,x_2)|\cot\alpha+C$, then~(\ref{defC4}) yields $0<\underline{v}(t+\sigma e^{\delta t},x_1,x_2)\le\delta$, whence
$$0<\underline{v}(t+\sigma e^{\delta t},x_1,x_2)<\overline{v}(t,x_1,x_2)=\underline{v}(t+\sigma e^{\delta t},x_1,x_2)+\delta e^{\delta(x_1+t)}\le 2\delta$$
since $x_1\le 0$ and $t\le T<0$. Thus, as above,~(\ref{overNbis}) holds from~(\ref{defdelta4}),~(\ref{overN}) and the positivity of~$\underline{v}_t$. Lastly, if $-C\le\zeta_2(t,x_1,x_2)-|\zeta_1(x_1,x_2)|\cot\alpha\le C$, then 
$$\baa{l}
f\big(\underline{v}(t+\sigma e^{\delta t},x_1,x_2)\big)-f(\overline{v}(t,x_1,x_2))\vspace{3pt}\\
\qquad\qquad\qquad\qquad\qquad=f\big(\underline{v}(t+\sigma e^{\delta t},x_1,x_2)\big)-f(\underline{v}(t+\sigma e^{\delta t},x_1,x_2)+\delta e^{\delta(x_1+t)})\vspace{3pt}\\
\qquad\qquad\qquad\qquad\qquad\ge-L\,\delta\,e^{\delta(x_1+t)}\eaa$$
since $L=\max_{[0,1]}|f'|$, while
$$\underline{v}_t(t+\sigma e^{\delta t},x_1,x_2)=-c\,\phi_{x_2}(\zeta_1(x_1,x_2),\zeta_2(t,x_1,x_2))\ge c\kappa$$
from~(\ref{defkappa4}). Therefore,
$$\overline{N}(t,x_1,x_2)\ge-L\,\delta\,e^{\delta(x_1+t)}+\sigma\,\delta\,c\,\kappa\,e^{\delta t}\ge\delta\,(\sigma\,c\,\kappa-L)\,e^{\delta t}\ge0$$
from~(\ref{sigma4}),~(\ref{overN}) and the fact that $x_1\le0$.\par
As a conclusion, $\overline{N}(t,x_1,x_2)\ge0$ for all $(t,x_1,x_2)\in(-\infty,T]\times\overline{H}$ such that $\overline{v}(t,x_1,x_2)<1$. The proof of Lemma~\ref{lemsuperv} is thereby complete.\hfill$\Box$\break

\noindent{\bf{Proof of Lemma~\ref{lemtransition}.}} We have to show that $u(t,x_1,x_2)\to1$ (resp. $0$) as $d((x_1,x_2),\Gamma_t)\to+\infty$ with $(x_1,x_2)\in\Omega^+_t$ (resp. $(x_1,x_2)\in\Omega^-_t$), uniformly in $t$.\hfill\break\par
{\it{Step 1: convergence to $1$ in $\Omega^+_t$.}} The inequality~(\ref{ineqleft}) and the fact that $\phi_f$ is decreasing imply that
\be\label{ineqleftbis}\baa{r}
\max\!\big(\phi_f(\!-\!x_1\sin(2\alpha)\!-\!x_2\cos(2\alpha)\!-\!c_ft),\phi_f(x_1\sin(2\alpha)\!-\!x_2\cos(2\alpha)\!-\!c_ft),\phi_f(x_2\!-\!c_ft)\big)\vspace{3pt}\\
\le u(t,x_1,x_2)\eaa
\ee
for all $(t,x_1,x_2)\in\R\times\R^2$. It immediately follows from the definitions~(\ref{defgammat-}),~(\ref{defgammat+}) and~(\ref{defomegapmt}) of $\Gamma_t$ and~$\Omega^+_t$, and from the fact that $\phi_f(-\infty)=1$, that
\be\label{uM1}
\lim_{M\to+\infty}\Big(\inf_{t\in\R,\,(x_1,x_2)\in\Omega^+_t,\ d((x_1,x_2),\Gamma_t)\ge M}u(t,x_1,x_2)\Big)=1.
\ee\par
{\it{Step 2: convergence to $0$ in~$\Omega^-_t$ for negative enough times.}} As far as the behavior of~$u$ in~$\Omega^-_t$ far away from $\Gamma_t$ is concerned, we will consider three cases: when~$t$ is very negative, when~$t$ is very positive and when~$t$ is in some bounded interval. Let us first consider the case when~$t$ is very negative. Let~$\epsilon>0$ be arbitrary. Remember that $T<0$ is given in Lemma~\ref{lemsuperv}. Since~$\phi_f(+\infty)=0$, it follows from the definitions~(\ref{defgammat-}) and~(\ref{defomegapmt}) of~$\Gamma_t$ and~$\Omega^-_t$ for $t\le 0$ that there is $M_1>0$ such that
\be\label{ineqr1}\baa{r}
\forall\,t\le T,\ \forall\,(x_1,x_2)\in\Omega^-_t,\quad\big(d((x_1,x_2),\Gamma_t)\ge M_1\big)\Longrightarrow\qquad\qquad\qquad\qquad\qquad\qquad\vspace{3pt}\\
\big(\max\!\big(\phi_f(-|x_1|\sin(2\alpha)\!-\!x_2\cos(2\alpha)\!-\!c_ft\!-\!c_f\sigma e^{\delta t}),\phi_f(x_2\!-\!c_ft\!-\!c_f\sigma e^{\delta t})\big)\le\displaystyle\frac{\epsilon}{3}\big).\eaa
\ee
Obviously, there is $T_1\le T$ such that
\be\label{ineqr2}
\delta\,e^{\delta(t-|x_1|)}\le\frac{\epsilon}{3}\ \hbox{ for all }t\le T_1\hbox{ and }(x_1,x_2)\in\R^2.
\ee
We now claim that there is $M_2>0$ such that
\be\label{ineqr3}\baa{r}
\forall\,t\le T,\ \forall\,(x_1,x_2)\in\Omega^-_t,\quad\big(d((x_1,x_2),\Gamma_t)\ge M_2\big)\Longrightarrow\qquad\qquad\qquad\qquad\vspace{3pt}\\
\big(\rho_1\,e^{-\omega_1\sqrt{(|x_1|\sin\alpha+x_2\cos\alpha)^2+(|x_1|\cos\alpha-x_2\sin\alpha+ct+c\sigma e^{\delta t})^2}}\le\displaystyle\frac{\epsilon}{3}\big).\eaa
\ee
Otherwise, there would exist a sequence $(t_n,x_{1,n},x_{2,n})_{n\in\N}$ in $\R\times\R^2$ such that
\be\label{contra}
t_n\le T,\ (x_{1,n},x_{2,n})\in\Omega^-_{t_n},\ d((x_{1,n},x_{2,n}),\Gamma_{t_n})\ge n\ \hbox{ for all }n\in\N
\ee
and the sequences
$$(y_n)_{n\in\N}\!:=\!(|x_{1,n}|\sin\alpha\!+\!x_{2,n}\cos\alpha)_{n\in\N}\hbox{ and }(z_n)_{n\in\N}\!:=\!(|x_{1,n}|\cos\alpha\!-\!x_{2,n}\sin\alpha\!+\!ct_n\!+\!c\sigma e^{\delta t_n})_{n\in\N}$$
are bounded. Therefore,
$$\left\{\baa{rcl}
|x_{1,n}| & \!\!\!=\!\!\! & z_n\cos\alpha+y_n\sin\alpha-ct_n\cos\alpha-(c\sigma\cos\alpha)e^{\delta t_n}=-ct_n\cos\alpha+O(1)\vspace{3pt}\\
x_{2,n} & \!\!\!=\!\!\! & \displaystyle-|x_{1,n}|\tan\alpha+\frac{y_n}{\cos\alpha}=ct_n\sin\alpha+O(1)\eaa\right.\hbox{as }n\to+\infty.$$
In other words, owing to the definitions~(\ref{defPlrt}) of the points $P^l_t$ and $P^r_t$ for $t\le 0$, this means that the sequence $(d((x_{1,n},x_{2,n}),\{P^l_{t_n},P^r_{t_n}\}))_{n\in\N}$ is bounded. But the points $P^l_{t_n}$ and $P^r_{t_n}$ lie on~$\Gamma_{t_n}$ for all $n\in\N$ (since $t_n\le T\le 0$). As a consequence, the sequence $(d((x_{1,n},x_{2,n}),\Gamma_{t_n}))_{n\in\N}$ is bounded, contradicting~(\ref{contra}). Finally,~(\ref{ineqr3}) holds for some $M_2>0$ and it follows from~(\ref{ineqright}),~(\ref{ineqr1}),~(\ref{ineqr2}) and~(\ref{ineqr3}) that
\be\label{T1M1}
\forall\,t\le T_1,\ \forall\,(x_1,x_2)\in\Omega^-_t,\quad\big(d((x_1,x_2),\Gamma_t)\ge\max(M_1,M_2)\big)\Longrightarrow\big(u(t,x_1,x_2)\le\epsilon\big).
\ee\par
{\it{Step 3: convergence to $0$ in $\Omega^-_t$ for bounded time intervals.}} We show in this step that
\be\label{tau0}
\forall\,\tau>0,\ \ \lim_{A\to+\infty}\Big(\sup_{|t|\le\tau,\,x_2\ge|\tan(2\alpha)|\,|x_1|+A}u(t,x_1,x_2)\Big)=0.
\ee
Indeed, if this were not true, there would exist a sequence $(t_n,x_{1,n},x_{2,n})_{n\in\N}$ in $\R\times\R^2$ such that
$$(t_n)_{n\in\N}\hbox{ is bounded},\ \lim_{n\to+\infty}\big(x_{2,n}-|\tan(2\alpha)|\,|x_{1,n}|\big)=+\infty\hbox{ and }\liminf_{n\to+\infty}u(t_n,x_{1,n},x_{2,n})>0.$$
Up to extraction of a subsequence, one can assume that $(t_n)_{n\in\N}$ converges to $t_{\infty}\in\R$. From standard parabolic estimates, the functions $u_n$ defined by
$$u_n(t,x_1,x_2)=u(t,x_1+x_{1,n},x_2+x_{2,n})\ \hbox{ for all }(t,x_1,x_2)\in\R\times\R^2$$
converge, up to extraction of a subsequence, to a solution $0\le u_{\infty}\le 1$ of~(\ref{eq}) in $\R^2$ such that~$u_{\infty}(t_{\infty},0,0)>0$. On the other hand,~(\ref{ineqright}) implies that
$$\baa{l}
u_n(t,x_1,x_2)\vspace{3pt}\\
\ \ \le\max\!\big(\phi_f(-|x_1\!+\!x_{1,n}|\sin(2\alpha)\!-\!(x_2\!+\!x_{2,n})\cos(2\alpha)\!-\!c_ft\!-\!c_f\sigma e^{\delta t}),\phi_f(x_2\!+\!x_{2,n}\!-\!c_ft\!-\!c_f\sigma e^{\delta t})\big)\vspace{3pt}\\
\quad\ +\rho_1\,e^{-\omega_1\sqrt{(|x_1+x_{1,n}|\sin\alpha+(x_2+x_{2,n})\cos\alpha)^2+(|x_1+x_{1,n}|\cos\alpha-(x_2\!+\!x_{2,n})\sin\alpha+ct+c\sigma e^{\delta t})^2}}+\delta\,e^{\delta(t-|x_1+x_{1,n}|)}\eaa$$
for all $n\in\N$, $t\le T$ and $(x_1,x_2)\in\R^2$. Since $\phi_f(+\infty)=0$, $\pi/4<\alpha<\pi/2$, and $x_{2,n}-|\tan(2\alpha)|\,|x_{1,n}|\to+\infty$ as $n\to+\infty$, the first term of the right-hand side converges to $0$ as $n\to+\infty$ locally uniformly in $(t,x_1,x_2)\in(-\infty,T]\times\R^2$. The second-term also converges to $0$ as in the proof of~(\ref{ineqr3}). Therefore, by passing to the limit as $n\to+\infty$, one infers that
$$u_{\infty}(t,x_1,x_2)\le\delta\,e^{\delta t}\hbox{ for all }t\le T\hbox{ and }(x_1,x_2)\in\R^2.$$
Let $\eta_0>0$ be such that $f\le 0$ on $[0,\eta_0]$. For any $\eta\in(0,\eta_0]$, there is $t_0\le T$ such that $0\le u_{\infty}(s,x_1,x_2)\le\delta\,e^{\delta s}\le\eta$ for all $s\le t_0$ and $(x_1,x_2)\in\R^2$, whence $0\le u_{\infty}(t,x_1,x_2)\le\eta$ for all $(t,x_1,x_2)\in[s,+\infty)\times\R^2$ from the maximum principle, and then for all $(t,x_1,x_2)\in\R\times\R^2$ since $s$ can be arbitrarily negative. Since $\eta>0$ can be arbitrarily small, it follows that $u_{\infty}\equiv 0$ in $\R\times\R^2$, which contradicts $u_{\infty}(t_{\infty},0,0)>0$. As a consequence,~(\ref{tau0}) is proved.\hfill\break\par
{\it{Step 4: convergence to $0$ in $\Omega^-_t$ for large enough times.}} In the beginning of the proof of Theorem~\ref{thex}, we introduced a $V$-shaped front $\phi(x_1,x_2-ct)$ solving~(\ref{eq}) in~$\R^2$ with $c=c_f/\sin\alpha$. Similarly, since $2\alpha-\pi/2\in(0,\pi/2)$, it follows from~\cite{hmr2,nt1} that there is a unique $V$-shaped front $\tilde{\phi}(x_1,x_2-\tilde{c}t)$ solving~(\ref{eq}) in~$\R^2$ with vertical speed
$$\tilde{c}=\frac{c_f}{\sin(2\alpha-\pi/2)}=\frac{c_f}{|\cos(2\alpha)|},$$
such that the function $\tilde{\phi}$ is of class $C^2(\R^2)$,  $0<\tilde{\phi}<1$ in $\R^2$,
\be\label{deftildephi}\left\{\baa{rcl}
\displaystyle\mathop{\liminf}_{A\to+\infty}\Big(\mathop{\inf}_{x_2\le|x_1|\,|\tan(2\alpha)|-A}\tilde{\phi}(x_1,x_2)\Big) = 1,\vspace{3pt}\\
\displaystyle\mathop{\limsup}_{A\to+\infty}\Big(\mathop{\sup}_{x_2\ge|x_1|\,|\tan(2\alpha)|+A}\tilde{\phi}(x_1,x_2)\Big) = 0,\eaa\right.
\ee
and
$$\tilde{\phi}(x_1,x_2)-\max\big(\phi_f(-x_1\sin(2\alpha)-x_2\cos(2\alpha)),\phi_f(x_1\sin(2\alpha)-x_2\cos(2\alpha))\big)\to0$$
as $x_1^2+x_2^2\to+\infty$ with
$$\tilde{\phi}(x_1,x_2)-\max\big(\phi_f(-x_1\sin(2\alpha)-x_2\cos(2\alpha)),\phi_f(x_1\sin(2\alpha)-x_2\cos(2\alpha))\big)\ge0$$
for all $(x_1,x_2)\in\R^2$. The goal of this step is to show that
\be\label{tildephi}
u(t,x_1,x_2)-\tilde{\phi}(x_1,x_2-\tilde{c}t)\to0\hbox{ as }t\to+\infty\hbox{ uniformly in }(x_1,x_2)\in\R^2.
\ee\par
Observe first that~(\ref{ineqleftbis}) implies that
$$u(0,x_1,x_2)\ge\max\big(\phi_f(-x_1\sin(2\alpha)-x_2\cos(2\alpha)),\phi_f(x_1\sin(2\alpha)-x_2\cos(2\alpha))\big)=:\underline{u}_0(x_1,x_2)$$
for all $(x_1,x_2)\in\R^2$. Let $\underline{u}$ be the solution of the Cauchy problem associated to~(\ref{eq}) in~$\R^2$ with initial condition $\underline{u}_0$ at time $t=0$. It follows from the maximum principle that $u(t,x_1,x_2)\ge\underline{u}(t,x_1,x_2)$ for all $(t,x_1,x_2)\in\R_+\times\R^2$ and from~\cite{hmr1,nt1} that
$$\underline{u}(t,x_1,x_2)-\tilde{\phi}(x_1,x_2-\tilde{c}t)\to0\hbox{ as }t\to+\infty\hbox{ uniformly in }(x_1,x_2)\in\R^2.$$
As a consequence,
\be\label{utildephi}
\liminf_{t\to+\infty}\Big(\inf_{(x_1,x_2)\in\R^2}(u(t,x_1,x_2)-\tilde{\phi}(x_1,x_2-\tilde{c}t))\Big)\ge0.
\ee\par
Let now $t_0\in(-\infty,T]$ be arbitrary. Since $u\le 1$, there holds
\be\label{defu0bar}\baa{r}
u(t_0,x_1,x_2)\le\min\!\big(\phi_f(-|x_1|\sin(2\alpha)\!-\!x_2\cos(2\alpha)\!-\!c_ft_0\!-\!c_f\sigma e^{\delta t_0})\!+\!\varsigma(x_1,x_2),1\big)\vspace{3pt}\\
=:\overline{u}_0(x_1,x_2)\eaa
\ee
for all $(x_1,x_2)\in\R^2$, where
$$\varsigma(x_1,x_2)=\big(u(t_0,x_1,x_2)-\phi_f(-|x_1|\sin(2\alpha)\!-\!x_2\cos(2\alpha)\!-\!c_ft_0\!-\!c_f\sigma e^{\delta t_0})\big)^+$$
and $s^+=\max(s,0)$ denotes the positive part of any real number $s$. Since $\phi_f(-\infty)=1$, $\phi_f(+\infty)=0$ and since $0\le u\le 1$ satisfies~(\ref{ineqleftbis}) and~(\ref{tau0}), one infers that
\be\label{R0}
\lim_{A\to+\infty}\Big(\sup_{|\,x_2-|\tan(2\alpha)|\,|x_1|\,|\ge A}\varsigma(x_1,x_2)\Big)=0.
\ee
On the other hand, it follows from~(\ref{ineqright}), from the fact that $\phi_f$ is decreasing and from the condition $t_0\le T$ that
$$\baa{rcl}
\varsigma(x_1,x_2) & \le & \max\!\big(\phi_f(-|x_1|\sin(2\alpha)\!-\!x_2\cos(2\alpha)\!-\!c_ft_0\!-\!c_f\sigma e^{\delta t_0}),\phi_f(x_2\!-\!c_ft_0\!-\!c_f\sigma e^{\delta t_0})\big)\vspace{3pt}\\
& & -\phi_f(-|x_1|\sin(2\alpha)\!-\!x_2\cos(2\alpha)\!-\!c_ft_0\!-\!c_f\sigma e^{\delta t_0})\vspace{3pt}\\
& & +\rho_1\,e^{-\omega_1\sqrt{(|x_1|\sin\alpha+x_2\cos\alpha)^2+(|x_1|\cos\alpha-x_2\sin\alpha+ct_0+c\sigma e^{\delta t_0})^2}}+\delta e^{\delta(t_0-|x_1|)}\vspace{3pt}\\
& \le & \phi_f(x_2\!-\!c_ft_0\!-\!c_f\sigma e^{\delta t_0})\vspace{3pt}\\
& & +\rho_1\,e^{-\omega_1\sqrt{(|x_1|\sin\alpha+x_2\cos\alpha)^2+(|x_1|\cos\alpha-x_2\sin\alpha+ct_0+c\sigma e^{\delta t_0})^2}}+\delta e^{\delta(t_0-|x_1|)}\eaa$$
for all $(x_1,x_2)\in\R^2$. Therefore, there holds
$$\forall\,A\ge 0,\ \ \lim_{\rho\to+\infty}\Big(\sup_{-A\le x_2-|\tan(2\alpha)|\,|x_1|\le A,\ x_1^2+x_2^2\ge\rho^2}\varsigma(x_1,x_2)\Big)=0,$$
since $\phi_f(+\infty)=0$. Together with~(\ref{R0}), this implies that
\be\label{R0bis}
\lim_{\rho\to+\infty}\Big(\sup_{x_1^2+x_2^2\ge\rho^2}\varsigma(x_1,x_2)\Big)=0.
\ee
Let now $\overline{u}$ be the solution of the Cauchy problem associated to~(\ref{eq}) in~$\R^2$ with initial condition~$\overline{u}_0$ at time $t=t_0$. It follows from the maximum principle that $u(t,x_1,x_2)\le\overline{u}(t,x_1,x_2)$ for all $(t,x_1,x_2)\in[t_0,+\infty)\times\R^2$. It also follows from the nonnegativity of $\varsigma$ and from~(\ref{R0bis}) that
$$\overline{u}(t,x_1,x_2)-\tilde{\phi}\big(x_1,x_2-\tilde{c}t-\frac{c_f\sigma e^{\delta t_0}}{|\cos(2\alpha)|}\big)\to0\hbox{ as }t\to+\infty,\hbox{ uniformly in }(x_1,x_2)\in\R^2,$$
see~\cite{hmr1,nt1}. As a consequence, denoting $\tilde{L}=\sup_{\R^2}|\tilde{\phi}_{x_2}|$, there holds
$$\limsup_{t\to+\infty}\Big(\sup_{(x_1,x_2)\in\R^2}(u(t,x_1,x_2)-\tilde{\phi}(x_1,x_2-\tilde{c}t))\Big)\le\frac{\tilde{L}c_f\sigma e^{\delta t_0}}{|\cos(2\alpha)|}.$$
Since $t_0\le T$ can be arbitrarily negative, one gets that
$$\limsup_{t\to+\infty}\Big(\sup_{(x_1,x_2)\in\R^2}(u(t,x_1,x_2)-\tilde{\phi}(x_1,x_2-\tilde{c}t))\Big)\le0.$$
Together with~(\ref{utildephi}), the desired claim~(\ref{tildephi}) follows.\hfill\break\par
{\it{Step 5: convergence to $0$ in $\Omega^-_t$.}} We here put together the conclusions of the steps~2,~3 and~4. Let $\epsilon>0$ and let $T_1\le T$, $M_1>0$ and $M_2>0$ be as in~(\ref{T1M1}). From~(\ref{deftildephi}),~(\ref{tildephi}) and from the definitions~(\ref{defgammat+}) and~(\ref{defomegapmt}) of~$\Gamma_t$ and~$\Omega^-_t$ for $t\ge 0$, there are $T_2>0$ and $M_3>0$ such that
$$\forall\,t\ge T_2,\ \forall\,(x_1,x_2)\in\Omega^-_t,\quad\big(d((x_1,x_2),\Gamma_t)\ge M_3\big)\Longrightarrow\big(u(t,x_1,x_2)\le\epsilon\big).$$
Lastly, it follows from~(\ref{tau0}) and from the definitions~(\ref{defgammat-}) and~(\ref{defgammat+}) of~$\Gamma_t$ that there exists $M_4>0$ such that
$$\forall\,T_1\le t\le T_2,\ \forall\,(x_1,x_2)\in\Omega^-_t,\quad\big(d((x_1,x_2),\Gamma_t)\ge M_4\big)\Longrightarrow\big(u(t,x_1,x_2)\le\epsilon\big).$$
With~(\ref{T1M1}), one infers that
$$\forall\,t\in\R,\ \forall\,(x_1,x_2)\in\Omega^-_t,\ \ \big(d((x_1,x_2),\Gamma_t)\ge\max(M_1,M_2,M_3,M_4)\big)\Longrightarrow\big(u(t,x_1,x_2)\le\epsilon\big).$$
This, together with~(\ref{uM1}) and the inequality $0\le u\le 1$, means that $u$ is a transition front connecting $0$ and $1$ for problem~(\ref{eq}) in $\R^2$. The proof of Lemma~\ref{lemtransition} is thereby complete.\hfill$\Box$


\end{document}